\documentclass{amsart}
\usepackage{amssymb,verbatim}

\newtheorem{theorem}{Theorem}
\newtheorem{lemma}{Lemma}
\newtheorem{proposition}{Proposition}
\newtheorem{corollary}{Corollary}

\theoremstyle{definition}
\newtheorem{definition}{Definition}
\newtheorem{example}{Example}

\theoremstyle{remark}
\newtheorem{remark}{Remark}

\numberwithin{equation}{section}
\numberwithin{theorem}{section}
\numberwithin{lemma}{section}
\numberwithin{proposition}{section}
\numberwithin{corollary}{section}
\numberwithin{definition}{section}
\numberwithin{example}{section}
\numberwithin{remark}{section}
\newcommand{\thref}[1]{Theorem~{\rm\ref{#1}}}
\newcommand{\prref}[1]{Proposition~{\rm\ref{#1}}}
\newcommand{\leref}[1]{Lemma~{\rm\ref{#1}}}
\newcommand{\coref}[1]{Corollary~{\rm\ref{#1}}}

\newcommand{\deref}[1]{Definition~{\rm\ref{#1}}}
\newcommand{\exref}[1]{Example~{\rm\ref{#1}}}

\newcommand{\reref}[1]{Remark~{\rm\ref{#1}}}

\newcommand{\seref}[1]{Section~{\rm\ref{#1}}}
\newcommand\lbb[1]{\label{#1}}

\def\as{associative}

\def\psalg{pseudo\-algebra}
\def\psalgs{pseudo\-algebras}

\def\tt{\otimes}                               
\def\bt{\boxtimes}
\def\<{\langle}
\def\>{\rangle}
\def\ti{\widetilde}
\def\wti{\widetilde}
\def\what{\widehat}
\def\ov{\overline}

\def\d{\partial}

\def\injto{\hookrightarrow}                    
\def\isoto{\xrightarrow{\sim}}                 
\def\st{\; | \;}                               
\def\symm{S}                                   

\def\mmod{\;\mathrm{mod}\;}

\def\di{{\mathrm{d}}}
\def\diz{\di_0}

\newcommand{\kk}{\mathbf{k}}       
\newcommand{\CC}{\mathbb{C}}       
\newcommand{\ZZ}{\mathbb{Z}}       

\newcommand{\fd}{{\mathfrak d}}

\newcommand{\fg}{\mathfrak{g}}

\def\al{\alpha}                         
\def\be{\beta}
\def\ga{\gamma}

\def\de{\delta}
\def\De{\Delta}
\def\ep{\varepsilon}
\def\io{\iota}
\def\la{\lambda}
\def\La{\Lambda}

\def\om{\omega}
\def\Om{\Omega}
\def\ph{\varphi}
\def\si{\sigma}

\def\th{\theta}

\def\g{{\mathfrak{g}}}      
\def\dd{{\mathfrak{d}}}

\def\gl{{\mathfrak{gl}}}
\def\gld{\gl\,\dd}
\def\sl{{\mathfrak{sl}}}
\def\sld{\sl\,\dd}

\def\Wd{W(\dd)}    
\def\Sd{S(\dd,\chi)}
\def\Hd{H(\dd,\chi,\om)}
\def\Kd{K(\dd,\th)}

\def\A{{\mathcal{A}}}
\def\L{{\mathcal{L}}}
\def\V{{\mathcal{V}}}
\def\W{{\mathcal{W}}}
\def\S{{\mathcal{S}}}
\def\s{{\mathfrak{s}}}

\def\N{\mathcal{N}}

\def\E{\mathcal{E}}

\def\O{\mathcal{O}}

\def\T{\mathcal{T}}           
\def\V{\mathcal{V}}           %
\def\TS{\mathcal{T}_\chi}     
\def\VS{\mathcal{V}_\chi}     %

\def\ue{U}                 

\DeclareMathOperator{\sgn}{sgn}
\DeclareMathOperator{\gr}{gr}

\DeclareMathOperator{\Span}{span}

\DeclareMathOperator{\Ind}{Ind}
\DeclareMathOperator{\ad}{ad}
\DeclareMathOperator{\coad}{coad}
\DeclareMathOperator{\tr}{tr}

\DeclareMathOperator{\sd}{\ltimes}
\DeclareMathOperator{\id}{id}
\DeclareMathOperator{\Id}{Id}
\DeclareMathOperator{\Div}{div}
\DeclareMathOperator{\symp}{S}     
\DeclareMathOperator{\fil}{F}      
\DeclareMathOperator{\coh}{H}      

\DeclareMathOperator{\Der}{Der}

\DeclareMathOperator{\Aut}{Aut}
\DeclareMathOperator{\Hom}{Hom}
\DeclareMathOperator{\End}{End}
\DeclareMathOperator{\Chom}{Chom}
\DeclareMathOperator{\Cend}{Cend}

\DeclareMathOperator{\gc}{gc}

\DeclareMathOperator{\Cur}{Cur}

\DeclareMathOperator{\Rad}{Rad}
\DeclareMathOperator{\Vir}{Vir}
\DeclareMathOperator{\sing}{sing}
\DeclareMathOperator{\coef}{coeff}

\begin{document}

\title[Irreducible Modules over Finite Simple Lie Pseudoalgebras I]
{Irreducible Modules over Finite Simple Lie Pseudoalgebras I. \\
Primitive Pseudoalgebras of Type $W$ and $S$}

\author[B.~Bakalov]{Bojko Bakalov}
\address{Department of Mathematics,
North Carolina State University,
Raleigh, NC 27695, USA}
\email{bojko\_bakalov@ncsu.edu}
\thanks{The first author was partially supported by the Miller Institute
for Basic Research in Science and by an FRPD grant from NCSU}

\author[A.~D'Andrea]{Alessandro D'Andrea}
\address{Dipartimento di Matematica,
Istituto ``Guido Castelnuovo'',
Universit\`a di Roma ``La Sapienza'',
00185 Rome, Italy}
\email{dandrea@mat.uniroma1.it}
\thanks{The second author was supported in part by
PRIN ``Spazi di Moduli e Teoria di Lie'' fundings from MIUR and
project MRTN-CT 2003-505078 ``LieGrits'' of the European Union}

\author[V.~G.~Kac]{Victor G.~Kac}
\address{Department of Mathematics, MIT, Cambridge, MA 02139, USA}
\email{kac@math.mit.edu}
\thanks{The third author was supported in part by NSF grants
DMS-9970007 and DMS-0201017.}

\date{October 7, 2004. {\it Revised}: June 19, 2005}


\maketitle
\tableofcontents


\section{Introduction}\lbb{sintro}
One of the algebraic structures that has emerged recently in the
study of the operator product expansions of chiral fields in
conformal field theory is that of a \emph{Lie conformal algebra}
\cite{K2}. 
Recall that this is a module~$L$ over the algebra of
polynomials $\CC [\partial]$ in the indeterminate $\partial$, endowed
with a $\CC$-linear map
\begin{equation*}
  L \otimes L \to \CC [\lambda] \otimes L \, , \quad
    a \otimes b \mapsto [a_{\lambda} b]\, ,
\end{equation*}
satisfying axioms similar to those of a Lie algebra 
(see \cite{DK,K2}).

Choosing a set of generators $\{ a^i \}_{i \in I}$ of the $\CC
[\partial]$-module $L$, we can write:
\begin{equation*}
  [a^i_{\lambda} a^j] =\sum_k Q^{ij}_k (\lambda ,\partial) a^k\, ,
\end{equation*}
where $Q^{ij}_k$ are some polynomials in $\lambda$ and $\partial$.
The commutators of the corresponding chiral fields $\{ a^i (z)
\}_{i \in I}$ then are:
\begin{equation*}
  [a^i (z) \, , \, a^j (w)] =\sum_k Q^{ij}_k
    (\partial_w,\partial_t) (a^k(t) \delta(z-w))|_{t=w}\, .
\end{equation*}
Letting $P^{ij}_k (x,y)= Q^{ij}_k (-x,x+y)$, we can rewrite this
in a more symmetric form
\begin{equation*}
  [a^i (z) ,a^j(w)] =\sum_k P^{ij}_k (\partial_z ,\partial_w)
     (a^k(w) \delta (z-w))\, .
\end{equation*}
We thus obtain an $H=\CC [\partial]$-bilinear map (i.e., a map of
$H \otimes H$-modules):
\begin{equation}
  \label{eq:1.1}
L \otimes L \to (H \otimes H) \otimes_H L \, , \quad
a \otimes b \mapsto [a * b]\, 
\end{equation}
(where $H$ acts on $H \otimes H$ via the comultiplication map
$\Delta (\partial) = \partial \otimes 1 +1 \otimes \partial$),
defined by
\begin{equation*}
  [a^i * a^j] = \sum_k P^{ij}_k
  (\partial \otimes 1  , 1 \otimes \partial)
  \otimes_H a^k\, .
\end{equation*}
Hence the notion of a $\lambda$-bracket $[a_{\lambda}b]$ is
equivalent to the notion of a $*$-bracket $[a*b]$, as introduced by
Beilinson and Drinfeld \cite{BD}.  For example, the Virasoro conformal
algebra $\Vir =\CC [\partial] \ell$ with $[\ell_{\lambda}\ell ] =
(\partial +2\lambda)\ell$ corresponds to the Virasoro $*$-bracket
\begin{equation}
  \label{eq:1.2}
  [\ell * \ell] = (1 \otimes \partial -\partial \otimes
  1) \otimes_{\CC [\partial]} \ell \, .
\end{equation}

A \emph{Lie pseudoalgebra} is a generalization of the notion of a
Lie conformal algebra for which $\CC [\partial]$ is replaced by
the Hopf algebra $H=U (\fd)$, where
$\fd$ is a finite-dimensional Lie algebra and $U
(\fd)$ is its universal enveloping algebra.
It is defined as an $H$-module $L$ endowed with an $H$-bilinear
map \eqref{eq:1.1} subject to certain skew-symmetry and Jacobi identity
axioms (see \cite{BD,BDK} and \seref{spsarep} below).
%
%
The name ``pseudoalgebra'' is motivated by the fact that this is
an algebra in a pseudotensor category, as introduced in \cite{L,BD}.  
Accordingly, the $*$-bracket is also called a
\emph{pseudobracket}.

In \cite{BDK} we gave a complete classification of \emph{finite}
(i.e.,~finitely generated as an $H$-module) simple Lie
pseudoalgebras.  In order to state the result, we introduce a
generalization of the Virasoro pseudoalgebra \eqref{eq:1.2}
defined for $H=\CC [\partial]$, to the case $H=U
(\fd)$, where $\fd$ is any
finite-dimensional Lie algebra.  This is the Lie pseudoalgebra $W
(\fd) = H \otimes \fd$ with the
pseudobracket
\begin{equation*}
  [(1\otimes a)* (1\otimes b)] = (1 \otimes 1)\otimes_H
  (1\otimes [a,b]) + (b \otimes 1) \otimes_H (1\otimes a)
  -(1\otimes a) \otimes_H (1\otimes b)\, .
\end{equation*}
It is shown in \cite{BDK} that all subalgebras of the Lie
pseudoalgebra $W(\fd)$ are simple and, along
with current Lie pseudoalgebras $\Cur \fg =H \otimes \fg$ with
pseudobracket
\begin{equation*}
  [(1 \otimes a) * (1 \otimes b)] = (1 \otimes 1)
  \otimes_H [a,b]\, ,
\end{equation*}
where $\fg$ is a simple finite-dimensional Lie superalgebra, they
form a complete list of finite simple Lie pseudoalgebras.

The notion of a Lie pseudoalgebra is intimately related to the
more classical notion of a differential Lie algebra.  
Let $L$ be a Lie pseudoalgebra, and let $Y$ be a commutative associative 
algebra with compatible left and right actions of the Hopf algebra $H$.
Then we define a Lie algebra $\A_Y L=Y \otimes_H L$ 
with the obvious left $H$-module structure and the following 
Lie bracket:
\begin{equation*}
[x\tt_H a, y\tt_H b] = \sum_i\, (x f_i)(y g_i) \tt_H c_i \,,
\quad\text{if}\quad
[a*b] = \sum_i\, (f_i\tt g_i) \tt_H c_i \,.
\end{equation*}
%
%

The main tool in the study of Lie pseudoalgebras and their
representations is the \emph{annihilation algebra} $\A_X L$,
where $X=H^*$ is the commutative associative algebra dual to the coalgebra
$H$.  In particular a module over a Lie pseudoalgebra $L$ is the
same as a ``conformal'' module over the extended annihilation Lie
algebra $\fd \ltimes \A_X L$ (see \cite{BDK} and \seref{spsanih} below).

The annihilation algebra of the Lie pseudoalgebra $W
(\fd)$ turns out to be isomorphic to the
linearly compact Lie algebra of all formal vectors fields on a
Lie group whose Lie algebra is $\fd$.  This leads to
a formalism of pseudoforms, similar to the usual formalism of
differential forms, which allows us to define three series of
subalgebras $S (\fd,\chi)$, $H
(\fd,\chi ,\omega)$ and $K
(\fd, \theta)$ of $W(\fd)$.  The
annihilation algebras of the simple Lie pseudoalgebras 
$W(\fd)$, $S (\fd,\chi ,\omega)$, $H(\fd, \chi ,\omega)$
and $K (\fd,\theta)$ are isomorphic to the four series of
Lie--Cartan linearly compact Lie algebras $W_N$, $S_N$, $P_N$
(which is an extension of $H_N$ by a $1$-dimensional center)  and
$K_N$, where $N=\dim \fd$.  
However, the Lie pseudoalgebras 
$S(\fd ,\chi)$, $H (\fd ,\chi ,\omega)$ and $K(\fd ,\theta )$
depend on certain parameters $\chi$, $\omega$ and $\theta$ due to
inequivalent actions of $\fd$ on the annihilation algebra
\cite{BDK}.
It is shown in \cite{BDK} that these series of subalgebras along
with their current generalizations, associated to subalgebras of
$\fd$, exhaust all subalgebras of~$W(\fd)$.

The main goal of the present paper is to give a complete list and
an explicit construction of all irreducible finite modules over
the Lie pseudoalgebras $W(\fd)$ and $S(\fd ,\chi)$.
A module will be called {\em irreducible\/} if it does not contain 
nontrivial proper submodules and in addition the action of the Lie 
\psalg\ on it is not identically zero (see \seref{spsarep} below).
Representation theory of the series $H(\fd, \chi ,\omega)$ and
$K(\fd ,\theta)$ will be treated in sequel papers.

The simplest example of a nontrivial $W (\fd)$-module is the
module $\Omega^0(\fd)=H$ (of rank $1$ over $H$) given by:
\begin{equation*}
  (f \otimes a) * g =- (f \otimes ga) \otimes_H 1\, , 
  \quad f,g \in H \, , \, a \in \fd \, .
\end{equation*}
The corresponding module over the annihilation Lie algebra is
just the representation of the Lie algebra of all formal vector
fields in the space of formal power series.  As in the latter
case, the $W (\fd)$-module $\Omega^0(\fd)$ is the first member of the
pseudo de Rham complex
\begin{equation*}
0 \to \Om^0(\dd) \xrightarrow{\di} \Om^1(\dd)
\xrightarrow{\di} \cdots \xrightarrow{\di} \Om^{N}(\dd) \,,
\end{equation*}
where $\Omega^n (\fd) =H \otimes \Om^n$,
$\Om^n = \bigwedge^n \fd^*$, and $N = \dim \fd$ (see \seref{sdrmw}).

The $W (\fd)$-modules $\Omega^n (\fd)$ of pseudodifferential
forms are special cases of \emph{tensor modules} $\T (U) = H
\otimes U$ over $W (\fd)$, associated to any $\gld$-module
$U$, given by:
\begin{equation}\lbb{eq:1.8}
\begin{split}
(1\tt \d_i)*(1\tt u) &= (1 \tt 1) \tt_H (1 \tt (\ad\d_i)u)
+ \sum_{j=1}^N\, (\d_j \tt 1) \tt_H (1 \tt e_i^j u)
\\
&- (1 \tt \d_i) \tt_H (1 \tt u) \,,
\end{split}
\end{equation}
where $\{ \partial_i \}$ is a basis of $\fd$ and $e^j_i
\partial_k = \de^j_k \partial_i$
(see \seref{stenw}).
Then $\Omega^n (\fd) = \T(\Om^n)$.

Furthermore, for a finite-dimensional $\fd$-module $\Pi$ we
define the \emph{twisting} of $\T (U)$ by $\Pi$ by
$\T(\Pi,U) = H \otimes (\Pi \otimes U)$ and by adding the term
$(1 \otimes 1) \otimes_H (1 \otimes \partial_i u)$ in the
right-hand side of \eqref{eq:1.8}. Then we have the
$\Pi$-twisted pseudo de Rham complex of $W (\fd)$-modules:
\begin{equation*}
0 \to \T(\Pi,\Om^0) \xrightarrow{\di_\Pi} \T(\Pi,\Om^1)
\xrightarrow{\di_\Pi} \cdots \xrightarrow{\di_\Pi} \T(\Pi,\Om^N) 
\end{equation*}
(see \seref{stwdrmw}).

The first main result of the present paper
(\thref{twmod}) states that:

\newcommand{\alphaparenlist}{
  \renewcommand{\theenumi}{\alph{enumi}}%
  \renewcommand{\labelenumi}{(\theenumi)}%
}
 \alphaparenlist 
\begin{enumerate}
\item 
The $W (\fd)$-module $\T (\Pi ,U)$ is irreducible if and only if
$\Pi$ and $U$ are irreducible and $U$ is not isomorphic to one of
the $\gld$-modules $\Om^n=\bigwedge^n \fd^*$ for 
$1 \le n \le N = \dim\dd$;

\item 
The $W (\fd)$-submodule $\di_{\Pi} \T (\Pi ,\Omega^n)$ of 
$\T (\Pi,\Omega^{n+1})$ is irreducible, provided that $\Pi$ is
irreducible, for all $0 \le n \leq N-1$;

\item 
The irreducible $W (\fd)$-modules listed in (a) and (b) exhaust
all irreducible finite $W (\fd)$-modules. (The isomorphic modules
among these are 
$\T (\Pi ,\Omega^0)$ $\simeq$ $\di_{\Pi} \T (\Pi ,\Omega^0)$.)

\end{enumerate}
The corresponding result for $S (\fd ,\chi)$ is \thref{tsmod}.
We also describe the structure of submodules of the $W
(\fd)$- and $S (\fd ,\chi)$-modules $\T (\Pi ,\Omega^n)$ 
(Lemmas \ref{lwirr2} and \ref{lsirr2}).

As in the Lie algebra case, the main part of
the problem is the computation of singular vectors.  However,
in the Lie pseudoalgebra framework the calculations are much
simpler.  In particular, we obtain simpler and more transparent
proofs of the results of Rudakov \cite{Ru1,Ru2}.

In the case when $\dd=\kk\d$ is $1$-dimensional, the Lie \psalg\
$W(\kk\d)$ is isomorphic to the Virasoro \psalg\ $\Vir$ with the
pseudobracket \eqref{eq:1.2}. Now \thref{twmod}
states that every irreducible $W(\kk\d)$-module is of the form
$\T(\Pi,U)$, where $\Pi$ is an irreducible $\kk\d$-module
and $U$ is an irreducible $\gl_1$-module not isomorphic to
$\Om^1 = \dd^*$. The module $\Pi$ is $1$-dimensional over $\kk$ and is
determined by the eigenvalue $\al\in\kk$ of $\d$. Similarly, $U$ 
is $1$-dimensional over $\kk$ and is determined by the eigenvalue 
of $\Id\in\gl_1$, which will be denoted by $\De-1$. If we denote
the corresponding $W(\kk\d)$-module $\T(\Pi,U)$ by $M(\al,\De)$,
then the module $M(\al,\De)$ is irreducible iff $\De\ne 0$,
and these are all nontrivial finite irreducible $W(\kk\d)$-modules.
Thus we recover the classification result of~\cite{CK}.

Note that the category of representations of a Lie \psalg\
is not semi\-simple in general, i.e., complete reducibility of modules 
does not hold. 
To study extensions of modules, as well as central extensions 
and infinitesimal deformations of Lie \psalgs, one defines
\emph{cohomology} of Lie \psalgs\ (see \cite{BKV,BDK}).
The cohomology of the Virasoro conformal algebra $\Vir$ was computed
in \cite{BKV}. The cohomology of $\Wd$ and its subalgebras will
be computed in a future publication.

\section{Basic Definitions}\lbb{sbdef}
In this section, we review some facts and notation from \cite{BDK},
which will be used throughout the paper.
We will work over an algebraically closed field $\kk$ of characteristic $0$.
Unless otherwise specified, all vector spaces, linear maps and
tensor products will be considered over $\kk$.
We will denote by $\ZZ_+$ the set of non-negative integers.

\subsection{Preliminaries on Hopf Algebras}\lbb{shopf}
Let $H$ be a Hopf algebra with a coproduct $\De$,
a counit $\ep$, and an antipode $S$.
We will use the following notation (cf.\ \cite{Sw}):
\begin{align}
\lbb{de1}
\De(h) &= h_{(1)} \tt h_{(2)}, \qquad\quad h\in H,
\\
\lbb{de2}
(\De\tt\id)\De(h) &= (\id\tt\De)\De(h) = h_{(1)} \tt h_{(2)} \tt h_{(3)},
\\
\lbb{de3}
(S\tt\id)\De(h) &= h_{(-1)} \tt h_{(2)},
\quad\text{etc.}
\end{align}
Note that notation \eqref{de2} uses coassociativity of $\De$.
The axioms of antipode and counit can be written
as follows:
\begin{align}
\lbb{antip}
h_{(-1)} h_{(2)} &= h_{(1)} h_{(-2)} = \ep(h),
\\
\lbb{cou}
\ep(h_{(1)}) h_{(2)} &= h_{(1)} \ep(h_{(2)}) = h,
\end{align}
while the fact that $\De$ is a homomorphism of algebras
translates as:
\begin{equation}
\lbb{deprod}
(fg)_{(1)} \tt (fg)_{(2)} = f_{(1)} g_{(1)} \tt f_{(2)} g_{(2)},
\qquad f,g\in H.
\end{equation}
Equations \eqref{antip} and \eqref{cou} imply the following
useful relations:
\begin{equation}
\lbb{cou2}
h_{(-1)} h_{(2)} \tt h_{(3)} = 1\tt h
= h_{(1)} h_{(-2)} \tt h_{(3)}.
\end{equation}

Let $X=H^* := \Hom_\kk(H,\kk)$ be the dual of $H$.
Recall that $H$ acts on $X$
by the formula ($h,f\in H$, $x,y\in X$):
\begin{equation}\lbb{hx}
\langle hx, f\rangle = \langle x, S(h)f\rangle,
\end{equation}
so that
\begin{equation}\lbb{hxy}
h(xy) = (h_{(1)}x) (h_{(2)}y).
\end{equation}
Moreover, $X$ is commutative when $H$ is cocommutative.
Similarly, one can define a right action of $H$ on $X$ by
\begin{equation}\lbb{xh}
\langle xh, f\rangle = \langle x, f S(h)\rangle,
\end{equation}
and then we have
\begin{equation}\lbb{xyh}
(xy)h = (x h_{(1)}) (y h_{(2)}).
\end{equation}
Associativity of $H$ implies that $X$ is an $H$-bimodule,
i.e.,
\begin{equation}\lbb{fxg}
f(xg) = (fx)g, \qquad f,g\in H, \; x\in X.
\end{equation}

Throughout the paper, $H=\ue(\dd)$ will be the universal enveloping algebra
of a finite-dimensional Lie algebra $\dd$.
In this case,
\begin{equation}\lbb{deasa}
\De(a)=a\tt1+1\tt a, \quad S(a)=-a, \quad\qquad a\in\dd \,;
\end{equation}
hence, $\De$ is cocommutative and $S^2=\id$.
Set $N=\dim\dd$ and fix a basis
$\{\d_i\}_{i=1,\dots,N}$ of $\dd$.
Then
\begin{equation}\lbb{dpbw}
\d^{(I)} = \d_1^{i_1} \dotsm \d_N^{i_N} / i_1! \dotsm i_N! \,,
\qquad I = (i_1,\dots,i_N) \in\ZZ_+^N \,,
\end{equation}
is a basis of $H$
(similar to the Poincar\'e--Birkhoff--Witt basis).
Moreover, it is easy to see that
\begin{equation}\lbb{dedn}
\De(\d^{(I)}) = \sum_{ J+K=I }
\d^{(J)} \tt \d^{(K)}.
\end{equation}
For a multi-index $I = (i_1,\dots,i_N)$, let
$|I|=i_1+\cdots+i_N$.
Recall that the canonical increasing filtration of $H=\ue(\dd)$ is given by
\begin{equation}\lbb{filued}
\fil^p \ue(\dd) = \Span_\kk\{ \d^{(I)} \st |I| \le p \} \,,
\qquad p=0,1,2,\dots
\end{equation}
and does not depend on the choice of basis of $\dd$.
This filtration is compatible with the structure of Hopf
algebra (see, e.g., \cite[Section~2.2]{BDK} for more details).
We have:
$\fil^{-1} H = \{0\}$, $\fil^0 H = \kk$, $\fil^1 H = \kk\oplus\dd$.

We define a filtration of $H \tt H$ in the usual way:
\begin{equation}\lbb{filhtth}
\fil^n (H\tt H) = \sum_{i+j=n} \fil^i H\tt\fil^j H \,.
\end{equation}
The following lemma, which is a reformulation of \cite[Lemma 2.3]{BDK},
plays an important role in the paper.
(This lemma holds for any Hopf algebra $H$.)

\begin{lemma}\lbb{lhhh}
{\rm(i)}
The linear maps
\begin{equation*}
H \tt H \to H \tt H \,, \qquad
f \tt g \mapsto (f \tt 1) \De(g)
\end{equation*}
and
\begin{equation*}
H \tt H \to H \tt H \,, \qquad
f \tt g \mapsto (1 \tt f) \De(g)
\end{equation*}
are isomorphisms of vector spaces.
These isomorphisms are compatible with the filtration \eqref{filhtth}.

{\rm(ii)}
For any $H$-module $V$, the linear maps
\begin{equation*}
H \tt V \to (H \tt H) \tt_H V \,, \qquad
h \tt v \mapsto (h \tt 1) \tt_H v
\end{equation*}
and
\begin{equation*}
H \tt V \to (H \tt H) \tt_H V \,, \qquad
h \tt v \mapsto (1 \tt h) \tt_H v
\end{equation*}
are isomorphisms of vector spaces.
In addition, we have{\rm:}
\begin{equation*}
(\fil^n H\tt\kk) \tt_H V = \fil^n (H\tt H) \tt_H V =
(\kk\tt\fil^n H) \tt_H V \,.
\end{equation*}
\end{lemma}

Let us define elements $x_I \in X$ by
$\langle x_I, \d^{(J)} \rangle = \de_I^J$,
where, as usual,
$\de_I^J =1$ if $I=J$ and $\de_I^J =0$ if $I\ne J$.
Then \eqref{dedn} implies $x_J x_K = x_{ J+K }$;
hence,
\begin{equation}\lbb{xi1}
x_I = (x^1)^{i_1} \dotsm (x^N)^{i_N} \,,
\qquad I = (i_1,\dots,i_N) \in\ZZ_+^N \,,
\end{equation}
where
\begin{equation}\lbb{xi2}
x^i = x_{\ep_i} \,,
\qquad \ep_i = (0,\dots,0,\underset{i}1,0,\dots,0) \,,
\;\; i=1,\dots,N \,.
\end{equation}
Therefore, $X$ can be identified with the ring
$\O_N = \kk[[t^1,\dots,t^N]]$
of formal power series in $N$ indeterminates.
We have a ring isomorphism
\begin{equation}\lbb{phxon}
\ph\colon X \isoto \O_N \,,
\quad \ph(x^i)=t^i \,, \; \ph(x_I)=t_I \,,
\end{equation}
where $t_I$ is given by a formula similar to \eqref{xi1}.

Let $\fil_p X = (\fil^p H)^\perp$ be the set of elements
from $X=H^*$ that vanish on $\fil^p H$.
Then $\{\fil_p X\}$ is a decreasing filtration of $X$ such that
$\fil_{-1} X = X$, $X/\fil_0 X \simeq\kk$, $\fil_0 X/\fil_1 X \simeq\dd^*$.
Under the isomorphism \eqref{phxon}, 
the filtration $\{\fil_p X\}$ becomes
\begin{equation}\lbb{filon}
\fil_p\O_N=(t^1,\dots,t^N)^{p+1}\,\O_N \,,
\qquad p=-1,0,1,\dots \;.
\end{equation}
This filtration has the following properties:
\begin{equation}\lbb{dfpx}
(\fil_n X) (\fil_p X) \subset \fil_{n+p+1} X \,,
\quad \dd(\fil_p X) \subset \fil_{p-1} X \,,
\quad (\fil_p X)\dd \subset \fil_{p-1} X \,.
\end{equation}
We can consider $x^i$ as elements of $\dd^*$; then $\{x^i\}$
is a basis of $\dd^*$ dual to the basis $\{\d_i\}$ of $\dd$,
i.e., $\<x^i,\d_j\>=\de^i_j$.

We define a topology of $X$ by considering $\{\fil_p X\}$
as a fundamental system of neighborhoods of $0$.
We will always consider $X$ with this topology, while $H$ and $\dd$ with
the discrete topology. Then $X$ is linearly compact
(see \cite[Chapter~6]{BDK}),
and the multiplication of $X$ and the (left and right) actions of $\dd$ on
it are continuous (see \eqref{dfpx}).

\begin{example}\lbb{edcommact}
When $\dd$ is commutative, its left and right actions on $\O_N$
coincide and are given by $\d_i \mapsto -\d/\d t^i$ for $i=1,\dots,N$.
\end{example}

The following lemma is well known (see also \cite[Section~6]{Re}).

\begin{lemma}\lbb{ldacton}
Let $c_{ij}^k$ be the structure constants of\/ $\dd$ in the basis $\{\d_i\}$,
so that $[\d_i,\d_j]=\sum\,c_{ij}^k\d_k$. Then we have the
following formulas for the left and right actions of\/
$\dd$ on $X${\rm:}
\begin{align*}
\d_i x^j &= -\de_i^j -\sum_{k<i}\, c_{ik}^j x^k \mod\fil_1 X \,,
\\
x^j \d_i &= -\de_i^j + \sum_{k>i}\, c_{ik}^j x^k \mod\fil_1 X \,.
\intertext{In particular,}
\d_i x^j - x^j \d_i &= -\sum_k\, c_{ik}^j x^k \mod\fil_1 X
\end{align*}
is the coadjoint action of\/ $\dd$ on
$\dd^* \simeq \fil_0 X / \fil_1 X$.
\end{lemma}
\begin{proof}
We will prove the first equality. The second one is proved in the same way,
while the third follows from the other two.
If we express $\d_i x^j$ in the basis $\{x_K\}$ of $X$, we have
\begin{equation*}
\d_i x^j = \sum_{K\in\ZZ_+^N} a_K x_K \iff
a_K = \< \d_i x^j , \d^{(K)} \> \,,
\end{equation*}
where $\d^{(K)}$ are from \eqref{dpbw}.
Since we are interested in $\d_i x^j \mmod\fil_1 X$,
we need to compute $a_K$ only for $|K| \le 1$,
i.e., only for $\d^{(K)} = 1$ or $\d^{(K)} = \d_k$.
Using \eqref{hx}, we obtain
\begin{align*}
\< \d_i x^j , 1 \> &= -\< x^j , \d_i \> = -\de_i^j \,,
\\
\< \d_i x^j , \d_k \> &= -\< x^j , \d_i \d_k \>
= -\< x^j , \d_k \d_i \>  -\< x^j , [\d_i, \d_k] \> \,.
\end{align*}
If $i \le k$, then $\d_i \d_k$ is (up to a constant)
an element of the basis \eqref{dpbw}
and $\d_i \d_k \ne \d_j$; hence, $\< x^j , \d_i \d_k \> = 0$.
If $i>k$, then by the same argument $\< x^j , \d_k \d_i \> = 0$,
while $\< x^j , [\d_i, \d_k] \> = c_{ik}^j$.
This completes the proof.
\end{proof}

\subsection{Pseudoalgebras and Their Representations}\lbb{spsarep}
In this subsection, we recall the definition of a \psalg\ from
\cite[Chapter~3]{BDK}.
Let $A$ be a (left) $H$-module. A {\em pseudoproduct\/} on $A$
is an $H$-bilinear map
\begin{equation}\lbb{psprod}
  A \tt A \to (H \tt H) \tt_H A \,, \quad
  a \tt b \mapsto a * b \,,
\end{equation}
where we use the
comultiplication $\Delta\colon H \to H \tt H$ to define
$(H\tt H) \tt_H A$.
A \emph{pseudoalgebra} is a (left) $H$-module $A$ endowed
with a pseudoproduct \eqref{psprod}.  The name is motivated by
the fact that this is an algebra in a pseudotensor category,
as introduced in \cite{L, BD} (see \cite[Chapter~3]{BDK}).

In order to define associativity of a pseudoproduct, we extend
it from $ A \tt A \to H^{\tt 2} \tt_H A$ to
$  (H^{\tt 2} \tt_H A) \tt A \to H^{\tt 3}
\tt_H A$ and to $A \tt (H^{\tt 2}
\tt_H A) \to H^{\tt 3}\tt_H A$ by letting:
\begin{align}
\lbb{psprod1}
(h\tt_{H}a)*b &= \sum \, (h \tt 1)\, (\Delta\tt\id)(g_i)
 \tt_H c_i \,,
\\
\lbb{psprod2}
a*(h\tt_{H}b) &= \sum \, (1 \tt h)\, (\id\tt\Delta)(g_i)
 \tt_H c_i \,,
\intertext{where $h\in H^{\tt 2}$, $a,b\in A$, and}
\lbb{psprod3}
a*b &= \sum \, g_i \tt_H c_i
\qquad\text{with \; $g_i\in H^{\tt2}$, $c_i \in A$.}
\end{align}
Then the associativity
property is given by the usual equality (in $H^{\tt 3} \tt_H A$):
\begin{equation}\lbb{psas}
  (a*b)*c=a*(b*c) \,.
\end{equation}

The main objects of our study are Lie pseudoalgebras.  The
corresponding pseudoproduct is conventionally called
{\em pseudobracket\/} and denoted by $[a*b]$.
A {\em Lie \psalg\/} is a (left) $H$-module equipped with
a pseudobracket satisfying the
following skew-commutativity and Jacobi identity axioms:
\begin{align}
  \lbb{psss}
  [b*a] &= -(\sigma \tt_H \id)\, [a*b] \,, \\
  \lbb{psjac}
  [[a*b]*c] &= [a*[b*c]] - ((\sigma \tt \id)\tt_H\id)\,
  [b*[a*c]] \,.
\end{align}
Here, $\sigma\colon H \tt H \to H \tt H$ is the permutation
of factors, and the compositions $[[a*b]*c]$, $[a*[b*c]]$
are defined using \eqref{psprod1}, \eqref{psprod2}.

\begin{remark}\lbb{assolie*}
Let $A$ be an \as\ \psalg\ with a pseudoproduct $a*b$.
Define a pseudobracket on $A$ as the commutator
\begin{equation}\lbb{aslie}
[a*b] = a*b - (\si\tt_H\id) \, (b*a) \,.
\end{equation}
Then, with this pseudobracket, $A$ is a Lie \psalg.
\end{remark}
\begin{example}\lbb{ecurps}
For any $\kk$-algebra $A$, let its
associated \emph{current} $H$-pseudoalgebra be $\Cur A =H\tt A$
with the pseudoproduct
\begin{equation}\lbb{cura}
  (f \tt a)*(g \tt b)=
  (f \tt g) \tt_H (1\tt ab) \,.
\end{equation}
Then the $H$-pseudoalgebra $\Cur A$ is Lie (or associative)
iff the $\kk$-algebra $A$ is.
\end{example}

%
%

The definitions of modules over Lie (or associative) \psalgs\
are obvious modifications of the above.
A {\em module\/} over a Lie \psalg\ $L$
is a left $H$-module $V$ together with an $H$-bilinear map
\begin{equation}\lbb{psprod4}
  L \tt V \to (H \tt H) \tt_H V \,, \quad
  a \tt v \mapsto a * v
\end{equation}
that satisfies ($a,b\in L$, $v\in V$)
\begin{equation}\lbb{psrep}
[a*b]*v = a*(b*v) - ((\si\tt\id)\tt_H\id) \, (b*(a*v)) \,.
\end{equation}
An $L$-module $V$ will be called {\em finite\/}
if it is finite (i.e., finitely generated) as an $H$-module.
The {\em trivial\/} $L$-module is the set $\{0\}$.

A subspace $W\subset V$ is an {\em $L$-submodule\/} if
it is an $H$-submodule and
$L*W \subset (H\tt H)\tt_H W$.
(Here $L*W$ is the linear span of all elements $a*w$,
where $a\in L$ and $w\in W$.)
A submodule $W\subset V$ is called {\em proper\/} if $W\ne V$.
An $L$-module $V$ is {\em irreducible\/} (or {\em simple})
if it does not contain any nontrivial proper $L$-submodules and
$L*V \ne \{0\}$.

Let $U$ and $V$ be two $L$-modules. A map
$\be\colon U\to V$ is a {\em homomorphism\/} of $L$-modules
if $\be$ is $H$-linear and it satisfies
\begin{equation}\lbb{psprod6}
\bigl( (\id\tt\id)\tt_H \be \bigr) (a*u)
= a * \be(u) \,, \qquad a\in L \,, \; u\in U \,.
\end{equation}

\begin{remark}\lbb{rlmod}
{\rm(i)}
Let $V$ be a module over a Lie \psalg\ $L$ and let
$W$ be an $H$-submodule of $V$. By  \leref{lhhh}(ii),
for each $a\in L$, $v\in V$, we can write
\begin{equation*}
a*v = \sum_{I\in\ZZ_+^N} (\d^{(I)} \tt 1) \tt_H v'_I \,,
\qquad v'_I \in V \,,
\end{equation*}
where the elements $v'_I$ are uniquely determined by $a$ and $v$.
Then $W\subset V$ is an $L$-submodule iff
it has the property that all $v'_I \in W$ whenever $v \in W$.
This follows again from \leref{lhhh}(ii).

{\rm(ii)}
Similarly, for each $a\in L$, $v\in V$, we can write
\begin{equation*}
a*v = \sum_{I\in\ZZ_+^N} (1 \tt \d^{(I)}) \tt_H v''_I \,,
\qquad v''_I \in V \,,
\end{equation*}
and $W$ is an $L$-submodule iff $v''_I \in W$ whenever $v \in W$.
\end{remark}
\begin{example}\lbb{elv0}
Let $L$ be a Lie \psalg, and let $V$ be an $L$-module,
which is finite dimensional $($over $\kk)$. Then the
action of $L$ on $V$ is trivial, i.e., $L*V = \{0\}$.
Indeed, since $\dim H = \infty$, every element $v\in V$ is torsion,
i.e., such that $hv=0$ for some nonzero $h\in H$.
Then the statement follows from \cite[Corollary~10.1]{BDK}.
\end{example}

\subsection{Annihilation Algebras of Lie Pseudoalgebras}\lbb{spsanih}
For a Lie $H$-\psalg\ $L$, we set $\A(L)=X\tt_H L$,
where as before $X=H^*$.
We define a Lie bracket on $\L=\A(L)$ by the formula
(cf.\ \cite[Eq.~(7.2)]{BDK}):
\begin{equation}\lbb{alliebr}
[x\tt_H a, y\tt_H b] = \sum\, (x f_i)(y g_i) \tt_H c_i \,,
\quad\text{if}\quad
[a*b] = \sum\, (f_i\tt g_i) \tt_H c_i \,.
\end{equation}
Then $\L$ is a Lie algebra, called the {\em annihilation algebra\/}
of $L$ (see \cite[Section~7.1]{BDK}).
We define a left action of
$H$ on $\L$ in the obvious way:
\begin{equation}\lbb{hactsonl}
h(x\tt_H a) = hx\tt_H a.
\end{equation}
In the case $H=\ue(\dd)$, the Lie algebra $\dd$ acts on $\L$ by
derivations. The semidirect sum $\ti\L = \dd\ltimes\L$ is called
the {\em extended annihilation algebra}.

Similarly, if $V$ is a module over a Lie \psalg\ $L$, we let $\A(V)=X\tt_H V$,
and define an action of $\L=\A(L)$ on $\A(V)$ by:
\begin{equation}\lbb{alav}
(x\tt_H a)(y\tt_H v) = \sum\, (x f_i)(y g_i) \tt_H v_i \,,
\quad\text{if}\quad
a*v = \sum\, (f_i\tt g_i) \tt_H v_i \,.
\end{equation}
We also define an $H$-action on $\A(V)$ similarly to \eqref{hactsonl}.
Then $\A(V)$ is an $\ti\L$-module \cite[Proposition~7.1]{BDK}.

When $L$ is a finite $H$-module, we can define a filtration on
$\L$ as follows (see \cite[Section~7.4]{BDK} for more details).
We fix a finite-dimensional vector
subspace $L_0$ of $L$ such that $L = HL_0$, and set
\begin{equation}\lbb{fill}
\fil_p\L = \{ x \tt_H a \in\L \st x\in\fil_p X \,, \; a\in L_0 \}
\,, \qquad p\ge -1 \,.
\end{equation}
The subspaces $\fil_p \L$ constitute a decreasing
filtration of $\L$, satisfying
\begin{equation}\lbb{filbr}
[\fil_n \L, \fil_p \L] \subset \fil_{n+p-\ell} \L \,,
\quad \dd (\fil_p \L) \subset \fil_{p-1} \L \,,
\end{equation}
where $\ell$ is an integer depending only on the choice of $L_0$.
Notice that the filtration just defined depends on the choice of
$L_0$, but the topology it induces does not \cite[Lemma~7.2]{BDK}.
We set $\L_p = \fil_{p+\ell} \L$, so that $[\L_n, \L_p] \subset \L_{n+p}$.
In particular, $\L_0$ is a Lie algebra.

%
%

We also define a filtration of $\ti\L$ by letting
$\fil_{-1}\ti\L=\ti\L$, $\fil_p\ti\L = \fil_p\L$ for $p\ge0$,
and we set $\ti\L_p = \fil_{p+\ell} \ti\L$.
An $\ti\L$-module $V$ is called {\em conformal\/}
if every $v\in V$ is killed by some $\L_p$; in other words, if $V$ is
a topological $\ti\L$-module when endowed with the discrete topology.


The next two results from \cite{BDK} will play a crucial role in our
study of representations (see \cite{BDK}, Propositions 9.1 and 14.2,
and Lemma 14.4).

\begin{proposition}\lbb{preplal2}
Any module $V$ over the Lie \psalg\ $L$ has a natural structure of a
conformal $\ti\L$-module, given by the action of\/ $\dd$ on $V$ and by
\begin{equation}\lbb{axm2}
(x\tt_H a) \cdot v
= \sum\, \< x, S(f_i {g_i}_{(-1)}) \> \, {g_i}_{(2)} v_i \,,
\quad\text{if}\quad
a*v = \sum\, (f_i\tt g_i) \tt_H v_i
\end{equation}
for $a\in L$, $x\in X$, $v\in V$.

Conversely, any conformal $\ti\L$-module $V$
has a natural structure of an $L$-module, given by
\begin{equation}\lbb{prpl2}
a*v = \sum_{I\in\ZZ_+^N} \bigl( S(\d^{(I)}) \tt1 \bigr)\tt_H
\bigl( (x_I\tt_H a) \cdot v \bigr) \,.
\end{equation}

Moreover, $V$ is irreducible
as an $L$-module iff it is irreducible as an $\ti\L$-module.
\end{proposition}

\begin{lemma}\lbb{lkey2}
Let $L$ be a finite Lie \psalg\
and $V$ be a finite $L$-module. For $p\ge -1-\ell$, let
\begin{equation*}\lbb{kernv}
\ker_p V
= \{ v \in V \st \L_p \, v = 0 \},
\end{equation*}
so that, for example, $\ker_{-1-\ell} V = \ker V$
and\/ $V = \bigcup \ker_p V$.
Then all vector spaces\/
$\ker_p V / \ker V$ are finite dimensional.
In particular, if\/ $\ker V=\{0\}$, then
every vector $v\in V$ is contained in a finite-dimensional subspace
invariant under~$\L_0$.
\end{lemma}

\section{Primitive Lie Pseudoalgebras of Type $W$ and $S$}
\lbb{sprim}
Here we introduce the main objects of our study: the {\em primitive\/}
Lie pseudoalgebras $\Wd$ and $\Sd$ and their annihilation algebras $\W$ and
$\S$ (see \cite[Chapter~8]{BDK}).

\subsection{Definition of $\Wd$ and $\Sd$}\lbb{subprim}
We define the Lie pseudoalgebra $\Wd$ as the free $H$-module
$H\tt\dd$ with the pseudobracket
\begin{equation}\lbb{wdbr*}
\begin{split}
[(f\tt a)*(g\tt b)]
&= (f\tt g)\tt_H(1\tt [a,b])
\\
&- (f\tt ga)\tt_H(1\tt b) + (fb\tt g)\tt_H(1\tt a) \,.
\end{split}
\end{equation}
There is a structure of a $\Wd$-module on $H$ given by:
\begin{equation}\lbb{wdac*}
(f\tt a)*g = -(f\tt ga)\tt_H 1 \,.
\end{equation}

Let $\chi$ be a trace form on $\dd$, i.e., a linear functional
from $\dd$ to $\kk$ that vanishes on $[\dd,\dd]$.
Define an $H$-linear map $\Div^\chi\colon\Wd\to H$ by the formula:
\begin{equation}\lbb{div*}
\Div^\chi \left( \sum\, h_i\tt \d_i \right)
= \sum\, h_i (\d_i + \chi(\d_i)) \,.
\end{equation}
Then
\begin{equation}\lbb{sd}
\Sd := \{ s \in\Wd \st \Div^\chi s = 0 \}
\end{equation}
is a subalgebra of the Lie \psalg\ $\Wd$.
It was shown in \cite[Proposition~8.1]{BDK} that
$\Sd$ is generated over $H$ by the elements
\begin{equation}\lbb{achib}
s_{ab} := (a+\chi(a))\tt b - (b+\chi(b))\tt a - 1\tt [a,b]
\qquad\text{for}\quad a,b\in\dd \,.
\end{equation}
Pseudobrackets of the elements $s_{ab}$ are explicitly
calculated in \cite[Proposition~8.1]{BDK}.
Notice that when $\dim\dd>2$, $\Sd$ is not free as an $H$-module,
because the elements $s_{ab}$ satisfy the relations
\cite[Eq.~(8.23)]{BDK}.

\begin{remark}\lbb{rdimd2}
If $\dim\dd=1$, then $\Sd=\{0\}$.
If $\dim\dd=2$, the Lie \psalg\ $\Sd$ is free as an $H$-module
and it is isomorphic to a primitive Lie \psalg\ of type $H$
(see \cite{BDK}, Section~8.6 and Example~8.1).
\end{remark}

Irreducible modules over primitive Lie \psalgs\ of type $H$
will be studied in a sequel paper.
{}From now on, whenever we consider the Lie \psalg\ $\Sd$,
we will assume that $\dim\dd>2$.

\subsection{Annihilation Algebra of $\Wd$}\lbb{subanw}
Let $\W = \A(\Wd)$ be the annihilation algebra of the
Lie \psalg\ $\Wd$. Since $\Wd=H \tt \dd$, we have
$\W = X \tt_H (H \tt \dd) \equiv X \tt \dd$, so we can
identify $\W$ with $X\tt\dd$. Then the Lie bracket \eqref{alliebr}
in $\W$ becomes ($x,y\in X$, $a,b\in\dd$):
\begin{equation}\lbb{Wbra}
[x \tt a, y \tt b] = xy \tt [a, b] - x(ya) \tt b + (xb)y \tt a \,,
\end{equation}
while the left action \eqref{hactsonl} of $H$ on $\W$ is given by:
$h(x \tt a) = hx \tt a$. The Lie algebra $\dd$ acts on $\W$ by derivations.
We denote by $\ti\W$ the extended annihilation algebra
$\dd\ltimes\W$, where
\begin{equation}\lbb{dactw}
[\d, x \tt a] = \d x \tt a \,, \qquad \d,a\in\dd, \; x\in X \,.
\end{equation}


We choose $L_0=\kk\tt\dd$ as a subspace of $\Wd$ such that
$\Wd=H L_0$, and we obtain the following filtration of $\W$:
\begin{equation}\lbb{wp}
\W_p = \fil_p \W = \fil_p X\tt_H L_0 \equiv \fil_p X\tt\dd \,.
\end{equation}
This is a decreasing filtration of $\W$, satisfying
$\W_{-1}=\W$ and \eqref{filbr} for $\ell=0$.
Note that $\W/\W_0 \simeq \kk\tt\dd \simeq \dd$ and
$\W_0/\W_1 \simeq \dd^*\tt\dd$.

Let us fix a basis $\{\d_i\}_{i=1,\dots,N}$ of $\dd$,
and let $x^i\in X$ be given by \eqref{xi2}.
We can view $x^i$ as elements of $\dd^*$;
then $\{x^i\}$ is a basis of $\dd^*$ dual to the basis $\{\d_i\}$ of $\dd$.
Let $e_i^j \in\gld$ be given by $e_i^j\d_k=\de^j_k\,\d_i$;
in other words, $e_i^j$ corresponds to $\d_i\tt x^j$
under the isomorphism $\gld\simeq\dd\tt\dd^*$.

\begin{lemma}\lbb{lwbra}
In the Lie algebra $\W=X\tt\dd$, we have the following\/{\rm:}
\begin{align*}
[x^j \tt \d_i, 1 \tt \d_k]
&= - \de^j_k\, 1 \tt \d_i \mod\W_0 \,,
\\
[x^j \tt \d_i, x^l \tt \d_k]
&= \de_i^l\, x^j \tt \d_k - \de^j_k\, x^l \tt \d_i \mod\W_1 \,.
\end{align*}
\end{lemma}
\begin{proof}
This follows from \eqref{Wbra} and \leref{ldacton}.
\end{proof}
\begin{corollary}\lbb{cwbra}
For $x\in \fil_0 X$, $a\in\dd$, the map
\begin{equation*}
x\tt a \mod \W_1 \mapsto -a \tt (x \mmod \fil_1 X)
\end{equation*}
is a Lie algebra isomorphism from
$\W_0/\W_1$ to\/ $\dd\tt\dd^*\simeq\gld$.
Under this isomorphism, the adjoint action of\/ $\W_0/\W_1$ on
$\W/\W_0$ coincides with the standard action of\/
$\gld$ on~$\dd$.
\end{corollary}
\begin{proof}
The above map takes $x^j \tt \d_i \mod\W_1$ to $-e_i^j\in\gld$.
\end{proof}

The action of $\Wd$ on $H$ induces a corresponding action
of the annihilation algebra $\W=\A(\Wd)$ on $\A(H)\equiv X$
given by \eqref{alav}:
\begin{equation}\lbb{wonx}
(x \tt a) y = - x(ya) \,,
\qquad x,y\in X , \; a\in\dd \,.
\end{equation}
Recall from \seref{shopf} that we have a ring isomorphism
$\ph\colon X \isoto \O_N$,
which is compatible with the corresponding filtrations
and topologies (see \eqref{phxon}, \eqref{filon}).
Since $\dd$ acts on $X$ by continuous derivations,
the Lie algebra $\W$ acts on $X$ by continuous derivations.
Hence, \eqref{wonx} defines a Lie algebra homomorphism
\begin{equation}\lbb{phwwn}
\ph\colon\W \to W_N \quad\text{such that}\quad
\ph(Ay) = \ph(A)\ph(y) \;\;\text{for}\;\; A\in\W \,, \; y\in X \,,
\end{equation}
where $W_N$ is the Lie algebra of continuous derivations of $\O_N$.

There is a natural filtration of $W_N$ given by
\begin{equation}\lbb{filpwn}
\fil_p W_N = \{ D\in W_N \st D (\fil_n \O_N)
\subset \fil_{n+p} \O_N
\;\;\text{ for all } n \}  , \qquad p\ge-1 \,.
\end{equation}
Explicitly, by \eqref{filon}, we have
\begin{equation}\lbb{filpwn2}
\fil_p W_N = \left\{ \sum_{i=1}^N f_i \frac\d{\d t^i}
\; \Big| \; f_i \in \fil_p \O_N \right\} \,.
\end{equation}
The filtration \eqref{filpwn} has the following important property
for $D\in W_N$:
\begin{equation}\lbb{filpwn4}
[D, \fil_p W_N] \subset \fil_{p+n} W_N
\iff D \in \fil_n W_N \,.
\end{equation}

\begin{proposition}\lbb{pwwn}
{\rm(i)}
We have{\rm:}
\begin{align*}
\ph(x \tt a) &= \ph(x) \ph(1 \tt a) \,, \qquad x\in X \,, \; a\in\dd \,,
\\
\ph(1\tt\d_i) &= -\frac\d{\d t^i} \mod\fil_0 W_N
\,, \qquad i=1,\dots,N \,.
\end{align*}

{\rm(ii)}
The homomorphism \eqref{phwwn} is an isomorphism and
$\ph(\W_p) = \fil_p W_N$ for all $p\ge-1$.
\end{proposition}
\begin{proof}
Part (i) follows from \eqref{phwwn} and \leref{ldacton}.
Part (ii) follows from (i) and \eqref{wp}, \eqref{filpwn2}.
\end{proof}

The adjoint action of the Euler vector field
\begin{equation}\lbb{euler}
E := \sum_{i=1}^N t^i \frac{\d}{\d t^i} \in \fil_0 W_N
\end{equation}
decomposes $W_N$ as a direct product of eigenspaces $W_{N;j}$ ($j \geq -1$),
on which the action of $E$ is multiplication by $j$.
One clearly has:
\begin{equation}\lbb{filpwn3}
\fil_p W_N = \prod_{j\geq p} W_{N;j}\,,
\quad \fil_p W_N / \fil_{p+1} W_N \simeq W_{N;p}\,.
\end{equation}
Notice that $W_{N;0} = \ker(\ad E)$ is a Lie algebra
isomorphic to $\gl_N$
and each space $W_{N;p}$ is a module over $W_{N;0}$.

\begin{definition}\lbb{deul}
The preimage $\E=\ph^{-1}(E) \in\W_0$ of the Euler vector field \eqref{euler}
under the isomorphism \eqref{phwwn}
will be called the {\em Euler element\/} of $\W$.
\end{definition}

By \prref{pwwn} and \coref{cwbra}, we have:
\begin{equation}\lbb{euler2}
\E = -\sum_{i=1}^N x^i \tt \d_i \mod\W_1
\,, \quad\text{i.e.,}\quad
\E\mmod\W_1 = \Id \in \gld \simeq \W_0/\W_1 \,.
\end{equation}

\subsection{The Normalizer $\N_{\W}$}\lbb{snw}
In this subsection, we study the normalizer of $\W_p$ $(p\ge0)$
in the extended annihilation algebra $\ti\W$.
These results will be used later in our classification of finite irreducible
$\Wd$-modules.

We denote by $\ad$ the adjoint action of $\dd$ on itself
(or on $H=\ue(\dd)$),
and by $\coad$ the coadjoint action of $\dd$ on $X=H^*$.
For $\d\in\dd$, we will also consider $\ad\d$ as an element of~$\gld$.
Note that, by \eqref{hx}, \eqref{xh}, we have
\begin{equation}\lbb{coaddx}
(\coad\d) x = \d x - x \d \,, \qquad \d\in\dd \,, \; x\in X \,.
\end{equation}
Since $\ad\d$ preserves the filtration \eqref{filued} of $H$,
it follows that $\coad\d$ preserves the filtration
$\{\fil_p X\}$ of $X$.

\begin{lemma}\lbb{lnw1}
{\rm{(i)}}
For $\d,a\in\dd$ and $x\in X$, the following formula holds in $\ti\W:$
\begin{equation*}
[\d + 1 \tt \d, x \tt a] = (\coad\d)x \tt a + x \tt [\d, a] \,.
\end{equation*}
In particular, the adjoint action of\/ $\d + 1 \tt \d \in\ti\W$ on
$\W\subset \ti \W$ preserves the filtration~$\{\W_p\}$.

{\rm{(ii)}}
The adjoint action of\/ $\d + 1 \tt \d$ on $\W/\W_0$
coincides with the standard action of\/ $\ad\d\in\gld$
on $\dd \simeq \W/\W_0$.
\end{lemma}
\begin{proof}
Part (i) follows from \eqref{Wbra}--\eqref{wp}, \eqref{coaddx},
and the above observation that $\coad\d$ preserves the filtration
$\{\fil_p X\}$ of $X$.
Part (ii) is obvious from (i).
\end{proof}

It is well known that all derivations of $W_N$ are inner.
Since $\W\simeq W_N$ and $\dd$ acts on $\W$ by derivations
(see \eqref{dactw}), there is an injective  Lie algebra homomorphism
\begin{equation}\lbb{gadw}
\ga\colon\dd\injto\W \quad\text{such that}\quad [\d,A]=[\ga(\d),A] \,,
\qquad \d\in\dd\subset\ti\W \,, \; A\in\W\subset\ti\W \,.
\end{equation}

\begin{definition}\lbb{dnw}
For $\d\in\dd$, let $\ti\d = \d-\ga(\d) \in\ti\W$,
where $\ga$ is given by \eqref{gadw}.
Let $\ti\dd = (\id-\ga)(\dd) \subset\ti\W$ and
$\N_\W = \ti\dd + \W_0 \subset\ti\W$.
\end{definition}
\begin{proposition}\lbb{pnw}
{\rm(i)}
The space $\ti\dd$ is a subalgebra of\/ $\ti\W$ centralizing $\W$, i.e.,
$[\ti\dd,\W] = \{0\}$.
The map $\d \mapsto \ti \d$ is a Lie algebra isomorphism from
$\dd$ to $\ti\dd$.

{\rm(ii)}
The space $\N_\W$ is a subalgebra of\/ $\ti\W$, and it decomposes
as a direct sum of Lie algebras,
$\N_\W = \ti\dd \oplus \W_0$.
\end{proposition}
\begin{proof}
It follows from \eqref{gadw} that
$[\ti\d,A] = 0$ for all $\d\in\dd$, $A\in\W$.
Then for $\d,\d'\in\dd$, we have
\begin{equation*}
[\d, \d'] = [\ti\d + \ga(\d), \ti\d' + \ga(\d')]
= [\ti\d, \ti\d'] + [\ga(\d), \ga(\d')] \,,
\end{equation*}
which implies $[\ti\d, \ti\d'] = \ti{[\d,\d']}$
since $\ga$ is a Lie algebra homomorphism.
This proves (i).
Part (ii) follows from (i) and \deref{dnw}.
\end{proof}
\begin{lemma}\lbb{lnw2}
For every $\d\in\dd$, the element $\d + 1 \tt \d - \ti\d \in\ti\W$
belongs to $\W_0$. Its image in $\W_0/\W_1$ coincides with
$\ad\d\in\gld \simeq \W_0/\W_1$.
%
\end{lemma}
\begin{proof}
First note that
$\d + 1 \tt \d - \ti\d = \ga(\d) + 1 \tt \d$
belongs to $\W$. By \eqref{gadw} and \leref{lnw1}(i), the adjoint
action of this element on $\W$ preserves the filtration $\{\W_p\}$.
Therefore, by \eqref{filpwn4}, $\ga(\d) + 1 \tt \d$ belongs to $\W_0$.
By \eqref{gadw} and \leref{lnw1}(ii), its image in
$\W_0/\W_1$ coincides with $\ad\d$.
%
\end{proof}

\begin{proposition}\lbb{normalizer}
For every $p\ge0$, the normalizer of\/ $\W_p$
in the extended annihilation algebra $\ti \W$
is equal to $\N_{\W}$. In particular, it is independent of $p$.
There is a decomposition as a direct sum of subspaces,
$\ti\W = \dd \oplus \N_{\W}$.
\end{proposition}
\begin{proof}
First, to show that $\ti\W = \dd \oplus \N_{\W}$, we have to check that
$\ti\W = \dd \oplus \ti\dd \oplus \W_0$ is a direct sum of subspaces.
This follows from \deref{dnw}, \leref{lnw2} and the
fact that $\ti\W = \dd \oplus \W$, $\W = (\kk\tt\dd) \oplus \W_0$
as vector spaces.

Next, it is clear that $\N_\W$ normalizes $\W_p$\,, because
$[\ti\dd,\W_p] = \{0\}$ and $[\W_0,\W_p] \subset \W_p$.
Assume that an element $\d\in\dd$ normalizes $\W_p$.
By \eqref{dactw}, we obtain that in this case
$\d(\fil_p X) \subset \fil_p X$.
However, one can deduce from \leref{ldacton} that
$\d(\fil_p X) = \fil_{p-1} X$, which is strictly larger than $\fil_p X$.
This contradiction shows that the normalizer of $\W_p$
is equal to $\N_\W$.
\end{proof}

In order to understand the irreducible representations of $\N_{\W}$,
we need the following lemma, which appeared
(in the more difficult super case) in \cite[Erratum]{CK}.

\begin{lemma}\lbb{trivialaction}
Let\/ $\g$ be a finite-dimensional Lie algebra, and let\/ $\g_0 \subset \g$
be either a simple Lie algebra or a $1$-dimensional Lie algebra.
Let $I$ be a subspace of the radical of\/ $\g$,
stabilized by $\ad \g_0$ and having the property that\/
$[\g_0, a] = 0$ for $a \in I$ implies $a=0$.
Then $I$ acts trivially on any irreducible finite-dimensional\/
$\g$-module $V$.
\end{lemma}
\begin{proof}
By Cartan--Jacobson's Theorem (see, e.g., \cite[Theorem~VI.5.1]{Se}),
every $a\in\Rad \g$ acts by scalar multiplication on $V$.
Let $J=\{ a \in I \st a(V) = 0\}$. Then $[\g_0, I] \subset J$.

Now, if $\g_0$ is simple, then $J$ is a $\g_0$-submodule of $I$
and, by complete reducibility, $I = J \oplus J^\perp$ as
$\g_0$-modules for some complement $J^\perp$. Hence,
$[\g_0,J^\perp] = 0$, so $J^\perp=0$ and $I=J$.

If instead $\g_0 = \kk e$ is $1$-dimensional, then $[e, I]\subset J$.
If $J\neq I$, then $\ad e \colon I\to J$ is not injective,
which is a contradiction. We conclude that $J=I$.
\end{proof}

An $\N_\W$-module $V$ will be called {\em conformal\/}
if it is conformal as a module over the subalgebra $\W_0\subset\N_\W$,
i.e., if every vector $v\in V$ is killed by some $\W_p$.

\begin{proposition}\lbb{Nreps}
The subalgebra $\W_1 \subset \N_\W$ acts trivially on any
irreducible finite-dimensional conformal $\N_\W$-module.
Irreducible finite-dimensional conformal $\N_\W$-modules
are in one-to-one correspondence with irreducible
finite-dimensional modules over the Lie algebra
$\N_\W/\W_1 \simeq \dd\oplus\gld$.
\end{proposition}
\begin{proof}
A finite-dimensional vector space $V$ is a conformal $\N_\W$-module
iff it is an $\N_\W$-module on which $\W_p$ acts trivially for some $p\ge0$,
i.e., iff it is a module over the finite-dimensional Lie algebra
$\g = \N_{\W}/\W_p = \ti\dd\oplus(\W_0/\W_p)$.

We will apply \leref{trivialaction} for $I=\W_1/\W_p$ and
$\g_0 = \kk\E \mmod\W_p \subset \W_0/\W_p$\,,
where $\E\in\W_0$ is the Euler element
(see \deref{deul}).
Note that $I \subset \Rad\g$ and $[\E,I] \subset I$,
because $[\W_i, \W_j] \subset \W_{i+j}$ for all $i,j$.
The adjoint action of $\E$ is injective on $I$, because
$\ad E$ is injective on
$\fil_1 W_N / \fil_p W_N = \prod_{j=1}^{p-1} W_{N;j}$
(see \eqref{filpwn3}).
We conclude that $I$ acts trivially on any finite-dimensional
conformal $\N_\W$-module. Hence, we can take $p=1$.
Then $\g = \ti\dd\oplus(\W_0/\W_1) \simeq \dd\oplus\gld$,
since $\ti\dd\simeq\dd$ and $\W_0/\W_1 \simeq\gld$.
\end{proof}

\subsection{Annihilation Algebra of $\Sd$}\lbb{subans}
Assume that $N=\dim\dd>2$. In this subsection, we study
the annihilation algebra $\S = \A(\Sd) := X\tt_H \Sd$ of the
Lie \psalg\ $\Sd$ defined in \seref{subprim}.
Our treatment here is more detailed than in \cite[Section~8.4]{BDK}.

We choose
\begin{equation}\lbb{sl0}
L_0=\Span_\kk \{ s_{ab} \st a,b \in\dd \} \subset \Sd
\end{equation}
as a subspace such that $\Sd=H L_0$, where
the elements $s_{ab}$ are given by \eqref{achib}.
We obtain a decreasing filtration of $\S$:
\begin{equation}\lbb{sp}
\S_p = \fil_{p+1} \S = \fil_{p+1} X\tt_H L_0 \,, \qquad p\ge -2 \,,
\end{equation}
satisfying
$\S_{-2}=\S$ and \eqref{filbr} for $\ell=1$.
Then $[\S_n, \S_p] \subset \S_{n+p}$ for all $n,p$.

The canonical injection of the subalgebra $\Sd$ into $\Wd$
induces a Lie algebra homomorphism $\io\colon\S\to\W$.
Explicitly, we have:
\begin{equation}\lbb{ioxs}
\begin{split}
\io(x \tt_H s) &= \sum\, x h_i\tt\d_i \in\W \equiv X\tt\dd
\\
&\text{for}\quad x\in X \,, \;
s=\sum\, h_i\tt\d_i \in\Sd \subset \Wd=H\tt\dd \,.
\end{split}
\end{equation}
Here, as before, we identify $\W=X\tt_H\Wd$ with $X\tt\dd$.

We define a map $\Div^\chi\colon\W\to X$ by the formula
(cf.\ \eqref{div*}):
\begin{equation}\lbb{divw}
\Div^\chi \left( \sum\, y_i\tt \d_i \right)
= \sum\, y_i (\d_i + \chi(\d_i)) \,.
\end{equation}
It is easy to see that 
\begin{equation}\lbb{divab}
\Div^\chi[A,B] = A (\Div^\chi B) - B (\Div^\chi A) \,,
\qquad A,B \in \W \,,
\end{equation}
where the action of $\W$ on $X$ is given by \eqref{wonx}.
This implies that
\begin{equation}\lbb{sw1}
\ov\S := \{ A \in\W \st \Div^\chi A = 0 \}
\end{equation}
is a Lie subalgebra of $\W$.
It was shown in \cite[Section~8.4]{BDK}
that $\ov\S$ is isomorphic to the Lie algebra
of divergence-zero vector fields
\begin{equation}\lbb{sw2}
S_N := \left\{ \sum_{i=1}^N f_i \frac{\d}{\d t^i} \in W_N
\;\Big|\;  \sum_{i=1}^N \frac{\d f_i}{\d t^i} = 0 \right\} \,.
\end{equation}

\begin{lemma}\lbb{lsw2}
If\/ $N=\dim\dd>2$, the map \eqref{ioxs} is an embedding of Lie algebras
$\io\colon \S \injto \ov\S$.
\end{lemma}
\begin{proof}
It follows from \eqref{div*}, \eqref{ioxs} and \eqref{divw} that
\begin{equation*}
x(\Div^\chi s) = \Div^\chi \io(x \tt_H s)
\,, \qquad x \in X \,, \; s\in\Wd \,.
\end{equation*}
Therefore, $\io(\S)$ is contained in $\ov\S$.
Next, note that for $N>2$, $\S$ is isomorphic to $S_N$ by
\cite[Theorem~8.2(i)]{BDK}. It is well known that the Lie algebra
$S_N$ is simple; hence, $\S$ is simple. Since $\io$ is a nonzero
homomorphism, it must be injective.
\end{proof}
\begin{remark}\lbb{rsn=2}
When $N=\dim\dd=2$, the Lie algebra $\S$ is isomorphic to $P_2$,
which is an extension of $S_2=H_2$ by a $1$-dimensional center
(cf.\ \reref{rdimd2}).
In this case, the homomorphism \eqref{ioxs} has a $1$-dimensional 
kernel.
\end{remark}

We will prove in \prref{psw2} below that, in fact, $\io(\S)=\ov\S$.
Recall that we have a Lie algebra isomorphism $\ph\colon\W\isoto W_N$,
given by \eqref{phwwn}. However, although $\ov\S \simeq S_N \subset W_N$,
the image $\ph(\ov\S) \subset W_N$ is not equal to $S_N$ in general.
Instead, we will show that the images of $\ph(\ov\S)$ and $S_N$ coincide in
the associated graded algebra of $W_N$ (see \prref{psw1} below).

\begin{lemma}\lbb{lsw3}
For every $p\ge-1$, we have
\begin{equation*}
\ph(\ov\S \cap\fil_p\W) \subset
(S_N \cap\fil_p W_N) + \fil_{p+1} W_N \,.
\end{equation*}
\end{lemma}
\begin{proof}
Take an element $A = \sum\, y_i\tt \d_i \in \fil_p\W$; then each
$y_i \in\fil_p X$. By \prref{pwwn}, we have
$\ph(A) = \sum\, f_i \, \ph(1\tt\d_i)$,
where $f_i = \ph(y_i) \in \fil_p\O_N$.
Since $\ph(1\tt\d_i) = -\d / \d t^i \mod \fil_0 W_N$, we have
$\ph(A) = -\sum\, f_i \, \d / \d t^i \mod \fil_{p+1} W_N$.
It follows from \eqref{divw} and \leref{ldacton} that
$\ph(\Div^\chi A) = -\sum\, \d f_i / \d t^i \mod \fil_p \O_N$.
If $A\in\ov\S\cap\fil_p\W$, then
$\sum\, \d f_i / \d t^i = 0 \mod \fil_p \O_N$.
Then there exist elements $\hat f_i \in\fil_p \O_N$ such that
$\hat f_i = f_i \mod \fil_{p+1} \O_N$ and
$\sum\, \d \hat f_i / \d t^i = 0$.
This means that
$\hat A := -\sum\, \hat f_i \, \d / \d t^i \in S_N \cap\fil_p W_N$
and $\ph(A) = \hat A \mod \fil_{p+1}W_N$.
\end{proof}

Consider the associated graded of $\W$,
\begin{equation}\lbb{grw1}
\gr \W := \bigoplus_{p=-1}^\infty \gr_p \W \,, \qquad
\gr_p \W := \fil_p \W / \fil_{p+1} \W \,.
\end{equation}
Note that, by \eqref{wp}, we have $\gr_p \W = (\gr_p X) \tt\dd$.
Similarly, we have
$\gr_p W_N = \sum_{i=1}^N \, (\gr_p \O_N) \, \d / \d t^i$.
The maps $\ph\colon X\to\O_N$ and $\ph\colon\W\to W_N$
(see \eqref{phxon}, \eqref{phwwn})
preserve the corresponding filtrations and induce maps
$\gr\ph\colon \gr X\to\gr\O_N$ and $\gr\ph\colon\gr\W\to \gr W_N$.
Note also that the map $\Div^\chi\colon \W\to X$ takes
$\fil_p\W$ to $\fil_{p-1} X$, and hence induces a map
$\gr\Div^\chi\colon \gr\W\to \gr X$ of degree $-1$.
The same is true for the map
$\Div\colon W_N\to\O_N$ given by
$\Div( \sum\, f_i \, \d / \d t^i ) := \sum\, \d f_i / \d t^i$.
{}From the proof of \leref{lsw3} we deduce:

\begin{corollary}\lbb{csw4}
The above maps satisfy
\begin{equation*}
(\gr\ph)\Bigl( \sum_{i=1}^N \bar y_i\tt \d_i \Bigr)
= -\sum_{i=1}^N (\gr\ph)(\bar y_i) \, \frac{\d}{\d t^i} \,,
\qquad \bar y_i \in \gr X
\end{equation*}
and
\begin{equation*}
\gr\ph \circ \gr\Div^\chi = \gr\Div \circ \gr\ph \,.
\end{equation*}
\end{corollary}

The Lie algebra $\S$ has a filtration \eqref{sp}, while
$\ov\S\subset\W$ has one obtained by restricting the
filtration \eqref{wp} of $\W$.
Using \leref{lsw3}, we can prove that $\io$ is
compatible with the filtrations.


\begin{proposition}\lbb{psw2}
Let\/ $\S$ be the annihilation algebra of\/ $\Sd$, and let\/
$\ov\S\subset\W$ be defined by \eqref{sw1}. Then for\/
$\dim\dd>2$, the map \eqref{ioxs} is an isomorphism of Lie algebras
$\io\colon \S \isoto \ov\S$ such that $\io(\S_p) = \ov\S\cap\W_p$
for all $p\ge-1$.
\end{proposition}
\begin{proof}
It is clear from definitions that
\begin{equation*}
\io(\S_p) = \fil_{p+1} X\tt_H \Span_\kk \{ s_{ab} \}
\subset \fil_p X\tt_H (\kk\tt\dd)
\equiv \fil_p X\tt\dd = \W_p \,.
\end{equation*}
In addition, $\io(\S) \subset \ov\S$ by \leref{lsw2};
hence, $\io(\S_p) \subset \ov\S\cap\W_p$.

Conversely, let $A\in\ov\S\cap\W_p$. By \leref{lsw3}, we can find
$\hat A \in S_N \cap\fil_p W_N$ such that
$\ph(A) = \hat A \mod \fil_{p+1}W_N$.
Any element of $S_N \cap\fil_p W_N$
can be written in the form
\begin{equation*}
\hat A = \sum_{i,j=1}^N
\frac{\d f_{ij}}{\d t^i} \frac{\d}{\d t^j}
- \frac{\d f_{ij}}{\d t^j} \frac{\d}{\d t^i} \,,
\qquad f_{ij} \in\fil_{p+1}\O_N \,.
\end{equation*}
Now consider the following element of $\S_p$:
\begin{equation*}
\tilde A := -\sum_{i,j=1}^N
y_{ij} \tt_H s_{\d_i \d_j} \,,
\qquad y_{ij} := \ph^{-1}(f_{ij}) \in\fil_{p+1} X \,.
\end{equation*}
Then we have $\tilde A \in\S_p$
and $\io(\tilde A) = A \mod \W_{p+1}$.

Let $A_1 = A - \io(\tilde A)$; then $A_1 \in \ov\S\cap\W_{p+1}$
and $A-A_1 \in\io(\S_p)$. By the above argument,
we can find an element $\tilde A_1 \in \S_{p+1}$ such that
$\io(\tilde A_1) = A_1 \mod \W_{p+2}$.
Let $A_2 = A_1 - \io(\tilde A_1)$; then $A_2 \in \ov\S\cap\W_{p+2}$
and $A_1-A_2 \in \io(\S_{p+1})$.
Continuing this way, we obtain a sequence of elements
$A_n \in \ov\S\cap\W_{p+n}$ such that
$A_n-A_{n+1} \in \io(\S_{p+n})$ for all $n\ge0$,
where $A_0 :=A$.
The sequence $A_n$ converges to $0$ in $\W$ and
$A-A_n \in \io(\S_p)$ for all $n\ge0$; therefore,
$A\in\io(\S_p)$.

This proves that $\io(\S_p) = \ov\S\cap\W_p$.
Taking $p=-1$, we get
$\io(\S) \supset \io(\S_{-1}) = \ov\S$,
because $\W_{-1}=\W\supset\ov\S$. Now \leref{lsw2} implies that
$\io$ is an isomorphism.
\end{proof}

Recall that any ring automorphism $\psi$ of $\O_N$ induces a
Lie algebra automorphism $\psi$ of $W_N = \Der\O_N$ such that
$\psi(Ay)=\psi(A)\psi(y)$ for $A\in W_N$, $y\in\O_N$.
Any $\psi\in\Aut\O_N$ preserves the filtration,
because $\fil_0 \O_N$ is the unique maximal ideal of $\O_N$ and
$\fil_p \O_N = (\fil_0 \O_N)^{p+1}$ for $p\ge0$ (see \eqref{filon}).
Then it follows from \eqref{filpwn} that $\psi$
preserves the filtration $\{\fil_p W_N \}$.

\begin{proposition}\lbb{psw1}
There exists a ring automorphism $\psi$ of\/ $\O_N$
such that the induced Lie algebra automorphism $\psi$ of\/ $W_N$
satisfies
$\ph(\ov\S) = \psi(S_N)$ and
\begin{equation}\lbb{sw5}
(\psi-\id)(\fil_p W_N) \subset \fil_{p+1} W_N
\,, \qquad p\ge-1 \,.
\end{equation}
\end{proposition}
\begin{proof}
In \cite[Remark 8.2]{BDK} the image $\ph(\ov\S)$ is described as
the Lie algebra of all vector fields annihilating a certain volume
form. But any two volume forms are related by a change of variables,
i.e., by a ring automorphism of $\O_N$,
and the subalgebra $S_N$ corresponds to the standard volume form
$\di t^1 \wedge\cdots\wedge \di t^N$. Hence, there exists an
automorphism $\psi$ of $\O_N$ such that $\ph(\ov\S) = \psi(S_N)$.
Due to \coref{csw4}, we can choose $\psi$ such that
\begin{equation*}
\psi(t^i)=t^i \mod \fil_1\O_N \,, \qquad i=1,\dots,N \,,
\end{equation*}
i.e., such that $\gr\psi=\id$.
Since the latter is equivalent to \eqref{sw5},
this completes the proof.
\end{proof}

\begin{corollary}\lbb{csw7}
The Lie algebra isomorphism $\psi^{-1}\ph\io \colon \S \isoto S_N$
maps $\S_p$ onto $S_N \cap\fil_p W_N$ for all $p\ge-1$.
In particular, $\S_{-2}=\S_{-1}=\S$.
\end{corollary}
\begin{proof}
By \prref{psw2}, $\io(\S_p) = \ov\S\cap\W_p$. Then under the isomorphism
$\ph\colon\W\to W_N$, we have $\ph\io(\S_p) = \ph(\ov\S)\cap\fil_p W_N$.
But, by \prref{psw1}, $\ph(\ov\S)=\psi(S_N)$ and
$\psi(\fil_p W_N) = \fil_p W_N$; hence,
$\ph\io(\S_p) = \psi(S_N \cap\fil_p W_N)$.
\end{proof}

Recall that $W_{N;p}$ is the subspace of $W_N$ on which
the adjoint action of the Euler vector field \eqref{euler}
is multiplication by $p$.
We let $S_{N;p} = S_N \cap W_{N;p}$.
Since $S_N$ is preserved by $\ad E$,
it admits a decomposition similar to \eqref{filpwn3}:
\begin{equation}\lbb{filpsn}
S_N \cap\fil_p W_N = \prod_{j\geq p} S_{N;j}\,,
\quad (S_N \cap\fil_p W_N) / (S_N \cap\fil_{p+1} W_N) \simeq S_{N;p}\,.
\end{equation}
The following facts about the Lie algebra $S_N\subset W_N$ are well known.

\begin{lemma}\lbb{lsw7}
{\rm(i)}
The Lie algebra $S_{N;0}$ is isomorphic to $\sl_N$.

{\rm(ii)}
For every $p\ge-1$, the $S_{N;0}$-module $S_{N;p}$ is isomorphic to the highest
component of the $\sl_N$-module $\kk^N \tt (\symp^{p+1}\kk^N)^*$.
In particular, $S_{N;p}$ has no trivial components
in its decomposition as a sum of irreducible $\sl_N$-modules.

{\rm(iii)}
The normalizer of\/ $S_N$ in $W_N$ is $S_N \oplus \kk E$.
%
\end{lemma}

\begin{definition}\lbb{deul2}
We let $\what E = \psi(E) \in W_N$ and
$\what\E = \ph^{-1}(\what E) \in \W$,
where $E$ is the Euler vector field \eqref{euler},
$\ph$ is from \eqref{phwwn} and $\psi$ is from \prref{psw1}.
\end{definition}

Combining the above results with \eqref{euler2},
we obtain the following corollary.

\begin{corollary}\lbb{csw6}
{\rm(i)}
The Lie algebra $\S_0 / \S_1$ is isomorphic to $\sld$.

{\rm(ii)}
For every $p\ge-1$, the $(\S_0 / \S_1)$-module
$\S_p / \S_{p+1}$ has no trivial $\sld$-components.

{\rm(iii)}
The normalizer of\/ $\ov\S$ in $\W$ is $\ov\S + \kk\what\E$.

{\rm(iv)}
$\what\E$ belongs to $\W_0$ and
its image in $\W_0/\W_1$ coincides with\/ $\Id\in\gld \simeq \W_0/\W_1$.
\end{corollary}

\subsection{The Normalizer $\N_{\S}$}\lbb{sns}
In this subsection, we study the normalizer of $\S_p$ $(p\ge0)$
in the extended annihilation algebra $\ti\S = \dd\ltimes\S$.
We will use extensively the results and notation of
Sections \ref{snw} and \ref{subans}, and we will identify
$\S$ with the subalgebra $\ov\S$ of $\W$ (see \prref{psw2}).

Recall that the filtration $\{\S_p\}$ of $\S$ has the properties:
$\S_{-2}=\S_{-1}=\S$ and $[\S_n,\S_p] \subset \S_{n+p}$
for all $n,p$. In addition, by \coref{csw6}, we have:
$\W_0 = \S_0 + \kk\what\E + \W_1$,
where the element $\what\E \in\W_0$
is from \deref{deul2}.

\begin{lemma}\lbb{lns1}
For every $\d\in\dd$, we have{\rm:}
$1 \tt \d - (\chi(\d)/N) \, \what\E \in \S+\W_1$.
\end{lemma}
\begin{proof}
As before, let $\{\d_i\}_{i=1,\dots,N}$ be a basis of $\dd$,
and let $x^i\in X$ be given by \eqref{xi2}.
Denote by $c_{ij}^k$ the structure constants of $\dd$
in the basis $\{\d_i\}$,
and let $\chi_i = \chi(\d_i)$ for short.
Using \eqref{achib}, \eqref{ioxs} and \leref{ldacton}, we compute
for $i<j$:
\begin{align*}
\io( x^i \tt_H s_{\d_i,\d_j} )
=
\chi_i x^i \tt \d_j &- \chi_j x^i \tt \d_i - x^i \tt [\d_i, \d_j]
+ x^i \d_i \tt \d_j - x^i \d_j \tt \d_i
\\
=
\chi_i x^i \tt \d_j &- \chi_j x^i \tt \d_i
- \sum_k \, c_{ij}^k x^i \tt \d_k
- 1 \tt \d_j
\\
&+ \sum_{k>i} \, c_{ik}^i x^k \tt \d_j
- \sum_{k>j} \, c_{jk}^i x^k \tt \d_i
\mod\W_1 \,.
\end{align*}
{}From here, we see that the element
$\io( x^i \tt_H s_{\d_i,\d_j} ) + 1 \tt \d_j$
belongs to $\W_0$. Next, using \coref{cwbra}, we find that the
image of this element in $\W_0 / \W_1 \simeq\gld$ has trace $\chi_j$.
Therefore, by \coref{csw6} (i), (iv),
\begin{equation*}
\io( x^i \tt_H s_{\d_i,\d_j} ) + 1 \tt \d_j - (\chi_j/N) \, \what\E
\in \S_0 + \W_1 \,,
\end{equation*}
which implies $1 \tt \d_j - (\chi_j/N) \, \what\E \in \S + \W_1$.
\end{proof}

\begin{lemma}\lbb{lns2}
For every $\d\in\dd$, we have{\rm:}
$\ga(\d) + 1 \tt \d - (\tr\ad(\d)/N) \, \what\E \in \S_0+\W_1$,
where $\ga$ is from \eqref{gadw}.
\end{lemma}
\begin{proof}
By \leref{lnw2}, $\ga(\d) + 1 \tt \d \in \W_0$
and its image in $\W_0/\W_1$ coincides with $\ad\d\in\gld$.
Now the statement follows from \coref{csw6} (i), (iv).
\end{proof}

\begin{definition}\lbb{dns}
For $\d\in\dd$, let
\begin{equation*}
\what\ga(\d) = \ga(\d) + \bigl( (\chi - \tr \ad)(\d)/N \bigr) \what\E
\in\W \,,
\end{equation*}
where $\ga$ is given by \eqref{gadw}.
Let $\what\d = \d - \what\ga(\d)$,
$\what\dd = (\id-\what\ga)(\dd) \subset\ti\W$, and
$\N_\S = \what\dd + \S_0 \subset\ti\W$.
\end{definition}

Note that
\begin{equation}\lbb{ns1}
\what\d = \ti\d - \bigl( (\chi - \tr \ad)(\d)/N \bigr) \what\E \,,
\qquad \d\in\dd \,,
\end{equation}
where $\ti\d$ is from \deref{dnw}.

\begin{proposition}\lbb{pns3}
{\rm(i)}
We have
$\what\ga(\dd) \subset \S$
and\/
$\what\dd \subset \ti\S$.

{\rm(ii)}
The map $\d\mapsto\what\d$ is a Lie algebra isomorphism
from $\dd$ to $\what\dd$.

{\rm(iii)}
The Lie algebra $\what\dd$ normalizes $\S_p$ for all $p\ge-1$.

{\rm(iv)}
The Lie algebra $\what\dd$ centralizes $\S_0/\S_1$.
\end{proposition}
\begin{proof}
(i)
Combining Lemmas \ref{lns1} and \ref{lns2}, we get
$\what\ga(\d) \in \S+\W_1$ for all $\d\in\dd$.
On the other hand, we deduce from
\eqref{gadw} and \coref{csw6}(iii) that
$\what\ga(\d)$ normalizes $\S$.
Hence, again by \coref{csw6}(iii), $\what\ga(\d) \in \S+\kk\what\E$.
However, the intersection $(\S + \W_1) \cap (\S + \kk\what\E)$
is equal to $\S$. This shows that $\what\ga(\d)\in \S$.

(ii)
Recall from \seref{snw} that $\d\mapsto\ti\d$ is a Lie algebra isomorphism
and $\ti\dd\subset\ti\W$ centralizes $\W$. Then part (ii)
follows from \eqref{ns1}
and the fact that $\chi - \tr \ad$ is a trace form on~$\dd$.

(iii) and (iv)
follow from \eqref{ns1}, \coref{csw6} (iii), (iv) and
$[\ti\dd,\W]=0$.
\end{proof}

It follows from \prref{pns3} that $\N_S$ is a Lie subalgebra of $\ti\S$,
isomorphic to the semidirect sum $\what\dd\ltimes\S_0$.

\begin{proposition}\lbb{normalizerS}
For every $p\ge0$, the normalizer of\/ $\S_p$
in the extended annihilation algebra $\ti\S$
is equal to $\N_\S$. In particular, it is independent of $p$.
There is a decomposition as a direct sum of subspaces,
$\ti\S = \dd \oplus \N_{\S}$.
\end{proposition}
\begin{proof}
The proof is similar to that of \prref{normalizer}.
\end{proof}

An $\N_\S$-module $V$ is called {\em conformal\/}
if it is conformal as a module over the subalgebra $\S_0\subset\N_\S$,
i.e., if every vector $v\in V$ is killed by some $\S_p$.

\begin{proposition}\lbb{NrepstypeS}
The subalgebra $\S_1 \subset \N_\S$ acts trivially on any
irreducible finite-dimensional conformal $\N_\S$-module.
Irreducible finite-dimensional conformal $\N_\S$-modules
are in one-to-one correspondence with irreducible
finite-dimensional modules over the Lie algebra
$\N_\S / \S_1 \simeq \dd \oplus \sld$.
\end{proposition}
\begin{proof}
As in \prref{Nreps}, the $\N_\S$-action factors via the
finite-dimensional Lie algebra $\g := \N_\S/\S_p$ for some $p\ge1$.
Recall that $[\S_i,\S_j] \subset \S_{i+j}$ for all $i,j$, so that
$I := \S_1/\S_p$ is contained in the radical of 
$\ov\g := \S_0/\S_p \subset \g$.
Moreover, the quotient $\ov\g/I \simeq \S_0/\S_1$ is isomorphic to the
simple Lie algebra $\sld$ by \coref{csw6}(i).
Therefore, $I$ coincides with the radical of $\ov\g$,
and we can lift $\ov\g/I$ to a subalgebra $\g_0$ of $\ov\g$
isomorphic to $\sld$.
Then $I$ is contained in the radical of $\g$, and the
adjoint action of $\g_0$ on $\g$ preserves it.
Moreover, by \coref{csw6}(ii), $I$ has no trivial $\g_0$-components.
We can now apply \leref{trivialaction} to deduce
that $I$ acts trivially on any irreducible finite-dimensional conformal
$\N_\S$-module. Therefore, the $\N_\S$-action factors via $\N_\S/\S_1$.
By \prref{pns3}(iv), $\what\dd$ centralizes $\S_0/\S_1$.
Hence, $\N_\S/\S_1$ is isomorphic to a direct sum of Lie algebras
$\what \dd \oplus (\S_0/\S_1) \simeq \dd \oplus \sld$.
\end{proof}

\section{Pseudo Linear Algebra}\lbb{spslin}
In this section, we generalize several linear algebra
constructions to the \psalg\ context.
We introduce an important class of $\Wd$-modules
called tensor modules.

\subsection{Pseudolinear Maps}\lbb{spslm}
The definition of a module over a pseudoalgebra motivates the following
definition of a pseudolinear map.
\begin{definition}[\cite{BDK}]\lbb{dchom}
Let $V$ and $W$ be two $H$-modules. An {\em $H$-pseudolinear map\/}
from $V$ to $W$ is a $\kk$-linear map
$\phi\colon V\to (H\tt H)\tt_H W$
such that
\begin{equation}\lbb{hclm}
\phi(hv) = ((1\tt h)\tt_H 1) \, \phi(v),
\qquad h\in H, v\in V.
\end{equation}
We denote the space of all such $\phi$ by $\Chom(V,W)$.
We will also use the notation $\phi*v \equiv \phi(v)$
for $\phi\in\Chom(V,W)$, $v\in V$.
We define a left action of $H$ on $\Chom(V,W)$ by:
\begin{equation}\lbb{hactchom}
(h\phi)(v) = ((h\tt1)\tt_H 1)\, \phi(v).
\end{equation}
When $V=W$, we set $\Cend V=\Chom(V,V)$.
\end{definition}
\begin{example}\lbb{ecendm}
Let $A$ be an $H$-\psalg, and let $V$ be an $A$-module. 
Then for every $a\in A$
the map $m_a\colon V\to (H\tt H)\tt_H V$ defined by $m_a(v) = a*v$
is an $H$-pseudolinear map. Moreover, we have $h \, m_a = m_{ha}$
for $h\in H$.
\end{example}
\begin{remark}\lbb{rchomfun}
Given two homomorphisms of left $H$-modules
$\be\colon V'\to V$ and $\ga\colon W\to W'$,
we define a homomorphism
\begin{equation}\lbb{chalbe1}
\Chom(\be,\ga) \colon \Chom(V,W) \to \Chom(V',W')
\end{equation}
by the formula
\begin{equation}\lbb{chalbe2}
\phi \mapsto ( (\id\tt\id)\tt_H\ga ) \circ\phi\circ\be \,.
\end{equation}
Then we can view $\Chom(-,-)$ as a bifunctor
from the category of left $H$-modules to itself,
contravariant in the first argument and covariant in the second one.
\end{remark}

Recall from \cite[Chapter~10]{BDK} that when $V$ is a finite
$H$-module, $\Cend V$
has a unique structure of an \as\ \psalg\ such that $V$ is a module
over it via $\phi*v=\phi(v)$. 
Denote by $\gc V$ the Lie \psalg\ obtained from $\Cend V$
by the construction of \reref{assolie*}. Then $V$ is also
a module over $\gc V$.

\begin{proposition}[\cite{BDK}]\lbb{pcend}
Let\/ $L$ be a Lie \psalg, and let\/ $V$ be a finite\/ $H$-module.
Then giving a structure of an\/ $L$-module on\/ $V$ is equivalent to giving
a homomorphism of Lie \psalgs\ from\/ $L$ to~$\gc V$.
\end{proposition}
\begin{proof}
If $V$ is a finite $L$-module, we define a map $\rho\colon L \to \gc V$ by
$a\mapsto m_a$, where $m_a$ is from \exref{ecendm}.
Then $\rho$ is a homomorphism of Lie \psalgs\
(cf.\ \cite[Proposition~10.1]{BDK}). Conversely, given a homomorphism
$\rho\colon L \to \gc V$, we define an action of $L$ on $V$ by
$a*v = \rho(a)*v$.
\end{proof}

In the case when $V$ is a free $H$-module of finite rank, one can give an
explicit description of $\Cend V$, and hence of $\gc V$, as follows
(see \cite[Proposition~10.3]{BDK}).
Let $V=H\tt V_0$, where $H$ acts trivially on $V_0$ and $\dim V_0 < \infty$.
Then $\Cend V$
is isomorphic to $H\tt H\tt\End V_0$, with $H$ acting by left
multiplication on the first factor,
and with the following pseudoproduct:
\begin{equation}\lbb{prodcend}
(f\tt a\tt A)*(g\tt b\tt B)
= (f \tt g a_{(1)}) \tt_H (1 \tt b a_{(2)} \tt AB) \,.
\end{equation}
The action of\/ $\Cend V$ on $V=H\tt V_0$ is given by{\rm:}
\begin{equation}\lbb{actcendv}
(f\tt a\tt A)*(h\tt v) = (f \tt ha) \tt_H (1 \tt Av) \,.
\end{equation}
The pseudobracket in $\gc V$ is given by:
\begin{equation}\lbb{brackgc}
\begin{split}
[(f\tt a\tt A)*(g\tt b\tt B)]
&= (f \tt g a_{(1)}) \tt_H (1 \tt b a_{(2)} \tt AB)
\\
&- (f b_{(1)} \tt g) \tt_H (1 \tt a b_{(2)} \tt BA) \,.
\end{split}
\end{equation}
The action of\/ $\gc V$ on $V$ is also given by \eqref{actcendv}.


\begin{remark}\lbb{rcend}
Let $L$ be a Lie 
\psalg. Let $V=H\tt V_0$ be a finite $L$-module,
which is free as an $H$-module.
For all $a\in L$, $v\in V_0$ we can write
\begin{equation}\lbb{modhom1}
a*(1\tt v) = \sum\, (f_i \tt g_i) \tt_H (1 \tt A_i v) \,,
\end{equation}
where $f_i, g_i \in H$, $A_i \in \End V_0$.
Then the homomorphism $L\to\gc V$ 
is given by
$a \mapsto \sum f_i \tt g_i \tt A_i$.
This follows from \eqref{actcendv} and the proof of \prref{pcend}.
\end{remark}


\begin{example}\lbb{ewdgc}
{\rm(i)}
The action \eqref{wdac*} of $\Wd$ on $H$ gives an embedding
of Lie \psalgs\ $\Wd\injto\gc H = H\tt H$,
$f\tt a \mapsto -f\tt a$ ($f\in H$, $a\in\dd\subset H$).

{\rm(ii)}
Consider the semidirect sum $H\rtimes\Wd$, where $H$ is regarded
as a commutative Lie \psalg\ and $\Wd$ acts on $H$ via \eqref{wdac*}.
Then we have an embedding $H\rtimes\Wd\injto\gc H$ given by
$g+f\tt a \mapsto g\tt1-f\tt a$ for $f,g\in H$, $a\in\dd\subset H$.
%
\end{example}
\begin{remark}\lbb{rcurgtogc}
For any Lie algebra $\fg$, we have a semidirect sum
$\Cur\fg\rtimes\Wd$, where $\Cur\fg$ is defined in \exref{ecurps}
and $\Wd$ acts on $\Cur\fg=H\tt\fg$ via
\begin{equation*}
(f\tt a)*(g\tt B) = -(f\tt ga)\tt_H(1\tt B) \,, \qquad
f,g\in H \,, \; a\in\dd \,, \; B\in\fg \,.
\end{equation*}
Let $V_0$ be a finite-dimensional $\fg$-module, and
let $\rho$ be the corresponding homomorphism $\fg\to\gl\, V_0$.
Then we have a homomorphism of Lie \psalgs\
$\Cur\fg\rtimes\Wd\to\gc(H\tt V_0)$, given by
\begin{equation}\lbb{curgtogc2}
g\tt B + f\tt a \mapsto g\tt1\tt \rho(B) - f\tt a\tt\Id \,.
\end{equation}
\end{remark}

\subsection{Duals and Twistings of Representations}\lbb{stwrep}
Let $L$ be a Lie $H$-\psalg, and let
$\Pi$ be any finite-dimensional $\dd$-module.
We consider $\Pi$ as an $L$-module equipped with the trivial action
of $L$ and with the action of $H=\ue(\dd)$ induced from
the action of $\dd$. In particular, $\kk$ has the trivial action
of both $L$ and~$H$.

\begin{lemma}[\cite{BDK}]\lbb{lactchom}
Let $L$ be a Lie \psalg, and let $V,W$ be finite $L$-modules.
Then the formula $(a\in L$, $v\in V$, $\phi\in\Chom(V,W))$
\begin{equation}\lbb{actchom}
(a*\phi)*v = a*(\phi*v) - ((\si\tt\id)\tt_H\id) \, (\phi*(a*v))
\end{equation}
provides $\Chom(V,W)$ with the structure of an $L$-module.
\end{lemma}

Note that if $\be\colon V'\to V$ and $\ga\colon W\to W'$
are homomorphisms of $L$-modules, the map \eqref{chalbe1}
is a homomorphism of $L$-modules. 

\begin{definition}\lbb{dutw}
{\rm(i)}
For any finite $L$-module $V$, the $L$-module
$D(V) = \Chom(V,\kk)$ is called the {\em dual\/} of $V$.
If $\be\colon V'\to V$ is a homomorphism of $L$-modules,
we define a homomorphism $D(\be) \colon D(V) \to D(V')$
as $D(\be) = \Chom(\be,\id)$
{\rm(}see \reref{rchomfun}{\rm)}.
Then $D$ is a contravariant functor from the
category of finite $L$-modules to itself.

{\rm(ii)}
For any finite $L$-module $V$ and any finite-dimensional $\dd$-module $\Pi$,
the $L$-module $T_\Pi(V) = \Chom(D(V),\Pi)$
is called the {\em twisting\/} of $V$ by $\Pi$.
If $\be\colon V\to V'$ is a homomorphism of $L$-modules,
we define a homomorphism $T_\Pi(\be) \colon T_\Pi(V) \to T_\Pi(V')$
as $T_\Pi(\be) = \Chom(D(\be),\id)$.
Then $T_\Pi$ is a covariant functor from the
category of finite $L$-modules to itself.
\end{definition}

Now let $V$ be a free $H$-module of finite rank,
$V=H\tt V_0$, where $H$ acts by left multiplication on the first factor
and $\dim V_0 < \infty$.
Then for any $H$-module $W$ we can identify $\Chom(V,W)$ with
$H\tt (W \tt V_0^*)$, where $H$ acts on the first factor.
Explicitly, by \leref{lhhh}(ii),
for any fixed $v\in V_0$, we can write
\begin{equation}\lbb{phiv1}
\phi(1\tt v) = \sum (h_i \tt 1) \tt_H w_i \,,
\end{equation}
where $h_i \in H$, $w_i \in W$. Then $\phi$ corresponds to
the $\kk$-linear map $V_0 \to H\tt W$,
$v \mapsto \sum h_i \tt w_i$.

In particular, we have isomorphisms
$D(V) \simeq H\tt V_0^*$ and $T_\Pi(V) \simeq H\tt (\Pi\tt V_0)$
as $H$-modules. Now we will describe the action of $L$ on them.

\begin{proposition}\lbb{ptwrep}
Let $V=H\tt V_0$ be a finite $L$-module, which is free as an $H$-module.
Let\/ $\{ v_i \}$ be a\/ $\kk$-basis of\/ $V_0$, and let\/ $\{ \psi_i \}$
be the dual basis of\/ $V_0^*$, so that\/ $\psi_i(v_j) = \de_{ij}$.
For a fixed $a\in L$, write
\begin{equation}\lbb{twrep1}
a*(1\tt v_i) = \sum_j\, (f_{ij} \tt g_{ij}) \tt_H (1 \tt v_j)
\end{equation}
where $f_{ij}, g_{ij} \in H$.
%
Then the action of\/ $L$ on\/ $D(V) \simeq H\tt V_0^*$ is given by
\begin{equation}\lbb{twrep2}
a*(1\tt\psi_k) =
-\sum_j\, \bigl( f_{jk} {g_{jk}}_{(-1)} \tt {g_{jk}}_{(-2)} \bigr)
\tt_H \bigl( 1 \tt \psi_j \bigr) \,.
\end{equation}
%
The action of\/ $L$ on\/ $T_\Pi(V) \simeq H \tt (\Pi\tt V_0)$ is given by
\begin{equation}\lbb{twrep3}
a*(1\tt u\tt v_i) =
\sum_j\, \bigl( f_{ij} \tt {g_{ij}}_{(1)} \bigr)
\tt_H \bigl( 1 \tt {g_{ij}}_{(-2)} u \tt v_j \bigr) \,.
\end{equation}
\end{proposition}

Both \eqref{twrep2} and \eqref{twrep3} can be easily derived
from the following lemma.

\begin{lemma}\lbb{ltwrep}
Under the assumptions of \prref{ptwrep}, the action of\/ $L$ on\/
$\Chom(V,\Pi) \simeq H \tt (\Pi\tt V_0^*)$ is given by
\begin{equation}\lbb{twrep5}
a*(1\tt u\tt\psi_k) =
-\sum_j\, \bigl( f_{jk} {g_{jk}}_{(-1)} \tt {g_{jk}}_{(-2)} \bigr)
\tt_H \bigl( 1 \tt {g_{jk}}_{(3)} u \tt \psi_j \bigr) \,.
\end{equation}
\end{lemma}
\begin{proof}
First, note that by \eqref{phiv1}, we have
\begin{equation}\lbb{twrep4}
(1\tt u\tt\psi_k)*(1\tt v_i) = (1\tt 1) \tt_H \de_{ki} u \,.
\end{equation}
We will compute $(a*(1\tt u\tt\psi_k))*(1\tt v_i)$ using \eqref{actchom}.
The first term in the right-hand side of \eqref{actchom} vanishes because
the action of $L$ on $\Pi$ is trivial. By \eqref{psprod2}, \eqref{twrep1}
and \eqref{twrep4}, the second term is equal to
\begin{align*}
- ((\si\tt\id)&\tt_H\id) \, \bigl( (1\tt u\tt\psi_k)* (a*(1\tt v_i)) \bigr)
\\
&= - \sum_j\, (\si\tt\id) \,
(1 \tt f_{ij} \tt g_{ij}) \tt_H \de_{kj} u
\\
&= - (f_{ik} \tt 1 \tt g_{ik}) \tt_H u
\\
&= - (f_{ik} {g_{ik}}_{(-1)} \tt {g_{ik}}_{(-2)} \tt 1)
\tt_H {g_{jk}}_{(3)} u \,,
\end{align*}
where we used \eqref{cou2} in the last equality.
We will obtain the same result
if we apply the right-hand side of \eqref{twrep5} to $1\tt v_i$
and use \eqref{psprod1} and~\eqref{twrep4}.
\end{proof}
\begin{example}\lbb{ewdpi}
Consider $H$ as a $\Wd$-module via \eqref{wdac*}.
Then $T_\Pi(H) = H\tt\Pi$ with the following action of $\Wd$:
\begin{equation}\lbb{wdpi}
(1\tt a)*(1\tt u) = (1\tt 1)\tt_H (1\tt au)-(1\tt a)\tt_H (1\tt u)
\end{equation}
for $a\in\dd$, $u\in\Pi$.
\end{example}
\begin{remark}\lbb{rwdpi}
There is an embedding of Lie \psalgs\
\begin{equation}\lbb{wdcurd}
\Wd\injto\Cur\dd\rtimes\Wd \,, \quad
1\tt a \mapsto 1\tt a + 1\tt a \,,
\end{equation}
where the first summand is in $\Cur\dd=H\tt\dd$,
and the second one is in $\Wd=H\tt\dd$.
By \reref{rcurgtogc}, the representation of $\dd$ on $\Pi$ gives rise to
a homomorphism 
$\Cur\dd\rtimes\Wd\to\gc(H\tt\Pi)$.
Composing \eqref{curgtogc2} with \eqref{wdcurd},
we obtain a homomorphism $\Wd\to\gc(H\tt\Pi)$, which corresponds to the
$\Wd$-module $T_\Pi(H)$ from \exref{ewdpi} (see \reref{rcend}).
\end{remark}

Next, we will describe explicitly the homomorphisms $D(\be)$
and $T_\Pi(\be)$ from \deref{dutw}.

\begin{proposition}\lbb{ptwrep2}
Let $V=H\tt V_0$ and\/ $V'=H\tt V'_0$ be finite free $H$-modules.
Let\/ $\{ v_i \}$ $($respectively $\{ v'_i \})$ be a\/ $\kk$-basis of\/
$V_0$ $($respectively $V'_0)$,
and let\/ $\{ \psi_i \}$ $($respectively $\{ \psi'_i \})$
be the dual basis of\/ $V_0^*$ $($respectively $(V'_0)^*)$.
For a homomorphism of\/ $H$-modules $\be\colon V\to V'$, write
\begin{align}
\lbb{twrep6}
\be(1\tt v_i) &= \sum_j\, h_{ij} \tt v'_j
\intertext{where $h_{ij} \in H$. Then we have{\rm:}}
\lbb{twrep7}
D(\be)(1\tt\psi'_k) &= \sum_j\, S(h_{jk}) \tt \psi_j
\intertext{and}
\lbb{twrep8}
T_\Pi(\be)(1\tt u\tt v_i) &=
\sum_j\, {h_{ij}}_{(1)} \tt {h_{ij}}_{(-2)} u \tt v'_j  \,.
\end{align}
\end{proposition}

By linearity, \prref{ptwrep2} follows from the following special case,
which we formulate as a lemma for future reference.

\begin{lemma}\lbb{ltwrep2}
Let $V=H\tt V_0$ and\/ $V'=H\tt V'_0$ be finite free $H$-modules. For fixed
$h \in H$, $B \in \Hom_\kk(V_0,V'_0)$, consider the
homomorphism of\/ $H$-modules $\be\colon V\to V'$ given by
\begin{align}
\lbb{twrep9}
\be(1\tt v) &= h \tt Bv
\,, \qquad v \in V_0 \,.
\intertext{Then we have{\rm:}}
\lbb{twrep10}
D(\be)(1\tt\psi') &= S(h) \tt (\psi' \circ B)
\,, \qquad \psi' \in (V'_0)^* = \Hom_\kk(V'_0,\kk)
\intertext{and}
\lbb{twrep11}
T_\Pi(\be)(1\tt u\tt v) &= h_{(1)} \tt h_{(-2)} u \tt Bv
\,, \qquad u \in \Pi \,, \; v \in V_0 \,.
\end{align}
\end{lemma}
\begin{proof}
The proof is straightforward from definition, and it is left to the reader.
\end{proof}



\subsection{Tensor Modules for $\Wd$}\lbb{stenw}
The adjoint representation of $\Wd=H\tt\dd$ gives rise to
the following homomorphism of Lie \psalgs\ $\Wd\to\gc(H\tt\dd)$
(see \eqref{wdbr*} and \reref{rcend}):
\begin{equation}\lbb{wdgcd}
1\tt a \mapsto 1 \tt 1 \tt \ad a - 1 \tt a \tt \Id + \ep_a \,,
\end{equation}
where the pseudolinear map $\ep_a$ is given by
\begin{equation}\lbb{epab}
\ep_a(g \tt b) = (b \tt g) \tt_H (1 \tt a) \,,
\qquad g \in H \,, \; b \in\dd \,.
\end{equation}
In \eqref{wdgcd} we have identified $\gc(H\tt\dd)$ with $H\tt H\tt\End\dd$\,;
in this identification
\begin{equation}\lbb{wdepa}
\ep_{\d_i} = \sum_{j=1}^N\, \d_j \tt 1 \tt e_i^j \,,
\end{equation}
where $\{\d_i\}_{i=1,\dots,N}$ is a basis of $\dd$
and $e_i^j(\d_k) = \de^j_k \d_i$.

\begin{lemma}\lbb{lwdincurgld}
The map
\begin{equation}\lbb{wdgld}
1\tt\d_i \mapsto
\Bigl( 1\tt\ad\d_i + \sum_{j=1}^N\, \d_j \tt e_i^j \Bigr) + 1\tt\d_i
\end{equation}
is an embedding of Lie \psalgs\ $\Wd\injto(\Cur\gld)\rtimes\Wd$.
\end{lemma}
\begin{proof}
The image of $\Wd$ under the embedding \eqref{wdgcd}
is contained in 
\begin{equation*}
H\tt\kk\tt\End\dd+H\tt\dd\tt\id \subset H\tt H\tt\End\dd \,,
\end{equation*}
which is isomorphic to $(\Cur\gld)\rtimes\Wd$
by \reref{rcurgtogc}.
\end{proof}

By \reref{rcurgtogc}, for any finite-dimensional $\gld$-module $V_0$,
we have a homomorphism of Lie \psalgs\
$(\Cur\gld)\rtimes\Wd\to\gc V$, where $V=H\tt V_0$.
After composing it with the embedding from \leref{lwdincurgld}, 
we obtain a homomorphism $\Wd\to\gc V$, i.e.,
a representation of $\Wd$ on $V$.
Explicitly, the action of $\Wd$ on $V$ is given by:
\begin{equation}\lbb{wdgcd2}
\begin{split}
(1\tt \d_i)*(1\tt v) &= (1 \tt 1) \tt_H (1 \tt (\ad\d_i)v)
+ \sum_{j=1}^N\, (\d_j \tt 1) \tt_H (1 \tt e_i^jv)
\\
&- (1 \tt \d_i) \tt_H (1 \tt v) \,.
\end{split}
\end{equation}

Now let $\Pi$ be any finite-dimensional $\dd$-module.
The twisting of $V$ by $\Pi$ is $T_\Pi(V) = H\tt (\Pi\tt V_0)$
with the following action of $\Wd$ (see \prref{ptwrep}):
\begin{equation}\lbb{wdgcd3}
\begin{split}
(1\tt \d_i)*(1\tt w) &= (1 \tt 1) \tt_H (1 \tt (\ad\d_i)w)
+ \sum_{j=1}^N\, (\d_j \tt 1) \tt_H (1 \tt e_i^j w)
\\
&- (1 \tt \d_i) \tt_H (1 \tt w)
+ (1 \tt 1) \tt_H (1 \tt \d_i w)
\end{split}
\end{equation}
for $w\in\Pi\tt V_0$, where $\dd$ acts on the factor $\Pi$
and $\gld$ acts on $V_0$.

\begin{definition}\lbb{dboxt}
Let $\fg_1$ and $\fg_2$ be Lie algebras, and let $U_i$ be a
$\fg_i$-module ($i=1,2$). Then we will denote by $U_1 \bt U_2$
the $(\fg_1\oplus\fg_2)$-module $U_1 \tt U_2$, where $\fg_1$
acts on the first factor and $\fg_2$ acts on the second one.
\end{definition}

The above formulas \eqref{wdgcd2}, \eqref{wdgcd3}
motivate the introduction of an important class of $\Wd$-modules.

\begin{definition}\lbb{dtenmodw}
{\rm(i)}
Let $W_0$ be a finite-dimensional $(\dd\oplus\gld)$-module.
The $\Wd$-module $H\tt W_0$, with the action of $\Wd$ given by \eqref{wdgcd3}
for $w\in W_0$, is called a {\em tensor module\/} and is
denoted as $\T(W_0)$.

{\rm(ii)}
Let $W_0 = \Pi \bt V_0$,
where $\Pi$ is a finite-dimensional $\dd$-module
and $V_0$ is a finite-dimensional $\gld$-module.
Then the tensor module $\T(W_0)$ will also be denoted as $\T(\Pi,V_0)$.

{\rm(iii)}
Occasionally, we will denote $\T(\Pi,V_0)$ also by $\T(\Pi, V_0, c)$,
where $V_0$ is viewed as a module over $\sld\subset\gld$,
and $c\in \kk$ denotes the scalar action of $\Id\in\gld$ on~$V_0$.

\end{definition}
\begin{remark}\lbb{rtenmodw}
By definition, we have $\T(\Pi,V_0) = T_\Pi(\T(\kk,V_0))$.
\end{remark}
\begin{remark}\lbb{rwdtenm}
Combining the embeddings \eqref{wdcurd} and \eqref{wdgld},
we get an embedding of Lie \psalgs\
$\Wd\injto\Cur(\dd\oplus\gld)\rtimes\Wd$,
\begin{equation}\lbb{wdgld2}
1\tt\d_i \mapsto
\Bigl( 1\tt\d_i + 1\tt\ad\d_i + \sum_{j=1}^N\, \d_j \tt e_i^j \Bigr)
+ 1\tt\d_i \,.
\end{equation}
Given a $(\dd\oplus\gld)$-module $W_0$, the $\Wd$-module obtained from it
by \reref{rcurgtogc} is exactly the tensor module $\T(W_0) = H\tt W_0$
corresponding to $W_0$.
\end{remark}

\section{Tensor Modules of de Rham Type}\lbb{sdrm}
%
Throughout this section, $\dd$ will be an $N$-dimensional Lie algebra.
We fix a basis $\{\d_i\}_{i=1,\dots,N}$ of $\dd$ with structure
constants $c_{ij}^k$: $[\d_i,\d_j]=\sum\,c_{ij}^k\d_k$.
Define elements $e^i_j\in\gld$ by $e^i_j(\d_k) = \de^i_k \d_j$.

\subsection{Forms with Constant Coefficients}\lbb{sfcc}
The material in this subsection is completely standard;
our purpose is just to fix the notation.
For $0\le n\le N$, let
\begin{equation}\lbb{omn}
\Om^n = {\textstyle\bigwedge}^n\dd^* \,, \quad
\Om = {\textstyle\bigwedge}^\bullet\dd^*
= \bigoplus_{n=0}^N \Om^n \,.
\end{equation}
Set $\Om^n = \{0\}$ if $n<0$ or $n>N$.
We will think of the elements of $\Om^n$ as skew-symmetric
{\em $n$-forms}, i.e., linear maps from $\bigwedge^n \dd$ to $\kk$.

Consider the cohomology complex of $\dd$
with trivial coefficients,
\begin{equation}\lbb{domcc}
0 \to \Om^0 \xrightarrow{\diz} \Om^1
\xrightarrow{\diz} \cdots \xrightarrow{\diz} \Om^{N} ,
\end{equation}
where the differential $\diz$
is given by the formula ($\al\in\Om^n$, $a_i\in\dd$):
\begin{equation}\lbb{d0al}
\begin{split}
(\diz&\al)(a_1 \wedge \dots \wedge a_{n+1})
\\
&= \sum_{i<j}\, (-1)^{i+j} \al([a_i,a_j] \wedge a_1 \wedge \dots \wedge
\what a_i \wedge \dots \wedge \what a_j \wedge \dots \wedge a_{n+1})
\end{split}
\end{equation}
if $n\ge1$, and $\diz\al = 0$ for $\al\in\Om^0=\kk$.
Here a hat over $a_i$ means that the term $a_i$ is omitted
in the wedge product.

For $a\in\dd$, define operators $\io_a\colon\Om^n\to\Om^{n-1}$ by
\begin{equation}\lbb{ioal}
(\io_a\al)(a_1 \wedge \dots \wedge a_{n-1})
=\al(a \wedge a_1 \wedge \dots \wedge a_{n-1}) \,,
\qquad\quad a_i\in\dd \,.
\end{equation}
For $A\in\gld$, denote by $A\cdot$ its action on $\Om$\,; explicitly,
\begin{equation}\lbb{acdotal}
(A\cdot\al)(a_1 \wedge \dots \wedge a_n)
= \sum_{i=1}^n\, (-1)^i \al(A a_i \wedge a_1 \wedge \dots \wedge \what a_i
\wedge \dots \wedge a_n) \,.
\end{equation}
Then we have the following Cartan formula for the coadjoint action of $\dd$\,:
\begin{equation}\lbb{cartan}
(\ad a)\cdot = \diz\io_a + \io_a\diz \,.
\end{equation}
This, together with $\diz^2=0$,
implies that $(\ad a)\cdot$ commutes with $\diz$.

\subsection{Pseudo de Rham Complex}\lbb{sdrmw}
Following \cite[Section~8.3]{BDK},
we define the spaces of {\em pseudoforms\/}
$\Om^n(\dd) = H\tt\Om^n$ and
$\Om(\dd) = H\tt\Om = \textstyle\bigoplus_{n=0}^N \Om^n(\dd)$.
They are considered as $H$-modules, where
$H$ acts on the first factor by left multiplication.
We can identify
$\Om^n(\dd)$ with the space of linear maps from $\bigwedge^n \dd$ to $H$,
and $H^{\tt2}\tt_H \Om^n(\dd)$ with
$\Hom(\bigwedge^n \dd, H^{\tt2})$.
Note that $\Om^n(\dd) = \{0\}$ if $n<0$ or $n>N$.

Let us consider $H=\ue(\dd)$ as a left $\dd$-module with respect
to the action $a\cdot h = -ha$, where
$ha$ is the product of $a\in\dd\subset H$ and $h\in H$ in $H$.
Consider the cohomology complex of $\dd$
with coefficients in $H$:
\begin{equation}\lbb{domd}
0 \to \Om^0(\dd) \xrightarrow{\di} \Om^1(\dd)
\xrightarrow{\di} \cdots \xrightarrow{\di} \Om^{N}(\dd) \,.
\end{equation}
Explicitly, the differential $\di$ is given by the formula
($\al\in\Om^n(\dd)$, $a_i\in\dd$):
\begin{equation}\lbb{dw}
\begin{split}
&\begin{split}
(\di \al)(a_1 & \wedge \dots \wedge a_{n+1})
\\
&= \sum_{i<j}\, (-1)^{i+j} \al([a_i,a_j] \wedge a_1 \wedge \dots \wedge
\what a_i \wedge \dots \wedge \what a_j \wedge \dots \wedge a_{n+1})
\\
&+ \sum_i\, (-1)^i \al(a_1 \wedge \dots \wedge \what a_i \wedge \dots
\wedge a_{n+1}) \, a_i
\qquad\text{if}\;\; n\ge1,
\end{split}
\\
&(\di \al)(a_1) = -\al a_1 \qquad\text{if}\;\; \al\in\Om^0(\dd)=H \,,
\end{split}
\end{equation}
where a hat over $a_i$ means that the term $a_i$ is omitted.
Notice that $\di$ is $H$-linear.

\begin{proposition}[\cite{BDK}]\lbb{pdomd1}
The $n$-th cohomology of the complex $(\Om(\dd), \di)$ is
trivial for $n\ne N=\dim\dd$ and\/ $1$-dimensional for $n=N$.
In particular, the sequence \eqref{domd} is exact.
\end{proposition}
\begin{proof}
By Poincar\'e duality 
$\coh^n(\dd,\ue(\dd)) \simeq \coh^{N-n}(\dd,\ue(\dd)^*)$.
But $\coh^n(\dd,\ue(\dd)^*) \simeq \coh_n(\dd,\ue(\dd))^*$
is trivial for $n>0$ and $1$-dimensional for $n=0$;
see, e.g., \cite{F}.
\end{proof}
\begin{definition}\lbb{dpsdr}
The sequence \eqref{domd} is called the {\em pseudo de Rham complex}.
\end{definition}

The following lemma provides another formula for the differential
\eqref{dw}, which will be useful later.

\begin{lemma}\lbb{ldw}
For every $\al\in\Om^n$, $n\ge0$, and $i=1,\dots,N$, consider the element
$1 \tt \io_{\d_i}\al \in \Om^{n-1}(\dd)$.
Then we have{\rm:}
\begin{equation}\lbb{dw2}
\di (1 \tt \io_{\d_i}\al)
= \sum_{k=1}^N \, \d_k \tt e^k_i \al
- \sum_{\substack{j,k,l=1 \\ k<l}}^N
1 \tt c^j_{kl} e^k_i e^l_j \al
- \sum_{\substack{k,l=1 \\ k<l}}^N
1 \tt c^k_{kl} e^l_i \al \,,
\end{equation}
where the action of\/ $\gld$ on\/ $\Om^n$ is given by \eqref{acdotal}.
\end{lemma}
\begin{proof}
For $n=0$, both sides of \eqref{dw2} are trivial; so we can assume that
$n\ge1$. Denote the three terms in the right-hand side of \eqref{dw2}
by $\be_1$, $\be_2$, $\be_3$.
Using \eqref{acdotal},
we compute for $1\le i_1 < \dotsm < i_n \le N$:
\begin{align*}
\be_1(\d_{i_1} &\wedge\cdots\wedge \d_{i_n})
= \sum_{k=1}^N \sum_{r=1}^n \, \d_k \tt
(-1)^r \de^k_{i_r} \,
\al(\d_i \wedge \d_{i_1} \wedge\cdots\wedge
\what \d_{i_r} \wedge\cdots\wedge \d_{i_n})
\intertext{and}
\be_2(\d_{i_1} &\wedge\cdots\wedge \d_{i_n})
= \sum_{\substack{j,k,l=1 \\ k<l}}^N \sum_{\substack{r,s=1 \\ r<s }}^n
1 \tt (-1)^{r+s+1} \de^k_{i_r} \de^l_{i_s} c^j_{kl} \,
\al(\d_i \wedge \d_j \wedge \d_{i_1} \wedge\cdots
\\
&\qquad\qquad\qquad\qquad\qquad\qquad\qquad\qquad
\wedge \what \d_{i_r} \wedge\cdots\wedge
\what \d_{i_s} \wedge\cdots\wedge \d_{i_n})
\\
&+ \sum_{\substack{j,k,l=1 \\ k<l}}^N \sum_{s=1}^n \,
1 \tt (-1)^{s+1} \de^k_j \de^l_{i_s} c^j_{kl} \,
\al(\d_i \wedge \d_{i_1} \wedge\cdots\wedge
\what \d_{i_s} \wedge\cdots\wedge \d_{i_n}) \,.
\intertext{Similarly,}
\be_3(\d_{i_1} &\wedge\cdots\wedge \d_{i_n})
= \sum_{\substack{k,l=1 \\ k<l}}^N \sum_{s=1}^n \,
1 \tt (-1)^s \de^l_{i_s} c^k_{kl} \,
\al(\d_i \wedge \d_{i_1} \wedge\cdots\wedge
\what \d_{i_s} \wedge\cdots\wedge \d_{i_n}) \,.
\end{align*}
These formulas, together with \eqref{dw}, \eqref{ioal}
and the equation
\begin{equation*}
\sum_{\substack{j,k,l=1 \\ k<l}}^N \de^k_{i_r} \de^l_{i_s} c^j_{kl} \, \d_j
= [\d_{i_r}, \d_{i_s}] \,, \qquad
r<s \,, \;\; i_r < i_s \,,
\end{equation*}
imply that $\di (1 \tt \io_{\d_i}\al) = \be_1-\be_2-\be_3$.
\end{proof}

Next, we introduce $H$-bilinear maps
\begin{align}
\lbb{*ion}
*_\io &\colon \Wd \tt \Om^n(\dd) \to H^{\tt2} \tt_H \Om^{n-1}(\dd) \,,
\\
\lbb{*n}
* &\colon \Wd \tt \Om^n(\dd) \to H^{\tt2} \tt_H \Om^n(\dd) \,,
\end{align}
by the formulas:
\begin{equation}\lbb{*iota}
(f\tt a)*_\io (g\tt\al) = (f\tt g) \tt_H \io_a\al \,,
\end{equation}
\begin{equation}\lbb{*cartan}
w*\ga = ((\id\tt\id)\tt_H \di)(w*_\io\ga) + w*_\io (\di\ga) \,,
\end{equation}
for $w=f\tt a\in\Wd$, $\ga=g\tt\al\in\Om^n(\dd)$.
Eq.\ \eqref{*cartan} is an analog of Cartan's formula \eqref{cartan}.
Explicitly, we have (see \cite[Eq.~(8.7)]{BDK}):
\begin{equation}\lbb{fa*om}
\begin{split}
(w&*\ga)(a_1 \wedge \dots \wedge a_n)
= -(f \tt g a) \, \al(a_1 \wedge \dots \wedge a_n)
\\
&+ \sum_{i=1}^n\, (-1)^i
(f a_i \tt g) \, \al(a \wedge a_1 \wedge \dots \wedge \what a_i
\wedge \dots \wedge a_n)
\\
&+ \sum_{i=1}^n\, (-1)^i
(f \tt g) \, \al([a,a_i] \wedge a_1 \wedge \dots \wedge \what a_i
\wedge \dots \wedge a_n) \in H^{\tt2}
\end{split}
\end{equation}
for $n\ge1$, and $w*\ga = -f\tt ga$ for $\ga=g\in\Om^0(\dd)=H$.
Note that the latter coincides with the action \eqref{wdac*}
of $\Wd$ on $H$.

\begin{theorem}\lbb{twdomd}
The maps \eqref{*n} provide each\/ $\Om^n(\dd)$
with a structure of a tensor\/ $\Wd$-module corresponding to
the\/ $(\dd\oplus\gld)$-module\/ $\kk\bt\Om^n$,
i.e., $\Om^n(\dd) = \T(\kk, \Om^n)$.
The differential\/ $\di\colon \Om^n(\dd) \to \Om^{n+1}(\dd)$
is a homomorphism of\/ $\Wd$-modules.
\end{theorem}
\begin{proof}
Comparing \eqref{fa*om} with \eqref{acdotal}, we obtain
for $\al\in\Om^n$:
\begin{align*}
(1\tt \d_k)&*(1\tt\al)
= - (1 \tt \d_k) \tt_H (1 \tt \al)
\\
&+ \sum_{j=1}^N\, (\d_j \tt 1) \tt_H (1 \tt e_k^j \al)
+ (1 \tt 1) \tt_H (1 \tt (\ad\d_k)\al) \,.
\end{align*}
But this is exactly \eqref{wdgcd2}; hence, $\Om^n(\dd) = \T(\kk, \Om^n)$.

To prove that $\di$ is a homomorphism, we have to check that
it satisfies \eqref{psprod6}.
This follows from \eqref{*cartan} and $\di^2=0$.
Indeed, replacing $\ga$ with $\di\ga$ in \eqref{*cartan}, we get
\begin{align*}
w*(\di\ga) &= ((\id\tt\id)\tt_H \di)(w*_\io\di\ga) \,,
\intertext{while applying $((\id\tt\id)\tt_H \di)$ to both sides of
\eqref{*cartan} gives}
((\id\tt\id)\tt_H \di)(w*\ga) &= ((\id\tt\id)\tt_H \di)(w*_\io\di\ga) \,.
\end{align*}
This completes the proof.
\end{proof}

\subsection{Twisting of the Pseudo de Rham Complex}\lbb{stwdrmw}
As before, let $\Pi$ be a finite-dimensional $\dd$-module, which we consider
as an $H$-module. We will apply the twisting functor $T_\Pi$
(see \deref{dutw}(ii)) to the pseudo de Rham complex \eqref{domd}.
Note that, by \thref{twdomd} and \reref{rtenmodw}, we have
$T_\Pi(\Om^n(\dd)) = \T(\Pi,\Om^n)$.
We obtain a complex of $\Wd$-modules
\begin{equation}\lbb{domd2}
0 \to \T(\Pi,\Om^0) \xrightarrow{\di_\Pi} \T(\Pi,\Om^1)
\xrightarrow{\di_\Pi} \cdots \xrightarrow{\di_\Pi} \T(\Pi,\Om^N) \,,
\qquad \di_\Pi \equiv T_\Pi(\di) \,,
\end{equation}
which we call the {\em $\Pi$-twisted pseudo de Rham complex}.

It follows from \eqref{dw} and \prref{ptwrep2} that the complex
\eqref{domd2} coincides with the cohomology complex of $\dd$
with coefficients in $H\tt\Pi$ considered with the action
\begin{equation}\lbb{dhpi}
a\cdot (h \tt u) = -ha \tt u + h \tt au \,,
\qquad a\in\dd \,, \; h\in H \,, \; u\in\Pi \,.
\end{equation}

\begin{lemma}\lbb{ldomd2}
The\/ $\dd$-module $H\tt\Pi$, equipped with the action \eqref{dhpi},
is isomorphic to $H\tt\Pi$ with\/ $\dd$ acting only on $H$ via
\begin{equation*}
a(h \tt u) = -ha \tt u \,,
\qquad a\in\dd \,, \; h\in H \,, \; u\in\Pi \,.
\end{equation*}
In other words, $H\tt\Pi$ is isomorphic to a direct sum of\/ $\dim\Pi$
copies of the\/ $\dd$-module~$H$.
\end{lemma}
\begin{proof}
Consider the linear map
\begin{equation*}
F\colon H\tt\Pi \to H\tt\Pi \,,
\qquad h\tt u \mapsto h_{(1)} \tt h_{(-2)} u \,.
\end{equation*}
{}From \eqref{cou2} it is easy to see that $F$ is a linear isomorphism and
\begin{equation*}
F^{-1}(h\tt u) = h_{(1)} \tt h_{(2)} u
\end{equation*}
(see \cite[Section~2.3]{BDK} for a similar argument).
Using \eqref{deprod} and \eqref{deasa}, we compute
\begin{align*}
F(-ha \tt u) &= -(ha)_{(1)} \tt (ha)_{(-2)} u
\\
&= -h_{(1)} a \tt h_{(-2)} u + h_{(1)} \tt a h_{(-2)} u
=a\cdot F(h \tt u) \,.
\end{align*}
This shows that $F$ is an isomorphism of the corresponding
$\dd$-modules.
\end{proof}

Now \prref{pdomd1} and \leref{ldomd2} immediately imply:

\begin{proposition}\lbb{pdomd2}
The sequence \eqref{domd2} is exact.
The image of\/ $\di_\Pi$ in\/ $\T(\Pi,\Om^N)$ has codimension~$\dim\Pi$.
\end{proposition}

Finally, let us give a formula for the differential $\di_\Pi$,
which is similar to \eqref{dw2}.

\begin{lemma}\lbb{ldw2}
For every $\al\in\Om^n$, $n\ge0$, $u\in\Pi$, and $i=1,\dots,N$, we have{\rm:}
\begin{equation}\lbb{dw3}
\begin{split}
\di_\Pi (1 \tt u \tt \io_{\d_i}\al)
&= \sum_{k=1}^N \, \d_k \tt u \tt e^k_i \al
- \sum_{k=1}^N \, 1 \tt \d_k u \tt e^k_i \al
\\
&- \sum_{\substack{j,k,l=1 \\ k<l}}^N
1 \tt u \tt c^j_{kl} e^k_i e^l_j \al
- \sum_{\substack{k,l=1 \\ k<l}}^N
1 \tt u \tt c^k_{kl} e^l_i \al \,,
\end{split}
\end{equation}
where the action of\/ $\gld$ on\/ $\Om^n$ is given by \eqref{acdotal}.
\end{lemma}
\begin{proof}
For a fixed $i$, extend by $H$-linearity the map $\io_{\d_i}$ to a map
from $\Om^n(\dd)$ to $\Om^{n-1}(\dd)$.
Note that, by \prref{ptwrep2} (or \leref{ltwrep2}), we have
\begin{equation*}
(T_\Pi(\io_{\d_i}))(1\tt u \tt \al) = 1 \tt u \tt \io_{\d_i} \al \,.
\end{equation*}
Consider the $H$-linear map
$\di\circ\io_{\d_i} \colon \Om^n(\dd) \to \Om^n(\dd)$.
This map sends $1\tt\al$ to the right-hand side of \eqref{dw2}.
Hence, by \prref{ptwrep2}, $(T_\Pi(\di\circ\io_{\d_i}))(1\tt u\tt\al)$
is given by the right-hand side of \eqref{dw3}.
On the other hand, we have
$T_\Pi(\di\circ\io_{\d_i}) = T_\Pi(\di) \circ T_\Pi(\io_{\d_i})$,
which gives $\di_\Pi (1 \tt u \tt \io_{\d_i}\al)$
when applied to $1\tt u \tt \al$.
\end{proof}

\section{Classification of Irreducible Finite $\Wd$-Modules}\lbb{swmod}

In this section we provide a complete classification of all irreducible
finite $\Wd$-modules. Our main result is \thref{twmod}.

\subsection{Singular Vectors and Tensor Modules}\lbb{swten}
Recall that the annihilation algebra $\W$ of $\Wd$ has a
decreasing filtration $\{\W_p\}_{p\ge-1}$ given by \eqref{wp}.
For a $\W$-module $V$, we denote by $\ker_p V$ the set of all
$v\in V$ that are killed by $\W_p$.
Then a $\W$-module $V$ is conformal iff $V=\bigcup \ker_p V$.
Recall also that the extended annihilation algebra is
defined as $\ti\W = \dd\ltimes\W$, where $\dd$ acts on $\W$
by \eqref{dactw}.
By \prref{preplal2}, any $\Wd$-module
has a natural structure of a conformal $\ti\W$-module and
vice versa.

For every $p\ge 0$, the normalizer of $\W_p$ in $\ti\W$
is equal to $\N_\W$ (see \deref{dnw} and \prref{normalizer}).
Therefore, each $\ker_p V$ is an $\N_\W$-module.
In fact, $\ker_p V$ is a module over the finite-dimensional
Lie algebra $\N_\W / \W_p = \ti\dd \oplus (\W_0 / \W_p)$.
The Lie algebra $\N_\W / \W_1$ is isomorphic to
the direct sum of Lie algebras $\dd\oplus\gld$.

\begin{definition}\lbb{dwsing}
For a $\Wd$-module $V$, a {\em singular vector\/}
is an element $v\in V$ such that $\W_1\cdot v = 0$.
The space of singular vectors in $V$ will be denoted by $\sing V$.
We will denote by
$\rho_{\sing} \colon \dd\oplus\gld \to \gl(\sing V)$
the representation obtained from the $\N_\W$-action on
$\sing V\equiv\ker_1 V$
via the isomorphism $\N_\W / \W_1 \simeq \dd\oplus\gld$.
\end{definition}

Recall that $\ker V \equiv \ker_{-1} V$ is the space of all
$v\in V$ such that $\W\cdot v = 0$. Then, obviously,
$\ker V \subset \sing V$.
Note also that $\ker V = \{0\}$ when $V$ is irreducible.

\begin{theorem}\lbb{twsing1}
For any nontrivial finite\/ $\Wd$-module $V$, we have\/ $\sing V \ne\{0\}$
and the space\/ $\sing V / \ker V$ is finite dimensional.
\end{theorem}
\begin{proof}
The second statement is a special case of \leref{lkey2}.
To show that $\sing V \ne\{0\}$, we can assume
without loss of generality that $\ker V = \{0\}$.
Since $V$ is a conformal $\W$-module, $\ker_p V\neq \{0\}$
for some $p\ge 1$. By \leref{lkey2}, the space $\ker_p V$
is finite dimensional. Let us choose a minimal $\N_\W$-submodule
$R$ of $\ker_p V$. Then $R$ is an irreducible $\N_\W$-module; hence,
by \prref{Nreps}, $\W_1$ acts trivially on $R$.
This means that $R\subset \sing V$.
\end{proof}
\begin{remark}\lbb{rwsing}
It follows from \eqref{wp} and \prref{preplal2} that
a vector $v\in V$ is singular if and only if
\begin{equation}\lbb{vsingf1}
(1\tt\d)*v \in (\fil^1 H \tt \kk) \tt_H V \,,
\qquad \d\in\dd \,,
\end{equation}
where $\fil^1 H = \kk\oplus\dd$. Similarly, by \leref{lhhh}(ii),
a vector $v\in V$ is singular if and only if
\begin{equation}\lbb{vsingf2}
(1\tt\d)*v \in (\kk \tt \fil^1 H) \tt_H V \,,
\qquad \d\in\dd \,.
\end{equation}
\end{remark}

As before, let $\{\d_i\}_{i=1,\dots,N}$ be a basis of $\dd$,
and let $x^i\in X$ be given by \eqref{xi2}.
We view $x^i$ as elements of $\dd^*$;
then $\{x^i\}$ is a basis of $\dd^*$ dual to the basis $\{\d_i\}$ of $\dd$.
Let $e_i^j \in\gld$ be given by $e_i^j\d_k=\de^j_k\,\d_i$;
in other words, $e_i^j$ corresponds to $\d_i\tt x^j$
under the isomorphism $\gld\simeq\dd\tt\dd^*$.

Note that, by \deref{dnw} and \coref{cwbra}, we have
\begin{equation}\lbb{rhosing}
\rho_{\sing}(\d) v = \ti\d \cdot v \,, \quad
\rho_{\sing}(e_i^j) v = -(x^j\tt\d_i) \cdot v \,,
\qquad \d\in\dd \,, \; v\in \sing V \,.
\end{equation}

\begin{lemma}\lbb{lwsing}
Let $V$ be a $\Wd$-module. Then for every singular vector
$v\in \sing V$, the action of\/ $\Wd$ on $v$ is given by
\begin{equation}\lbb{Wactionsing}
\begin{split}
(1 \tt \d_i) * v
= \sum_{j=1}^N  & \, (\d_j \tt 1) \tt_H \rho_{\sing}(e_i^j) v
- (1 \tt 1) \tt_H \d_i v
\\
+ & \, (1 \tt 1) \tt_H \rho_{\sing}(\d_i + \ad \d_i) v  \,.
\end{split}
\end{equation}
\end{lemma}
\begin{proof}
Since $\W_1 \cdot v =0$, it follows from \prref{preplal2} that
for $\d\in\dd$
\begin{equation*}
(1\tt\d) * v = (1 \tt 1) \tt_H (1\tt\d) \cdot v
- \sum_j\, (\d_j \tt 1) \tt_H (x^j\tt\d) \cdot v \,,
\end{equation*}
while \leref{lnw2} implies
\begin{equation*}
(\d + 1 \tt \d - \ti\d) \cdot v = \rho_{\sing}(\ad \d) v  \,.
\end{equation*}
Combining the above equations with \eqref{rhosing}
proves \eqref{Wactionsing}.
\end{proof}
\begin{corollary}\lbb{cwsing}
Let $V$ be a $\Wd$-module and let $R$ be a nontrivial
$(\dd\oplus\gld)$-submodule of\/ $\sing V$.
Denote by $HR$ the $H$-submodule of\/ $V$ generated by $R$.
Then $HR$ is a $\Wd$-submodule of\/ $V$.
In particular, if\/ $V$ is irreducible, then $V=HR$.
\end{corollary}
\begin{proof}
It follows from \eqref{Wactionsing} that
$\Wd * R \subset (H\tt H)\tt_H HR$. Then, by $H$-bilinearity,
$\Wd * HR \subset (H\tt H)\tt_H HR$, which means that
$HR$ is a $\Wd$-submodule of $V$.
\end{proof}

Let $R$ be a finite-dimensional $(\dd\oplus\gld)$-module,
with an action denoted as $\rho_R$. Let $V=H\tt R$ be the free $H$-module
generated by $R$, where $H$ acts by left multiplication on the first
factor. We define a pseudoproduct
\begin{equation}\lbb{wsing1}
\begin{split}
(1 \tt \d_i) * (1 \tt u)
= \sum_{j=1}^N  & \, (\d_j \tt 1) \tt_H (1 \tt \rho_R(e_i^j) u)
- (1 \tt 1) \tt_H (\d_i \tt u)
\\
+ & \, (1 \tt 1) \tt_H (1 \tt \rho_R(\d_i + \ad \d_i) u)  \,,
\qquad u\in R \,,
\end{split}
\end{equation}
and then extend it by $H$-bilinearity to a map
$* \colon \Wd\tt V \to (H\tt H)\tt_H V$.

\begin{lemma}\lbb{lwsing2}
Let\/ $R$ be a finite-dimensional $(\dd\oplus\gld)$-module with
an action $\rho_R$. Then formula \eqref{wsing1} defines a structure of a
$\Wd$-module on $V=H\tt R$. We have $\kk\tt R \subset\sing V$ and
\begin{equation}\lbb{wsing7}
\rho_{\sing}(A)(1\tt u) = 1\tt \rho_R(A)u \,,
\qquad A \in \dd\oplus\gld \,, \; u \in R \,.
\end{equation}
\end{lemma}
\begin{proof}
The fact that $V$ is a $\Wd$-module can be proved by a straightforward
computation, using \eqref{psrep} and \eqref{wdbr*}. Instead, we will show
that $V$ is a tensor module (see \deref{dtenmodw}).
Let us compare \eqref{wsing1} to \eqref{wdgcd3},
keeping in mind that, by definition,
\begin{equation*}
(1 \tt 1) \tt_H (\d_i \tt u)
= (\d_i \tt 1) \tt_H (1 \tt u)
+ (1 \tt \d_i) \tt_H (1 \tt u) \,.
\end{equation*}
We see that $V=H\tt R$ coincides with the tensor module $\T(R)$,
where $R$ is equipped with the following modified action of $\dd\oplus\gld$:
\begin{align}
\lbb{wsing2}
\d u &= \bigl( \rho_R(\d) + \tr(\ad\d) \bigr) u \,,
&&\d\in\dd \,, \; u\in R \,,
\\
\lbb{wsing3}
A u &= \bigl( \rho_R(A) - \tr A \bigr) u \,,
&&A\in\gld \,, \; u\in R \,.
\end{align}
The fact that $\kk\tt R \subset\sing V$ follows from \reref{rwsing},
and \eqref{wsing7} follows from comparing
\eqref{Wactionsing} with \eqref{wsing1}.
This completes the proof.
\end{proof}
\begin{definition}\lbb{dvmodw}
{\rm(i)}
Let $R$ be a finite-dimensional $(\dd\oplus\gld)$-module with
an action $\rho_R$. Then the $\Wd$-module $H\tt R$, with the action of
$\Wd$ given by \eqref{wsing1},
will be denoted as~$\V(R)$.

{\rm(ii)}
Let $R = \Pi\bt U$,
where $\Pi$ is a finite-dimensional $\dd$-module
and $U$ is a finite-dimensional $\gld$-module.
Then the module $\V(R)$ will also be denoted as~$\V(\Pi,U)$.

\end{definition}
\begin{remark}\lbb{rvrind}
If $R$ is a finite-dimensional $(\dd\oplus\gld)$-module, we can define
an action of $\N_\W$ on it by letting $\W_1$ act as zero. 
Then as a $\wti\W$-module, $\V(R)$ is isomorphic to the induced module
$\Ind^{\wti\W}_{\N_\W} R$. This follows from the fact that
$\wti\W = \dd \oplus \N_\W$ as a vector space (see \prref{normalizer}).
\end{remark}


For a Lie algebra $\fg$ and a trace form $\chi$ on $\fg$,
we denote by $\kk_\chi$ the $1$-dimensional $\fg$-module
such that each $a\in\fg$ acts as the scalar $\chi(a)$.
Then \eqref{wsing2} and \eqref{wsing3} are equivalent to:
\begin{equation}\lbb{wsing4}
\V(R) = \T(R \tt (\kk_{\tr\ad} \bt \kk_{-\!\tr} )) \,, \quad
\T(R) = \V(R \tt (\kk_{-\!\tr\ad} \bt \kk_{\tr} )) \,,
\end{equation}
This can also be written as follows (cf.\ \deref{dtenmodw}(ii)):
\begin{equation}\lbb{wsing5}
\V(\Pi,U) = \T(\Pi\tt\kk_{\tr\ad} , U\tt\kk_{-\!\tr} ) \,, \quad
\T(\Pi,U) = \V(\Pi\tt\kk_{-\!\tr\ad} , U\tt\kk_{\tr} ) \,.
\end{equation}

\begin{theorem}\lbb{irrfacttens}
Let $V$ be an irreducible finite $\Wd$-module,
and let $R$ be an irreducible $(\dd \oplus \gld)$-submodule of\/ $\sing V$.
Then $V$ is a homomorphic image of\/ $\V(R)$.
In particular, every irreducible finite $\Wd$-module
is a quotient of a tensor module.
\end{theorem}
\begin{proof}
By \coref{cwsing}, we have $V=HR$. Consider the natural projection
\begin{equation*}
\pi\colon \V(R)=H\tt R \to HR=V \,, \qquad
h\tt u \mapsto hu \,.
\end{equation*}
Note that $\pi$ is $H$-linear.
Comparing \eqref{Wactionsing} with \eqref{wsing1}, we see that
\begin{equation*}
\bigl( (\id\tt\id)\tt_H \pi \bigr)
\bigl( (1 \tt \d_i) * (1 \tt u) \bigr)
= (1 \tt \d_i) * u \,, \qquad i=1,\dots,N \,, \; u\in R \,.
\end{equation*}
By $H$-bilinearity, this leads to
\begin{equation*}
\bigl( (\id\tt\id)\tt_H \pi \bigr) (a*v)
= a * \pi(v) \,, \qquad a\in\Wd \,, \; v\in\V(R) \,,
\end{equation*}
which means that $\pi$ is a homomorphism of $\Wd$-modules
(cf.\ \eqref{psprod6}).
\end{proof}

\subsection{Filtration of Tensor Modules}\lbb{swfilten}
Let $\V(R)$ be a tensor $\Wd$-module, as defined in \deref{dvmodw}(i).
Recall the canonical increasing filtration $\{\fil^p H\}$ of $H$
given by \eqref{filued}. We introduce an increasing filtration
of $\V(R)=H\tt R$ as follows:
\begin{equation}\lbb{wfilten1}
\fil^p \V(R) = \fil^p H \tt R \,, \qquad p=-1,0,\dots \;.
\end{equation}
Note that $\fil^{-1} \V(R) = \{0\}$, $\fil^0 \V(R) = \kk\tt R$.

The associated graded space of $\V(R)$ is
\begin{equation}\lbb{wfilten2}
\gr \V(R) = \bigoplus_{p\ge0} \gr^p \V(R) \,, \qquad
\gr^p \V(R) = (\fil^p H \tt R) / (\fil^{p-1} H \tt R) \,.
\end{equation}
We have isomorphisms of vector spaces:
\begin{equation}\lbb{wfilten3}
\gr^p \V(R) \simeq \gr^p H \tt R \simeq \symp^p \dd \tt R\,,
\end{equation}
where $\symp^p \dd$ is the $p$-th symmetric power of the vector space $\dd$.

Next, we study the action of the extended annihilation algebra $\ti\W$
on the filtration~\eqref{wfilten1}.

\begin{lemma}\lbb{lwfilten1}
For every $p\ge0$, we have{\rm:}

{\rm(i)}
$\dd \cdot \fil^p \V(R) \subset \fil^{p+1} \V(R)$,

{\rm(ii)}
$\N_\W \cdot \fil^p \V(R) \subset \fil^p \V(R)$,

{\rm(iii)}
$\W_1 \cdot \fil^p \V(R) \subset \fil^{p-1} \V(R)$.
\end{lemma}
\begin{proof}
Part (i) is obvious from definitions, since
\begin{equation*}
\d \cdot (h \tt u) = \d h \tt u \,, \qquad
\d\in\dd \,, \; h\in H \,, \; u\in R \,.
\end{equation*}

We will prove parts (ii) and (iii) by induction on~$p$.
For $p=0$, we have $\fil^0 \V(R) = \kk\tt R \subset \sing \V(R)$;
hence, (ii) and (iii) hold by the definition of a singular vector.
Now assume that (ii) is satisfied for some $p\ge0$.
Then it is enough to show that
\begin{equation*}
A \cdot (\d v) \in \fil^{p+1} \V(R) \qquad\text{for all}\quad
A\in\N_\W \,, \; \d\in\dd \,, \; v\in\fil^p \V(R) \,.
\end{equation*}
Note that, since $\ti\W=\dd+\N_\W$ (see \prref{normalizer}),
statements (i) and (ii) imply $\ti\W \cdot v \subset \fil^{p+1} \V(R)$.
Then we have
\begin{equation*}
A \cdot (\d v) = \d \cdot (A \cdot v) + [A,\d] \cdot v
\in \dd\cdot(\N_\W\cdot v) + \ti\W \cdot v \subset \fil^{p+1} \V(R) \,,
\end{equation*}
by part (i) and the inductive assumption. This proves (ii).

Similarly, assume that (iii) holds for some $p\ge0$. Then we want to
show that
\begin{equation*}
B \cdot (\d v) \in \fil^p \V(R) \qquad\text{for all}\quad
B\in\W_1 \,, \; \d\in\dd \,, \; v\in\fil^p \V(R) \,.
\end{equation*}
We have
\begin{equation*}
B \cdot (\d v) = \d \cdot (B \cdot v) + [B,\d] \cdot v
\in \dd\cdot(\W_1\cdot v) + \N_\W \cdot v \subset \fil^p \V(R) \,,
\end{equation*}
by (i), (ii) and the inductive assumption, because
$[B,\d] \in\W_0 \subset \N_\W$.
This completes the proof.
\end{proof}

\leref{lwfilten1}(ii) implies that the Lie algebra $\N_\W$ acts on the
associated graded space $\gr \V(R)$. By \leref{lwfilten1}(iii),
the same is true for the Lie algebra
$\N_\W / \W_1 = \ti\dd \oplus (\W_0/\W_1) \simeq \dd\oplus\gld$.
This action is described in the next two lemmas.

\begin{lemma}\lbb{lwfilten2}
For every $\d\in\dd$, $h\in \fil^p H$, $u\in R$, we have
\begin{equation*}
\ti\d \cdot (h \tt u) = h \tt \rho_R(\d) u \mod \fil^{p-1} \V(R) \,.
\end{equation*}
\end{lemma}
\begin{proof}
The proof is by induction on $p$ and is similar to that of
\leref{lwfilten1}(ii). First, for $p=0$, we have
$\fil^0 H = \kk$ and $1 \tt u \in \sing\V(R)$. Hence,
$\ti\d \cdot (1 \tt u) = 1 \tt \rho_R(\d) u$
by \eqref{rhosing}, \eqref{wsing7}.

Now assume the statement holds for $h\in \fil^p H$, and consider
$\ti\d \cdot (\d' h \tt u)$ for $\d'\in\dd$.
Note that, by \prref{pnw}(i), we have:
$[\ti\d, \d'] = [\ti\d, \ti\d' + \ga(\d') ]
= [\ti\d, \ti\d']  \in\ti\dd$.
{}From the inductive assumption, we get
$[\ti\d, \d'] \cdot (h \tt u) \in \fil^p \V(R)$. Therefore,
\begin{align*}
\ti\d \cdot (\d' h \tt u)
= \d' \cdot (\ti\d \cdot (h \tt u)) \mod \fil^p \V(R)
= \d' h \tt \rho_R(\d) u \mod \fil^p \V(R)
\end{align*}
by the inductive assumption.
\end{proof}
\begin{lemma}\lbb{lwfilten3}
The action of\/ $\gld \simeq \W_0/\W_1$ on the space
$\gr^p \V(R) \simeq \symp^p \dd \tt R$
is given by
\begin{equation*}
A \cdot (f \tt u) = Af \tt u + f \tt \rho_R(A) u \,,
\qquad A\in\gld \,, \; f\in \symp^p \dd \,, \; u\in R \,,
\end{equation*}
where $Af$ is the standard action of\/ $\gld$ on\/ $\symp^p \dd$.
\end{lemma}
\begin{proof}
The proof uses the same argument as in
Lemmas \ref{lwfilten1}(ii) and \ref{lwfilten2},
and the fact that via the isomorphisms $\W_0/\W_1 \simeq\gld$ and
$\W/\W_0 \simeq\dd$ the adjoint action $[A,\d]$
becomes the standard action of $\gld$ on $\dd$
(see \coref{cwbra}).
\end{proof}

When $R=\Pi\bt U$, the above two lemmas can be summarized as follows.

\begin{corollary}\lbb{cwfilten4}
We have\/
$\gr^p \V(\Pi,U) \simeq \Pi\bt (\symp^p \dd \tt U)$
as\/ $(\dd\oplus\gld)$-modules.
\end{corollary}


\subsection{Submodules of Tensor Modules}\lbb{swirten}
Let $\T(R)=H\tt R$ be a tensor module (see \deref{dtenmodw}).
We will assume that $R$ is an irreducible finite-dimensional
$(\dd\oplus\gld)$-module. Then $R=\Pi\bt U$, where $\Pi$ (respectively $U$)
is an irreducible finite-dimensional module over $\dd$ (respectively $\gld$).
In this case, $\T(R) = \T(\Pi,U)$.

As usual, we fix a basis $\{\d_i\}$ of $\dd$.
Recall that the action of $1\tt\d_i \in \Wd$ on an element
$1\tt u \in\kk\tt R \subset\T(R)$ is given by \eqref{wdgcd3}.
For us, it will be convenient to rewrite \eqref{wdgcd3} as follows
(making use of \eqref{deasa}):
\begin{equation}\lbb{wirten2}
\begin{split}
(1\tt \d_i)&*(1\tt u)
= (1 \tt 1) \tt_H (1 \tt (\d_i+\ad\d_i)u)
\\
&+ \sum_{j=1}^N\, (1 \tt 1) \tt_H (\d_j  \tt e_i^j u)
- \sum_{j=1}^N\, (1 \tt \d_j) \tt_H (1  \tt (e_i^j+\de_i^j) u)
\,.
\end{split}
\end{equation}
Introduce the following notation:
\begin{equation}\lbb{sdu2}
s(\d_i, u) = \sum_{j=1}^N\, \d_j \tt e_i^j u  \,,
\qquad u \in R \,, \; i=1,\dots,N \,.
\end{equation}
By linearity, we define $s(\d, u)$ for all $\d\in\dd$.
Then $s(\d, u)$ does not depend on the choice of basis
$\{\d_i\}$ of $\dd$ (cf.\ \eqref{epab}, \eqref{wdepa}).

Since $\T(R)=H\tt R$, any element $v\in\T(R)$
can be written uniquely in the form
\begin{equation}\lbb{wirten3}
v = \sum_{I\in\ZZ_+^N} \d^{(I)} \tt v_I
\,, \qquad v_I \in R \,,
\end{equation}
where $\d^{(I)} \in H$ are given by \eqref{dpbw}.
Note that the above sum is finite, i.e., $v_I \ne 0$ only for finitely
many $I$.
{}From \eqref{wirten2}--\eqref{wirten3} and $H$-bilinearity, we find
\begin{equation}\lbb{wirten5}
\begin{split}
(1\tt \d_i)*v
&= \sum_I (1 \tt \d^{(I)}) \tt_H (1 \tt (\d_i+\ad\d_i)v_I)
\\
&+ \sum_I\, (1 \tt \d^{(I)}) \tt_H s(\d_i,v_I)
\\
&- \sum_I \sum_{j=1}^N\, (1 \tt \d^{(I)} \d_j) \tt_H
(1  \tt (e_i^j+\de_i^j) v_I)
\,.
\end{split}
\end{equation}

\begin{definition}\lbb{dcoeff}
The nonzero elements $v_I$ in the expression \eqref{wirten3}
are called {\em coefficients\/} of $v\in\T(R)$.
%
For a submodule $M\subset\T(R)$, we denote by $\coef M$
the subspace of $R$
linearly generated by all coefficients 
of elements of~$M$.
\end{definition}

Recall that $\T(R)$ has a filtration given by
$\fil^p \T(R) = \fil^p H \tt R$ (cf.\ \eqref{wfilten1} and
\eqref{wsing4}). We have: $\fil^{-1} \T(R) = \{0\}$,
$\fil^0 \T(R) = \kk\tt R$
and $\fil^1 \T(R) = (\kk+\dd)\tt R$.

\begin{lemma}\lbb{noconstants}
For any nontrivial proper $\Wd$-submodule $M$ of\/ $\T(R)$,
we have $M \cap\fil^0 \T(R)  = \{0\}$.
\end{lemma}
\begin{proof}
Let $M_0$ be the set of all $u\in R$ such that $1\tt u \in M$.
By \eqref{wirten2} and \reref{rlmod}(i),
we have:
\begin{equation*}
(\d_i+\ad\d_i)u \in M_0 \,, \quad
(e_i^j+\de_i^j) u \in M_0 \quad\text{for all}\quad
i,j=1,\dots,N \,, \; u\in M_0 \,.
\end{equation*}
This means that $M_0$ is a $(\dd \oplus \gld)$-submodule of $R$.
Since $R$ is irreducible, either $M_0 = \{0\}$ or $M_0 = R$.
In the latter case, we obtain $M \supset H \tt M_0 = \T(R)$,
which is a contradiction. Therefore, $M \cap (\kk \tt R) = \{0\}$.
\end{proof}
\begin{corollary}\lbb{cwirr0}
If\/ $\sing\T(R) = \fil^0 \T(R)$, then the tensor\/ $\Wd$-module\/
$\T(R)$ is irreducible.
\end{corollary}
\begin{proof}
If $M\subset\T(R)$ is a nontrivial proper submodule,
then by \thref{twsing1} it contains a nonzero singular vector.
This contradicts \leref{noconstants}.
\end{proof}

\coref{cwirr0} will play a crucial role in our
classification of irreducible finite $\Wd$-modules.

\begin{lemma}\lbb{lcoeff}
For any nontrivial proper\/ $\Wd$-submodule\/ $M$ of\/ $\T(R)$,
we have\/ $\coef M = R$.
\end{lemma}
\begin{proof}
Pick a nonzero element $v\in M$ and write in the form \eqref{wirten3}.
Then, for fixed $I$, the coefficient multiplying $1 \tt \d^{(I)}$
in the right-hand side of \eqref{wirten5} equals
\begin{equation}\lbb{wirten4}
s(\d_i,v_I) + 1  \tt \d_i v_I
+ \text{terms in}\; \kk\tt (\gld+\kk)(\coef M) \,.
\end{equation}
By \reref{rlmod}(ii), this is an element of $M$.
Hence, for each coefficient $v_I$ of $v$, we have
$e_i^j v_I \in \coef M$.
Then from \eqref{wirten4} we also get $\d_i v_I \in \coef M$.
Therefore, $\coef M$ is a nontrivial $(\dd \oplus \gld)$-submodule
of $R$. But $R$ is irreducible; hence, $\coef M=R$.
\end{proof}
\begin{lemma}\lbb{lsmdv}
Let\/ $M$ be a nontrivial proper\/ $\Wd$-submodule of\/ $\T(R)$.
Then for every\/ $\d \in \dd$ and\/ $u \in R$, there is
a unique element\/ $s_M(\d, u) \in M$ such that
\begin{equation}\lbb{wcsv8}
s_M(\d, u) = s(\d, u) \mod\fil^0\T(R) \,,
\end{equation}
where\/ $s(\d, u)$ is given by\/ \eqref{sdu2}.
The element\/ $s_M(\d, u)$ depends linearly on both\/ $\d$ and~$u$.
\end{lemma}
\begin{proof}
Uniqueness follows from \leref{noconstants}. {}From uniqueness we deduce
that $s_M(\d, u)$ depends linearly on $\d$ and $u$.
Then, to prove existence, it is enough to consider the case
$\d=\d_i$ and $u=v_I$ for some $v\in M$
(because $R=\coef M$ by \leref{lcoeff}).
In this case $s_M(\d_i, v_I)$ is exactly the element~\eqref{wirten4}.
\end{proof}

Elements $s_M(\d, u)$ will be used in the next subsection
to determine all singular vectors of~$M$.

\subsection{Computation of Singular Vectors}\lbb{swcsv}
In this subsection, we continue to use the notation of \seref{swirten}.
Our goal is to find all singular vectors of $\T(R)=\T(\Pi,U)$.
Given a nontrivial proper $\Wd$-submodule $M$ of $\T(R)$,
we also find all singular vectors of $M$.
These results will be used in  \seref{swirr}
to classify irreducible finite $\Wd$-modules.

First, we consider the case when the $\gld$-action on $R$ is trivial.

\begin{proposition}\lbb{om0irr}
For any irreducible finite-dimensional\/ $\dd$-module\/ $\Pi$,
we have{\rm:}

{\rm(i)}
$\sing\T(\Pi,\kk) = \fil^0\T(\Pi,\kk);$

{\rm(ii)}
$\T(\Pi,\kk)$ is an irreducible\/ $\Wd$-module.
\end{proposition}
\begin{proof}
Pick a singular vector $v\in\T(\Pi, \kk)$,
and write it in the form \eqref{wirten3}.
Then, by~\eqref{wirten5},
\begin{equation*}
(1 \tt \d_i) * v
= \sum_I\, (1 \tt \d^{(I)}) \tt_H (1 \tt \d_i v_I)
- \sum_I\, (1 \tt \d^{(I)} \d_i) \tt_H (1 \tt v_I) \,.
\end{equation*}
Now \reref{rwsing} implies that $v_I=0$ whenever $|I| \ge 1$.
This proves (i).

(ii) follows from (i) and \coref{cwirr0}.
\end{proof}


\begin{lemma}\lbb{lwcsv4}
For any irreducible finite-dimensional\/ $(\dd\oplus\gld)$-module\/ $R$,
the tensor\/ $\Wd$-module\/ $\T(R)$ satisfies
\begin{equation*}
\fil^0 \T(R) \subset \sing\T(R) \subset \fil^1 \T(R) \,.
\end{equation*}
\end{lemma}
\begin{proof}
First of all, by \prref{om0irr}(i), we can assume that
the $\gld$-action on $R$ is nontrivial.
Since $\fil^0 \T(R) = \kk\tt R$,
the first inclusion follows from \eqref{wirten2} and \reref{rwsing}.
To prove the second one,
pick a nonzero singular vector $v$ and write in the form \eqref{wirten3}.
Then $(1\tt \d_i)*v$ is given by formula \eqref{wirten5}.
The coefficient multiplying $1 \tt \d^{(I)}$ in \eqref{wirten5}
is given by \eqref{wirten4}.
By \reref{rwsing}, this coefficient must vanish whenever $|I|>1$.
Hence, $s(\d_i,v_I) = 0$ for all $i$, which implies $(\gld)v_I = 0$.
Therefore, $v_I=0$. This proves that
$\sing\T(R) \subset \fil^1 \T(R)$.
\end{proof}
\begin{lemma}\lbb{lwcsv7}
An element
\begin{equation}\lbb{wcsv1}
v = \sum_{k=1}^N\, \d_k \tt v^k \in \dd\tt R \subset \T(R)
\end{equation}
is a singular vector iff it satisfies the equations
\begin{equation}\lbb{wcsv3}
(e_i^j+\de_i^j) v^k + (e_i^k+\de_i^k) v^j = 0
\qquad\text{for all}\quad i,j,k=1,\dots,N \,.
\end{equation}
In this case, for the action\/ $\rho_{\sing}$ of\/ $\gld$ on\/ $v$,
we have $($see \eqref{sdu2}$)${\rm:}
\begin{equation}\lbb{wcsv4}
\rho_{\sing}(e_i^k) v = -s(\d_i, v^k) \mod\fil^0\T(R) \,.
\end{equation}
\end{lemma}
\begin{proof}
As a special case of \eqref{wirten5}, we have:
\begin{equation*}
\begin{split}
(1 \tt \d_i)*v
&= \sum_{k=1}^N\, (1 \tt \d_k) \tt_H (1 \tt (\d_i+\ad\d_i) v^k)
+ \sum_{k=1}^N\, (1 \tt \d_k) \tt_H s(\d_i, v^k)
\\
&- \sum_{k,j=1}^N\, (1 \tt \d_k \d_j) \tt_H (1 \tt (e_i^j+\de_i^j) v^k) \,.
\end{split}
\end{equation*}
For $k\le j$, the coefficient multiplying $1 \tt \d_k \d_j$
is up to a sign equal to
\begin{equation}\lbb{wcsv5}
1 \tt (e_i^j+\de_i^j) v^k + 1 \tt (e_i^k+\de_i^k) v^j \,.
\end{equation}
By \reref{rwsing}, $v$ is a singular vector iff
this coefficient vanishes for all $j,k$. This proves \eqref{wcsv3}.
On the other hand, the coefficient multiplying $1 \tt \d_k$
is equal to $s(\d_i, v^k)$ modulo $\fil^0 \T(R)$.
Then \eqref{wcsv4} follows from~\eqref{Wactionsing}.
\end{proof}

Our next result describes all singular vectors in
a tensor $\Wd$-module $\T(R)$.

\begin{theorem}\lbb{twsing2}
For any irreducible finite-dimensional\/ $(\dd\oplus\gld)$-module\/ $R$,
we have{\rm:}
\begin{equation}\lbb{wcsv6}
\sing\T(R) = \fil^0 \T(R)
+ \{ s(\d, u) \st \d \in \dd \,, \; u \in R_0 \} \,,
\end{equation}
where $s(\d, u)$ is defined by \eqref{sdu2}
and\/ $R_0$ is the subspace of all\/ $u\in R$ satisfying the equations
\begin{equation}\lbb{wcsv7}
(e_i^j+\de_i^j) e_l^k u + (e_i^k+\de_i^k) e_l^j u = 0 \,,
\qquad  i,j,k,l=1,\dots,N \,.
\end{equation}
The subspace\/ $R_0$ is either\/ $\{0\}$ or the whole\/ $R$.
\end{theorem}
\begin{proof}
When $R=\Pi\bt\kk$, \eqref{wcsv6} follows from \prref{om0irr}(i),
because in this case all $s(\d, u)=0$.
Let us assume that the $\gld$-action on $R$ is nontrivial,
and denote the space in the right-hand side of \eqref{wcsv6} by $S$.
Notice that $s(\d, u)$ is a singular vector iff $u\in R_0$,
because for $v=s(\d_l, u)$ we have $v^k = e_l^k u$
and \eqref{wcsv3} becomes \eqref{wcsv7}.
Hence, Lemmas \ref{lwcsv4} and \ref{lwcsv7} imply $S \subset \sing\T(R)$.
{}From these lemmas, we also deduce that the action
$\rho_{\sing}$ of $\gld$ on $\sing\T(R)$ maps $\sing\T(R)$ into~$S$.

Consider the finite-dimensional $\gld$-module $\sing\T(R) / \fil^0 \T(R)$.
We claim that its decomposition as a direct sum of irreducibles
does not contain the trivial $\gld$-module.
Indeed, let $v\in\sing\T(R)$ be such that
$\rho_{\sing}(e_i^k) v \in\fil^0 \T(R)$ for all $i,k$.
We want to show that $v \in\fil^0 \T(R)$.
Without loss of generality,
we can assume that $v \in\dd\tt R$.
Then, by \eqref{wcsv4}, all $s(\d_i, v^k) = 0$,
and from \eqref{sdu2}, $e_i^j v^k = 0$ for all $i,j,k$.
This implies $v^k = 0$ for all $k$, and hence $v=0$, 
which proves our claim.

Therefore, the $\gld$-action on $\sing\T(R) / \fil^0 \T(R)$
is surjective and $\sing\T(R)=S$.
Next, it is clear that for $R=\Pi\bt U$ we have $R_0 = \Pi\bt U_0$, where
$U_0$ is the subspace of all $u\in U$ satisfying \eqref{wcsv7}.
We claim that $U_0$ is a $\gld$-submodule of $U$.
Since $U$ is an irreducible $\gld$-module, this would imply
that either $U_0=\{0\}$ or $U_0=U$.
Now if $u\in R_0$, then $v=s(\d_l, u)$ is a singular vector
for all $l$. By \eqref{wcsv4}, all $s(\d_i, v^k)$ are singular vectors
too. Hence, $v^k = e_l^k u$ belongs to $R_0$ for all $k,l$.
Therefore, $(\gld)R_0 \subset R_0$, which implies
$(\gld)U_0 \subset U_0$.
\end{proof}
\begin{corollary}\lbb{cwcsv2}
If\/ the\/ $\Wd$-module\/ $\T(R)$ is not irreducible, then
equations\/ \eqref{wcsv7} are satisfied for all~$u\in R$.
\end{corollary}
\begin{proof}
This follows from \coref{cwirr0} and \thref{twsing2}.
\end{proof}

Next, we find all singular vectors in
a nontrivial proper $\Wd$-submodule $M$ of $\T(R)$.
Recall the elements $s_M(\d, u) \in M$, constructed
in \leref{lsmdv}.

\begin{theorem}\lbb{twcsv5}
For any nontrivial proper\/ $\Wd$-submodule\/ $M$ of\/ $\T(R)$,
we have{\rm:}

{\rm(i)}
$\sing M = M \cap \fil^1 \T(R)
= \{ s_M(\d, u) \st \d\in\dd \,, \; u \in R \} ;$

{\rm(ii)}
$\sing\T(R) = \fil^0 \T(R) \oplus \sing M$
as\/ $(\dd\oplus\gld)$-modules.
\end{theorem}
\begin{proof}
(i)
Note that $\sing M \subset M \cap \fil^1 \T(R)$ by \leref{lwcsv4}.
Conversely, pick $v' \in M \cap \fil^1 \T(R)$, and write
$v' = v + 1\tt u$ with $v\in\dd\tt R$, $u \in R$.
Since $1\tt u \in\sing\T(R)$, the vector $v'$ is singular if
and only if $v$ is. In the proof of \leref{lwcsv7} we saw that
the coefficient multiplying $1 \tt \d_k \d_j$ in $(1 \tt \d_i)*v$
is up to a sign equal to \eqref{wcsv5}.
By \reref{rwsing}, $(1 \tt \d_i)*v'$ has the same coefficient.
Now \reref{rlmod}(ii) implies that the elements \eqref{wcsv5} belong to $M$.
Hence, they vanish by \leref{noconstants}. Then, by \leref{lwcsv7},
$v$ is a singular vector, and $v' \in\sing M$.

This proves the first equality in (i).
The second equality follows from the first one,
\leref{lsmdv} and \thref{twsing2}.

(ii)
The sum is direct because of \leref{noconstants}.
The equality follows from part (i), \leref{lsmdv} and \thref{twsing2}.
\end{proof}
\begin{corollary}\lbb{cwcsv3}
Let\/ $R$ be an irreducible finite-dimensional\/ $(\dd\oplus\gld)$-module,
and let\/ $M,M'$ be two nontrivial proper\/
$\Wd$-submodules of\/ $\T(R)$. Then\/
$\sing M = \sing M'$.
\end{corollary}
\begin{proof}
Consider the canonical projection of $(\dd\oplus\gld)$-modules
\begin{equation*}
\pi\colon \sing\T(R) \to \sing\T(R) / \fil^0 \T(R) 
\subset \gr^1 \T(R) \,.
\end{equation*}
By \thref{twcsv5}(ii), the restriction of $\pi$
to $\sing M$ is an isomorphism. On the other hand, combining
\eqref{wsing5} with \coref{cwfilten4}, we obtain isomorphisms
of $(\dd\oplus\gld)$-modules
\begin{equation*}
\fil^0 \T(R) \simeq \ti\Pi \bt \ti U \,, \quad
\gr^1 \T(R) \simeq \ti\Pi \bt (\dd \tt \ti U) \,,
\end{equation*}
where
\begin{equation*}
R=\Pi\bt U \,, \quad
\ti\Pi = \Pi\tt\kk_{-\!\tr\ad} \,, \quad
\ti U = U\tt\kk_{\tr} \,.
\end{equation*}
The $\gld$-module $\ti U$ is irreducible.
Say that $\Id\in\gld$ acts as the scalar $c$ on $\ti U$.
Then it acts as $c+1$ on $\dd \tt \ti U$. It follows that
$\sing M$ is precisely the set of all vectors $v\in\sing\T(R)$
such that $\Id\cdot v = (c+1)v$. The same is true for $M'$ instead of~$M$.
\end{proof}

\subsection{Irreducible Finite $\Wd$-Modules}\lbb{swirr}
This subsection contains our main results about irreducible
finite $\Wd$-modules.
As before, let $\Pi$ $($respectively $U)$ be an irreducible
finite-dimensional representation of $\dd$ $($respectively $\gld)$.
First, we determine which tensor $\Wd$-modules are irreducible.

\begin{theorem}\lbb{twirten}
The tensor\/ $\Wd$-module\/ $\T(\Pi, U)$ is irreducible
if and only if, as a\/ $\gld$-module, $U$ is not isomorphic to\/
$\bigwedge^n \dd^*$  for any~$n\ge1$.
\end{theorem}
\begin{proof}
Assume that $\T(\Pi, U)$ is not irreducible.
Then, by \coref{cwcsv2}, equations \eqref{wcsv7} are satisfied for all
$u\in R$. In the special case $i=j=k=l$ they give
\begin{equation}\lbb{eiei}
(e_i^i+1) e_i^i u = 0 \qquad\text{for all}\quad u \in R=\Pi\bt U \,.
\end{equation}
We claim that the $\gld$-module $U$ is isomorphic to 
$\Om^n := \bigwedge^n \dd^*$ for some $n$.
To prove this, first note that the matrix $\Id\in\gld$ acts as
a scalar on $U$, and the module $U$ remains irreducible
when restricted to $\sld$. Then $U$ has a highest weight vector $v$,
and the $\gld$-module $U$ is uniquely determined by its highest weight,
i.e., by the eigenvalues $\la_i$ of $e_i^i$ on $v$. Furthermore,
all $\la_i - \la_{i+1}$ are non-negative integers
(see, e.g., \cite[Chapter~VII]{Se}). 
But by \eqref{eiei}, all $\la_i = 0$ or $-1$; hence
the $N$-tuple $(\la_1,\dots,\la_N)$ has the form
$(0,\dots,0,-1,\dots,-1)$. The module $\Om^n$ has such a highest weight,
where the number of $-1$'s is $n$.
Therefore, $U \simeq \Om^n$, and the case $n=0$ 
is excluded by \prref{om0irr}(ii).

Next, recall the $\Pi$-twisted pseudo de Rham complex \eqref{domd2},
and introduce the shorthand notation
\begin{equation}\lbb{tnin}
T^n := \T(\Pi, \Om^n) \,, \quad
I^n := \di_\Pi (T^{n-1}) \subset T^n \,.
\end{equation}
Since the differential $\di_\Pi$ is a homomorphism of $\Wd$-modules,
the image $I^n$ is a submodule of $T^n$ for each $n=1,\dots,N$.
It is easy to see from \prref{pdomd2} and \leref{ldw2} that
$I^n$ is nontrivial and proper. Therefore, the tensor modules
$T^n$ are not irreducible for~$n\ge1$.
\end{proof}
\begin{corollary}\lbb{cwirr1}
Let\/ $R$ be an irreducible finite-dimensional\/ $(\dd\oplus\gld)$-module.
Then the\/ $\Wd$-module\/ $\T(R)$ is irreducible if and only if\/
$\sing\T(R) = \fil^0 \T(R)$.
\end{corollary}
\begin{proof}
In one direction, the statement is exactly \coref{cwirr0}.
The opposite direction follows from
\thref{twsing2} and the proof of \thref{twirten}.
\end{proof}

Our next goal is to study the submodules $I^n$ of $T^n$
(see \eqref{tnin}).

\begin{lemma}\lbb{lwirr1}
For\/ $1\le n\le N$,
the\/ $\Wd$-submodule\/ $I^n \subset T^n$ 
has the following properties{\rm:}

{\rm(i)}
$I^n$ is nontrivial and proper\/$;$

{\rm(ii)}
$\sing I^n = \di_\Pi(\kk\tt\Pi\tt\Om^{n-1});$

{\rm(iii)}
$I^n$ is generated by\/ $\sing I^n$ as an\/ $H$-module$;$

{\rm(iv)}
Any nontrivial proper submodule\/ $M$ of\/ $T^n$ contains $I^n;$

{\rm(v)}
$I^n$ is an irreducible $\Wd$-module.
\end{lemma}
\begin{proof}
(i) is easy to see from \prref{pdomd2} and \leref{ldw2}.

(ii) Formula \eqref{dw3} and \leref{lsmdv} imply
\begin{equation*}
s_{I^n}(\d_i, u\tt\al) = \di_\Pi (1 \tt u \tt \io_{\d_i}\al)
\,, \qquad u\in\Pi \,, \; \al\in\Om^n \,, \; i=1,\dots,N \,.
\end{equation*}
Then (ii) follows from \thref{twcsv5}(i) and the fact that
$\Om^{n-1}$ is linearly spanned by all~$\io_{\d_i}\al$.

(iii) is obvious from (ii) and the $H$-linearity of~$\di_\Pi$.

(iv) By \coref{cwcsv3}, we have $\sing M = \sing I^n$.
Then $M \supset H(\sing M) = H(\sing I^n)$, which is
equal to $I^n$ by part~(iii).

(v) is obvious from~(iv).
\end{proof}

Note that, from the exactness of the complex \eqref{domd2},
we have $I^1 \simeq T^0 = \T(\Pi,\Om^0)=\T(\Pi,\kk)$.

\begin{lemma}\lbb{lwirr2}
For\/ $1\le n\le N-1$, $I^n$ is the unique nontrivial proper
$\Wd$-submodule of\/ $T^n$.
\end{lemma}
\begin{proof}
If $M$ is a nontrivial proper submodule of $T^n$, it contains $I^n$.
The image $\di_\Pi M$ is a submodule of $I^{n+1}$;
hence, $\di_\Pi M$ is either $\{0\}$ or the whole $I^{n+1}$.
But the kernel of $\di_\Pi$ is equal to $I^n$,
by the exactness of the complex \eqref{domd2}.
We obtain that either $M=I^n$ or $M=T^n$.
\end{proof}

Now we can classify all irreducible finite $\Wd$-modules.

\begin{theorem}\lbb{twmod}
Any irreducible finite\/ $\Wd$-module is isomorphic to
one of the following\/{\rm:}

{\rm(i)}
Tensor modules\/ $\T(\Pi, U)$,
where\/ $\Pi$ is an irreducible finite-dimensional\/ $\dd$-module,
and\/ $U$ is an irreducible finite-dimensional\/ $\gld$-module
not isomorphic to\/ $\Om^n=\bigwedge^n \dd^*$ for any\/
$1\le n \le\dim\dd;$

{\rm(ii)}
Images\/ $\di_\Pi \T(\Pi, \Om^n)$,
where\/ $\Pi$ is an irreducible finite-dimensional\/ $\dd$-module,
and\/ $1\le n \le\dim\dd-1$ $($see~\eqref{domd2}$)$.
\end{theorem}
\begin{proof}
Let $V$ be an irreducible finite $\Wd$-module.
Then, by \thref{irrfacttens} and \eqref{wsing4},
$V$ is a quotient of some tensor module $\T(R) = \T(\Pi,U)$.
If $U \not\simeq \Om^n$
as a $\gld$-module for any $n\ge1$,
then $\T(R)$ is irreducible by \thref{twirten},
and in this case $V\simeq\T(R)$.

Assume that $U \simeq \Om^n$ for some $n\ge1$;
then $\T(R) \simeq \T(\Pi,\Om^n) = T^n$ (see \eqref{tnin}).
Now if $n\le N-1$, $N=\dim\dd$, \leref{lwirr2} implies
that $V \simeq T^n/I^n$. By the exactness of \eqref{domd2},
we get $V \simeq I^{n+1} = \di_\Pi \T(\Pi, \Om^n)$.

Finally, it remains to consider the case when $V$ is a quotient
of $T^N$. Then $V \simeq T^N/M$, where $M\supset I^N$
due to \leref{lwirr1}(iv). Now \prref{pdomd2} implies that
$V$ is finite dimensional; hence, $\Wd$ acts trivially on it
by \exref{elv0}. So in this case $V$ cannot be irreducible.
\end{proof}
\begin{theorem}\lbb{twnoiso}
The irreducible finite\/ $\Wd$-modules listed in \thref{twmod}
satisfy{\rm:}

{\rm(i)}
$\sing\T(\Pi, U) \simeq
( \Pi\tt\kk_{-\!\tr\ad} ) \bt ( U \tt\kk_{\tr} )$
as\/ $(\dd\oplus\gld)$-modules$;$

{\rm(ii)}
$\sing(\di_\Pi \T(\Pi, \Om^n)) \simeq
( \Pi\tt\kk_{-\!\tr\ad} ) \bt ( \Om^n \tt\kk_{\tr} )$
as\/ $(\dd\oplus\gld)$-modules.
\\
In particular, no two of them are isomorphic to each other.
\end{theorem}
\begin{proof}
First, note that if $\be\colon V\to V'$ is a homomorphism
of $\Wd$-modules, then its restriction to $\sing V$ is
a homomorphism of $(\dd\oplus\gld)$-modules $\sing V \to \sing V'$.
In particular, if $V$ and $V'$ are isomorphic, then
$\sing V \simeq \sing V'$.

(i)
If $\T(R) = \T(\Pi,U)$ is irreducible, then by \coref{cwirr1},
$\sing\T(R) = \fil^0 \T(R) = \kk \tt R$.
Now (i) follows from \eqref{wsing7} and~\eqref{wsing5}.

{\rm(ii)}
By \leref{lwirr1}(ii), we have:
$\sing(\di_\Pi \T(\Pi, \Om^n)) = \di_\Pi(\fil^0 \T(\Pi, \Om^n))$.
But 
$\fil^0 \T(\Pi, \Om^n) \simeq
( \Pi\tt\kk_{-\!\tr\ad} ) \bt ( \Om^n \tt\kk_{\tr} )$
is an irreducible $(\dd\oplus\gld)$-module.
Therefore, $\di_\Pi$ is an isomorphism from $\fil^0 \T(\Pi, \Om^n)$
onto $\sing(\di_\Pi \T(\Pi, \Om^n))$.
\end{proof}
\begin{remark}\lbb{rwhom}
Let $R$ and $R'$ be two non-isomorphic irreducible finite-dimensional
$(\dd\oplus\gld)$-modules.
Using \thref{twcsv5}(ii) and the same argument as in the proof of
\thref{twnoiso}, one can show that the only nonzero homomorphisms
of $\Wd$-modules $\T(R) \to \T(R')$ are, up to a constant factor, 
the differentials~$\di_\Pi$.
\end{remark}

\section{Classification of Irreducible Finite $\Sd$-Modules}\lbb{ssmod}

In this section we adapt the classification results of \seref{swmod}
to the case of the Lie pseudoalgebra $\Sd$.
Our main result is \thref{tsmod}.

\subsection{Singular Vectors}\lbb{ssings}

Recall that the annihilation algebra $\S$ of $\Sd$ possesses a
decreasing filtration $\{\S_p\}_{p\geq -1}$ by subspaces of finite
codimension, as given by \eqref{sp}. This filtration is compatible
with that of $\W$ given by \eqref{wp}, \eqref{filpwn}, in the sense
made clear by \leref{lsw3} and \prref{psw2}. In particular, we know
that $\S\subset \W \simeq W_N$ is a graded Lie algebra isomorphic to
$S_N$, the grading being given by the eigenspace decomposition
with respect to the adjoint action of the element $\what \E \in
\W$ described in \deref{deul2}. We denote the $i$-eigenspace by
$\s_i$; hence, $\S_p = \prod_{i\geq p} \s_i$.
Note that we can do the same with the extended annihilation
algebra $\wti \S$, as $\what \E$ commutes with $\wti\dd$. 
We denote the corresponding eigenspaces by $\wti\s_i$;
then $\wti\dd \subset \wti\s_0$ 
and $\wti\s_i = \s_i$ for $i\ne0$.

In analogy with the case of $\Wd$, for an $\S$-module $V$, we
denote by $\ker_p V$ the space of all elements $v\in V$ that are
killed by $\S_p$. We denote by $\ker V$ the space $\ker_{-1} V$
killed by $\S=\S_{-1}$. Then the module $V$ is {\em conformal\/} iff
$V = \bigcup_p \ker_p V$.
Any $\Sd$-module has a natural structure of a conformal module over the
extended annihilation algebra $\wti \S = \dd \sd \S$
(see \prref{preplal2}). The normalizer of $\S_p$ in $\wti\S$ 
was computed in \seref{sns}. It is independent of $p\geq 0$, 
and is denoted by $\N_\S$. Note that
$\N_\S = \prod_{i\geq 0} \wti\s_i$
and $\wti\S=\wti\s_{-1} \oplus \N_\S$ as a vector space. 

Each $\ker_p V$ is a module over the finite-dimensional quotient
$\N_\S/\S_p$; moreover, the Lie algebra $\N_\S/\S_1 = \what \dd
\oplus (\S_0/\S_1)$ is isomorphic to the direct sum $\dd \oplus
\sld$.

\begin{definition}\lbb{dssing}
A {\em singular vector\/} in an $\Sd$-module $V$ is an element $v\in V$
such that $\S_1 \cdot v = 0$. The space of singular vectors
$\ker_1 V$ is also denoted by $\sing V$. 
\end{definition}
\begin{theorem}\lbb{tssing1}
For any nontrivial finite\/ $\Sd$-module\/ $V$, 
we have\/ $\sing V \neq \{0\}$, 
and the space\/ $\sing V/\ker V$ is finite dimensional.
\end{theorem}
\begin{proof}
The proof is the same as that of \thref{twsing1}, making use of
\prref{NrepstypeS} instead of \prref{Nreps}.
\end{proof}

Recall that $\Sd$ is generated over $H$ by the elements $s_{ab}$
defined in \eqref{achib}. It will be convenient to introduce the
notation
\begin{equation}\lbb{bard}
\bar\d = \d+\chi(\d) \,, \qquad \d\in\dd \,,
\end{equation}
and
\begin{equation}\lbb{sij}
s_{ij} \equiv s_{\d_i\d_j} 
= \bar\d_i \tt \d_j - \bar\d_j \tt \d_i - 1 \tt [\d_i,\d_j] \,,
\end{equation}
where, as before, $\{\d_i\}$ is a fixed basis of $\dd$.

\begin{remark}\lbb{rsings}
By \eqref{sp} and \prref{preplal2}, a vector $v\in V$ is
singular if and only if 
\begin{equation*}
s_{ij}* v \in (\fil^2 H \tt \kk) \tt_H V \,, 
\qquad i,j=1,\dots,N \,,
\end{equation*}
or, equivalently, 
\begin{equation*}
s_{ij}* v \in (\kk\tt \fil^2 H) \tt_H V \,,
\qquad i,j=1,\dots,N \,.
\end{equation*}
\end{remark}

\subsection{Tensor Modules for $\Sd$}\lbb{stens2}
Let $R$ be a finite-dimensional $(\dd \oplus \sld)$-module, with an
action denoted by $\what\rho_R$. Then the isomorphism 
$\N_\S/\S_1 \simeq \dd \oplus \sld$ can be employed
to make $R$ an $\N_\S$-module with a trivial action of $\S_1$. 
For example, the action of the subalgebra $\what\dd\subset\N_\S$
is given by:
\begin{equation}\lbb{hatdr}
\what\d \cdot u = \what\rho_R(\d) u \,,
\qquad \d\in\dd \,, \; u \in R \,.
\end{equation}
Consider the induced $\wti\S$-module $V=\Ind_{\N_\S}^{\wti \S} R$. 
Since as a vector space $\wti\S = \dd \oplus \N_\S$
(see \prref{normalizerS}), 
as an $H$-module $V$ is isomorphic to the free module $H \tt R$.

\begin{definition}\lbb{dtens}
The $\Sd$-module $\Ind_{\N_\S}^{\wti \S} R$ 
constructed above will be denoted by $\VS(R)$, 
and will be called a {\em tensor module\/} for the Lie pseudoalgebra $\Sd$.
If $R$ is an irreducible $(\dd \oplus \sld)$-module
isomorphic to $\Pi \bt U$, then we will also write 
$\VS(R) = \VS(\Pi, U)$. 
\end{definition}

The name tensor module is justified by the
fact that, as we will show in \thref{prrestriction} below, 
the $\Sd$-module $\VS(R)$ is
the restriction to $\Sd\subset \Wd$ of a tensor module for $\Wd$
(see also \reref{rvrind}).

\begin{theorem}\lbb{prstensoruniversal}
Let\/ $V$ be an\/ $\Sd$-module, and let\/ $R$ be a\/ 
$(\dd \oplus \sld)$-submodule of\/ $\sing V$. Then\/ $HR$ is
an $\Sd$-submodule of\/ $V$, and there is a natural surjective
homomorphism\/ $\VS(R) \to HR$.
In particular, every irreducible
finite\/ $\Sd$-module is a quotient of a tensor module.
\end{theorem}
\begin{proof}
Since $\wti\S = \dd \oplus \N_\S$ as a vector space, and 
$R\subset\sing V$, it follows that $HR$ is preserved
by the action of $\wti\S$. Then by \prref{preplal2}, $HR$ is an 
$\Sd$-submodule of $V$. The existence of a natural surjective
homomorphism $\VS(R) \to HR$ follows from the definition of $\VS(R)$.
Finally, if $V$ is irreducible and finite, then by \thref{tssing1}, 
$\sing V \ne\{0\}$, and we have $H(\sing V) = V$.
\end{proof}
\begin{lemma}\lbb{lswtranslation}
The unique injection\/
$\iota\colon \Sd \to \Wd$ induces an injective Lie
algebra homomorphism\/ 
$\iota_*\colon \N_\S/\S_1 \to \N_\W/\W_1$. 
The homomorphism\/ $\iota_*$ satisfies\/
$\iota_* (\what \dd) \subset \wti \dd \oplus \kk \Id$.
More precisely $($see \eqref{ns1}$)$,
\begin{equation}\lbb{io*d}
\iota_* (\what \d) 
= \wti \d + \frac{1}{N}(\tr\ad - \chi)(\d) \Id
\,, \qquad \d\in\dd \,.
\end{equation}
Furthermore, $\iota_*$ embeds\/ $\S_0/\S_1 \simeq \sld$ as the Lie
subalgebra\/ $\sld \subset \gld \simeq \W_0/\W_1$.
\end{lemma}
\begin{proof}
By \prref{psw2}, the induced Lie algebra homomorphism
$\io\colon\S\to\W$ is injective and satisfies $\io(\S_1) \subset \W_1$.
Hence, $\io_*$ is injective. The rest of the lemma follows from
\eqref{ns1} and \coref{csw6}.
\end{proof}
\leref{lswtranslation} shows that 
$\N_\W/\W_1 = \iota_*(\N_\S/\S_1) \oplus \kk \Id$.
Hence, every $\N_\S/\S_1$-module can be extended to an
$\N_\W/\W_1$-module by imposing the element $\Id$ to act as
multiplication by a scalar $c \in \kk$. These are the only
possible extensions if the action of $\N_\S/\S_1$ is irreducible.
\begin{theorem}\lbb{prrestriction}
Every tensor module for\/ $\Sd$ can be obtained as the restriction
of a tensor module for\/ $\Wd$. 
More precisely, for every\/ $c \in \kk$,
$\VS(R) = \VS(\Pi, U)$ is isomorphic to the restriction of\/ 
$\V(\Pi \tt \kk_{c(\chi - \tr\ad)/N}, U, c)$.
\end{theorem}
\begin{proof}
Note that, as an $H$-module, 
$V = \V(\Pi \tt \kk_{c(\chi - \tr\ad)/N}, U, c)$ 
can be identified with $H\tt R$.
Moreover, since $R\subset\sing V$, we have $\W_1\cdot R = \{0\}$.
Then $R$ becomes an $\N_\S/\S_1$-module
via the embedding $\iota_*$ from \leref{lswtranslation}.

We identify each of the Lie algebras $\wti\dd$ and $\what\dd$ with $\dd$.
It follows from \eqref{io*d} that if $\wti\dd$ acts on
$\Pi \tt \kk_{c(\chi -\tr\ad) /N}$ and $\Id$ acts as $c$,
then $\what \dd$ acts on $\Pi$. 
Similarly, $\sld$ acts on $U$, so the action of 
$\dd\oplus \sld \simeq \N_\S/\S_1$ on $R$ 
is isomorphic to $\Pi \bt U$. 

Then, by the definition of $\VS(R)$, there is a natural surjective
homomorphism of $\Sd$-modules
\begin{equation*}
\pi\colon \VS(R) = \VS(\Pi, U) \to HR=V \,.
\end{equation*}
The homomorphism $\pi$ takes an element
$u \in R \subset \VS(R)$ to the element $u \in R \subset V$.
But $V$ is free as an $H$-module; hence, $\pi$ is injective
and $\VS(R) \simeq V$.
\end{proof}

We will denote the restriction of $\T(\Pi,U)$ to $\Sd$ by $\TS(\Pi,U)$, 
and similarly for $\T(\Pi,U,c)$ (cf.\ \deref{dtenmodw}).
Note that by \eqref{wsing5}, we have
\begin{equation}
\T(\Pi \tt \kk_{\chi + c(\chi - \tr\ad)/N}, U, c) 
= \V(\Pi \tt \kk_{(N+c)(\chi - \tr\ad)/N}, U, N+c) \,.
\end{equation}
Then \thref{prrestriction} implies
\begin{equation}\lbb{tspuc}
\VS(\Pi, U) \simeq \TS(\Pi \tt \kk_{\chi + c(\chi - \tr\ad)/N}, U, c) \,,
\qquad c\in\kk \,.
\end{equation}
Observe that $\chi = \tr \ad$ is the only case for which the
restriction of $\T(\Pi, U, c)$ to $\Sd$ is independent of~$c$.

\begin{example}\lbb{etpi0n}
Note that $\T(\Pi,\Om^n)=\T(\Pi,\Om^n,-n)$, because 
the element $\Id\in\gld$ acts on $\Om^n = \bigwedge^n\dd^*$ 
as $-n$. Then it follows from \eqref{tspuc} that
\begin{equation}\lbb{tspuc1}
\TS(\Pi, \Om^n) 
\simeq \VS(\Pi \tt \kk_{-\chi + n(\chi - \tr\ad)/N}, \Om^n) \,.
\end{equation}
In particular, we have
\begin{equation}\lbb{tspuc0}
\VS(\Pi, \kk) \simeq \TS(\Pi \tt \kk_{\chi}, \Om^0) 
\simeq \TS(\Pi \tt \kk_{\tr\ad}, \Om^N) \,.
\end{equation}
\end{example}

One can use \eqref{tspuc} for $c=0$ to write an explicit expression for
the action of $\Sd$ on its tensor module $\VS(\Pi, U)$.
First, we note that if we identify $\T(\Pi \tt \kk_{\chi}, U, 0)$
with $H\tt R$, then $\d\in\dd$ acts on $R$ as $\bar\d$
(see \eqref{bard}). Then we use \eqref{sij} and \eqref{wirten2} to compute
$s_{ij}*(1\tt u)$ for $u\in R$.
The full expression is too cumbersome to write here. Because of \eqref{sij},
it is a sum of three terms. The third term is just a direct application of
\eqref{wirten2} for $[\d_i,\d_j]$. The second term is obtained from
the first one by switching the roles of $i$ and $j$. 
Finally, using $H$-bilinearity, \eqref{wirten2}, and \eqref{sdu2},
we find that the first term is equal to:
\begin{equation}\lbb{saction2}
\begin{split}
(\bar\d_i \tt \d_j) &* (1 \tt u) 
= (\bar\d_i \tt 1) \tt_H (1 \tt (\bar\d_j+\ad\d_j)u)
\\
&+ (\bar\d_i \tt 1) \tt_H s(\d_j,u)
- \sum_{k=1}^N\, (\bar\d_i \tt \d_k) \tt_H (1 \tt (e_j^k+\de_j^k)u) \,.
\end{split}
\end{equation}

Recall that $\De(\d)=\d\tt1+1\tt\d$ for $\d\in\dd$. This implies
$\De(\bar\d)=\bar\d\tt1+1\tt\d$ and
\begin{equation}\lbb{bardh}
\begin{split}
(\bar\d \tt g) \tt_H (h \tt u)
= (1 \tt g) \tt_H (\bar\d h \tt u)
- (&1 \tt g\d) \tt_H (h \tt u) \,,
\\
& \d\in\dd \,, \; g,h\in H \,, \; u\in R \,.
\end{split}
\end{equation}
Applying \eqref{bardh}, we rewrite \eqref{saction2} as follows:
\begin{equation}\lbb{saction3}
\begin{split}
(\bar\d_i \tt \d_j) &* (1 \tt u) 
= (1 \tt 1) \tt_H (\bar\d_i \tt (\bar\d_j+\ad\d_j)u)
\\
&- (1 \tt \d_i) \tt_H (1 \tt (\bar\d_j+\ad\d_j)u)
+ (1 \tt 1) \tt_H \bar\d_i s(\d_j,u)
\\
&- (1 \tt \d_i) \tt_H s(\d_j,u)
- \sum_{k=1}^N\, (1 \tt \d_k) \tt_H (\bar\d_i \tt (e_j^k+\de_j^k)u)
\\
&+ \sum_{k=1}^N\, (1 \tt \d_k\d_i) \tt_H (1 \tt (e_j^k+\de_j^k)u) \,. 
\end{split}
\end{equation}
%
%
When the action of $\gld$ on $R$ is trivial,
things can be rearranged in a more elegant form as follows:
\begin{equation}\lbb{sactiontrivial}
\begin{split}
s_{ij}* (1 \tt u) 
&= (1 \tt 1) \tt_H ({\d_i} \tt {\d_j}u - {\d_j} \tt
{\d_i}u - 1 \tt [\d_i, \d_j]u)
\\
&+ (1 \tt \d_i) \tt_H ({\d_j} \tt u - 1 \tt {\d_j}u)
\\
&- (1 \tt \d_j) \tt_H ({\d_i} \tt u - 1 \tt {\d_i}u) \,.
\end{split}
\end{equation}

Even though the expression \eqref{saction3} is not very inspiring, it will
turn out to be useful. We state as a lemma the properties that we are
going to need later. Before that let us introduce the notation
(see \eqref{sdu2}):
\begin{equation}\lbb{aru}
a_{ij}(u) = \d_i s(\d_j,u) - \d_j s(\d_i,u) 
= \sum_{k=1}^N \left( \d_i \d_k \tt e_j^k u -
\d_j \d_k \tt e_i^k u \right) \,.
\end{equation}
Note that $a_{ii}(u)=0$.

\begin{lemma}\lbb{lsacton}
Consider the tensor\/ $\Sd$-module\/ $\VS(R)$. Then the action
of the elements\/ $s_{ij} \in\Sd$, defined in \eqref{sij}, 
on an element\/ $1\tt u \in \kk\tt R \subset\VS(R)$ has the form
\begin{equation}\lbb{saction4}
\begin{split}
s_{ij}* (1 \tt u) 
&= (1 \tt 1) \tt_H A_{ij} (u)
+ \sum_{k=1}^N\, (1 \tt \d_k) \tt_H A_{ij}^k (u)
\\
&+ \sum_{\substack{ k,l=1 \\ k\le l }}^N\, 
(1 \tt \d_k\d_l) \tt_H A_{ij}^{kl} (u) \,,
\end{split}
\end{equation}
where
\begin{align*}
A_{ij}^{kl} (u) &\in \kk \tt (\kk+\sld) u \,,
\\
A_{ij}^k (u) &\in \dd \tt (\kk+\sld) u + \kk \tt (\kk+\what\dd+\sld) u \,,
\intertext{and}
A_{ij} (u) &\in a_{ij}(u)
+ \d_i \tt \what\d_j u - \d_j \tt \what\d_i u 
\\
&\qquad + \dd \tt (\kk+\sld) u + \kk \tt (\kk+\what\dd+\sld) u \,.
\end{align*}
\end{lemma}

\subsection{Filtration of Tensor Modules}\lbb{ssfil}
In analogy with tensor modules for $\Wd$, one can define an
increasing filtration $\fil^p \VS(R)$ of $\VS(R)= H \tt R$ by
\begin{equation}\lbb{filvsr}
\fil^p \VS(R) = \fil^p H \tt R \,, \qquad p \geq -1 \,. 
\end{equation}
Note that $\fil^{-1} \VS(R) = \{0\}$ and
$\fil^0 \VS(R) = \kk\tt R$.
The associated graded space of $\VS(R)$ is
\begin{equation}\lbb{grvsr}
\gr \VS(R) = \bigoplus_{p \geq -1} \gr^p \VS(R) \,, \qquad
\gr^p \VS(R) = \fil^p \VS(R)/\fil^{p-1} \VS(R) \,.
\end{equation}
The proof of the following lemma is completely similar to
that of \leref{lwfilten1}, so we omit it.

\begin{lemma}\lbb{lftms1}
The action of\/ $\wti\S$ on\/ $\VS(R)$ satisfies{\rm:}

{\rm(i)}
$\dd \cdot \fil^p \VS(R) \subset \fil^{p+1}\VS(R),$

{\rm(ii)}
$\N_\S \cdot \fil^p \VS(R) \subset \fil^p \VS(R),$

{\rm(iii)}
$\S_1 \cdot \fil^p \VS(R) \subset \fil^{p-1} \VS(R).$
\end{lemma}

\leref{lftms1} implies that each $\gr^p \VS(R)$ is a module
over the Lie algebra $\N_\S / \S_1 \simeq \dd\oplus\sld$.
This module is described in the next lemma.

\begin{lemma}\lbb{lftms2}
We have
\begin{equation*}
\gr^p \VS(\Pi, U) \simeq 
(\Pi\tt \kk_{p(\tr\ad -\chi)/N}) \bt (\symp^p \dd \tt U)
\end{equation*}
as\/ $(\dd\oplus\sld)$-modules.
\end{lemma}
\begin{proof}
Let us extend the $\sld$-action on $U$ to an action of $\gld$
by letting $\Id$ act as $0$. Then, by \thref{prrestriction},
$\VS(\Pi,U)$ is the restriction to $\Sd$ of the tensor $\Wd$-module
$\V(\Pi,U,0)$. Moreover, the filtration \eqref{filvsr} coincides
with the one defined in \seref{swfilten}. The structure of a
$(\dd\oplus\gld)$-module on $\gr^p \V(\Pi,U,0)$ is described in
\coref{cwfilten4}. Note that this describes the action of $\wti\dd$.
Using \eqref{io*d}, we find that $\what\dd$ acts as 
$\Pi\tt \kk_{p(\tr\ad -\chi)/N}$, 
because $\Id$ acts as $p$ on $\symp^p \dd \tt U$.
\end{proof}

The grading of $\wti \S$ can be used to endow $\VS(R)$ with
a graded module structure as follows. Recall that
$\VS(R)=\Ind_{\N_\S}^{\wti \S} R$ and
$\wti\S=\s_{-1} \oplus \N_\S$ as a vector space.
Therefore, as a vector space, $\VS(R) = \ue(\s_{-1}) \tt R$.
However, the Lie algebra $\s_{-1}$ is commutative, because
the degree $-1$ part in $S_N$ is commutative
and because the isomorphism $\S \simeq S_N$
is compatible with the grading (see \coref{csw7}).
Then $\ue(\s_{-1})$ is the symmetric algebra generated by $\s_{-1}$,
and we grade $\VS(R)$ by letting $\s_{-1}$ have degree $-1$ and $R$ 
have degree~$0$.

By definition, the above grading of $\VS(R)$ is compatible with 
the grading of $\wti\S$. It is also compatible with 
the filtration \eqref{filvsr}. 
%

\subsection{Submodules of Tensor Modules}\lbb{sstm}
In what follows, $\VS(R)$ will be a tensor module for $\Sd$. 
We will assume that 
$R = \Pi \bt U$, where $\Pi$ (respectively $U$) 
is an irreducible finite-dimensional representation of $\dd$ 
(respectively $\sld$).


Recall that every element $v \in \VS(R)$ can be expressed uniquely 
as a finite sum 
\begin{equation}\lbb{vvi2}
v = \sum_{I\in\ZZ_+^N} \d^{(I)} \tt v_I
\,, \qquad v_I \in R \,.
\end{equation}
As in \seref{swirten}, nonzero
elements $v_I$ are called {\em coefficients\/} of $v$, and we denote by
$\coef M$ the subspace of $R$ linearly spanned by coefficients
of elements $v\in M$.

\begin{lemma}\lbb{lm0triv}
For any nontrivial proper\/ $\Sd$-submodule\/ $M\subset\VS(R)$, we
have\/ $M \cap \fil^0 \VS(R) = \{0\}$.
\end{lemma}
\begin{proof}
The action of $\dd \oplus \sld$ preserves both $\kk \tt R$
(by the definition of $\VS(R)$) and $M$ (because it is an $\Sd$-submodule).
Thus, their intersection $M \cap \fil^0 \VS(R)$ is a 
$(\dd \oplus \sld)$-submodule of $R$. Irreducibility of $R$
implies that this intersection is trivial.
\end{proof}
\begin{lemma}\lbb{lscoef}
For any nontrivial proper\/ $\Sd$-submodule\/ $M\subset\VS(R)$, we have\/
$\coef M = R$.
\end{lemma}
\begin{proof}
Take an element $v\in M$ and write it in the form \eqref{vvi2}.
For $i\ne j$, we compute
$s_{ij}*v$ using $H$-bilinearity and \leref{lsacton}.
Denote by $m$ the coefficient of $1 \tt \d^{(I)}$ in the expression for
$s_{ij}*v$; then, by \reref{rlmod}(ii), $m\in M$.
By \leref{lsacton}, we have
\begin{equation*}
m = a_{ij}(v_I) \mod \fil^1 H \tt R \,.
\end{equation*}
Note that $e_i^k v_I$ for $k\neq i$, $e_j^k v_I$ for $k \neq j$, and 
$(e_j^j-e_i^i)v_I$ are coefficients of $a_{ij}(v_I)$ (see \eqref{aru}).
Hence, they are coefficients of $m$,
and we conclude that $(\sld)v_I \subset \coef M$.

In order to show that $\what\dd \cdot v_I \subset \coef M$, we look at
the degree one part of the above element $m$. Again by \leref{lsacton},
it is equal to 
\begin{equation*}
\d_i \tt \what\d_j v_I - \d_j \tt \what\d_i v_I 
+ \dd \tt (\kk + \sld)(\coef M) \,.
\end{equation*}
Hence, the action of $\what\dd$ preserves $\coef M$.
Then $\coef M$ is a $(\dd \oplus \sld)$-submodule of $R$,
which is irreducible. This shows that $\coef M=R$, as
it cannot be $\{0\}$.
\end{proof}
\begin{lemma}\lbb{laij}
Let\/ $M$ be a nontrivial proper\/ $\Sd$-submodule of\/ $\V(R)$. Then
for every\/ $i\ne j$ and\/ $u \in R$, there exists an element\/
$m \in M \cap \fil^2 \VS(R)$ such that
\begin{equation*}
m = a_{ij}(u) \mod \fil^1 \VS(R) \,,
\end{equation*}
where\/ $a_{ij}(u)$ is given by \eqref{aru}.
\end{lemma}
\begin{proof}
As $\coef M = R$, it is enough to prove the statement when $u$
is a coefficient $v_I$ of some element $v\in M$. Then $m$ is the element
considered in the proof of \leref{lscoef}.
\end{proof}

Our next result describes which tensor $\Sd$-modules are irreducible.

\begin{theorem}\lbb{trirr2}
Let\/ $\Pi$ $($respectively\/ $U)$
be an irreducible finite-dimensional\/ module over\/ $\dd$
$($respectively\/ $\sld)$.
Then the\/ $\Sd$-module\/ $\VS(\Pi, U)$ is irreducible if and only if\/
$U$ is not isomorphic to\/
$\Om^n = \bigwedge^n \dd^*$ for any\/ $n\ge0$. 
\end{theorem}
\begin{proof}
Recall from \thref{twirten} that the tensor $\Wd$-module $\T(\Pi,\Om^n)$
is not irreducible for $n\ge1$. Then its restriction $\TS(\Pi,\Om^n)$
is not irreducible either. 
It follows from \eqref{tspuc0} and \prref{pdomd2}
that $\VS(\Pi,\Om^0)$ is not irreducible as well.

Conversely, assume that $\VS(\Pi, U)$ is not irreducible,
and let $M$ be a nontrivial proper submodule.
Apply $s_{ij}$ on an element $m\in M$ satisfying
the conditions of \leref{laij} (for the same $i,j$).
Then for $i<j$
the coefficient multiplying $1 \tt \d_i^2 \d_j^2$ in the
expression for $s_{ij} * m$ is equal to (see \leref{lsacton}):
\begin{equation*}
1 \tt \bigl( (e_i^i-e_j^j)^2 - (e_i^j e_j^i+ e_j^i e_i^j) \bigr)u \,.
\end{equation*}
By \reref{rlmod}(ii), this is an element of $M$. \leref{lm0triv}
implies that this element vanishes for all $u\in R=\Pi\bt U$.

Note that $h=e_i^i-e_j^j$, $e=e_i^j$, and $f=e_j^i$ are standard
generators of a subalgebra of $\sld$ isomorphic to $\sl_2$. We know that
$h^2 - (ef +fe)$ acts trivially on $R$.
The element $h^2$ is a linear combination of $h^2 -(ef +fe)$ and
of the Casimir element; hence, it acts on any irreducible $\sl_2$-submodule
$W\subset U$ as a scalar. This means
that $\dim W = 0$ or $1$; hence, for $i<j$ the eigenvalues of $e_i^i-e_j^j$ 
on weight vectors can only be $0$ or $1$. 

Recall that every irreducible $\sld$-module $U$ has a 
highest weight vector, and $U$ is uniquely determined by its highest weight
(see, e.g., \cite[Chapter~VII]{Se}). 
Let us denote by $\la_{ij}$ the 
eigenvalue of $e_i^i-e_j^j$ on the highest weight vector of $U$.
Then $\la_{ij}+\la_{jk}=\la_{ik}$ but all $\la_{ij} = 0$ or $1$.
This implies that $\la_{i,i+1} = \de_{i,n}$ for some $n$; 
in other words, $U$ is the $n$-th fundamental representation, 
which is isomorphic to~$\Om^{N-n}$.
\end{proof}
\begin{corollary}\lbb{cnotomfr}
Let\/ $V$ be a finite\/ $\Sd$-module, and let $R \subset \sing V$ be an
irreducible\/ $(\dd \oplus \sld)$-submodule. 
Assume that\/ $R \simeq \Pi\bt U$ with\/
$U \not \simeq \Om^n$ for any\/ $n$.
Then as an\/ $H$-module, $HR \simeq H \tt R$.
\end{corollary}
\begin{proof}
By the definition of $\VS(R)$, there is a natural surjective
homomorphism of $\Sd$-modules $\VS(R) \to HR \subset V$. However,
by \thref{trirr2}, the tensor module $\VS(R)$ is irreducible.
Therefore, $HR \simeq \VS(R) = H \tt R$ is free as an $H$-module. 
\end{proof}

\subsection{Computation of Singular Vectors}\lbb{sssing}
In this section, we will compute singular vectors for all tensor
$\Sd$-modules of the form $\VS(\Pi,\Om^n)$, where
$\Pi$ is an irreducible finite-dimensional 
representation of $\dd$ and $\Om^n = \bigwedge^n \dd^*$.
This will be used in \seref{sifsm} for the
classification of all irreducible quotients of tensor modules.

\begin{proposition}\lbb{pss1}
For\/ $V = \VS(\Pi, \Om^n)$, we have\/ 
$\fil^0 V \subset \sing V \subset \fil^2 V$. 
Furthermore, if\/ $V = \VS(\Pi, \kk)$, then\/
$\sing V \subset \fil^1 V$.
\end{proposition}
\begin{proof}
Let us consider first the case $n=0$, i.e., $V = \VS(\Pi, \kk)$.
Let $v\in V$ be a singular vector written in the form \eqref{vvi2}.
Then using \eqref{sactiontrivial} and $H$-bilinearity, we get:
\begin{equation*}
\begin{split}
s_{ij}*v &=  \sum_I (1 \tt \d^{(I)}) \tt_H ( \d_i \tt 
\d_j v_I -  \d_j \tt  \d_i v_I - 1 \tt [\d_i,\d_j] v_I)
\\
&+ \sum_I (1 \tt \d^{(I)} \d_i) \tt_H ( \d_j \tt v_I - 1 \tt \d_j v_I)
\\
&- \sum_I (1 \tt \d^{(I)} \d_j) \tt_H ( \d_i \tt v_I - 1 \tt \d_i v_I) \,.
\end{split}
\end{equation*}
Assume that $v\not\in\fil^1 V$, and choose $I$ so that
$|I|$ is maximal among those for which $v_I \neq 0$. Then,
by \reref{rsings}, the
element multiplying $1 \tt \d^{(I)} \d_i$ in the above expression
must vanish. Hence, $v_I=0$, which is a contradiction.

Now let us assume that $n \neq 0, N$, i.e., $\Om^n \not\simeq\kk$.
We proceed as above: consider a singular vector 
$v = \sum \d^{(I)} \tt v_I$
and use \eqref{saction4} to compute $s_{ij} * v$. Then the 
coefficient of $1 \tt \d^{(I)}$ is equal to $a_{ij}(v_I)$ mod $\fil^1 V$. 
If $|I|>2$ and $v_I \neq 0$, this contradicts \reref{rsings}.
\end{proof}

Recall that the tensor $\Sd$-module $V=\VS(R)$ has a filtration 
$\{ \fil^p V \}$, given by \eqref{filvsr}.
If $v\in V$ is a nonzero singular vector, we can find a unique $p\ge0$
such that $v \in \fil^p V \setminus \fil^{p-1} V$.
Note that, by \leref{lftms1}, both $\fil^p V$ and $\gr^p V$ are 
$(\dd\oplus\sld)$-modules, and the natural projection 
$\pi^p \colon \fil^p V \to \gr^p V$ is a homomorphism.
Therefore, the restriction
\begin{equation}\lbb{pips}
\pi^p \colon \sing V \cap \fil^p V \to \gr^p V
\end{equation}
is a homomorphism of $(\dd\oplus\sld)$-modules.
Since $\pi^p(v) \ne0$, we obtain that 
the $\sld$-modules $\sing V$ and $\gr^p V$ contain an isomorphic 
irreducible summand.
We will utilize these remarks, together with the next lemma,
to study singular vectors.

\begin{lemma}\lbb{lss3}
Let\/ $V=\VS(\Pi,U)$, and let\/ $U'$ be an irreducible\/ 
$\sld$-submodule of\/ $\gr^p V$.
Assume that\/ $U' \not\simeq \Om^m$ for any\/ $m$, 
and that\/ $\dim U' > \dim U$. 
Then the submodule\/ $\pi^p(\sing V \cap \fil^p V) \subset \gr^p V$ 
does not intersect\/ $U'$.
\end{lemma}
\begin{proof}
Let 
\begin{equation*}
R=\Pi\bt U \,, \quad
R' = (\Pi\tt \kk_{p(\tr\ad -\chi)/N}) \bt U' \,.
\end{equation*}
By \leref{lftms2}, $R'$ is an irreducible $(\dd\oplus\sld)$-submodule 
of $\gr^p V$. If we assume that $\pi^p(\sing V \cap \fil^p V)$ 
intersects $U'$, then $\sing V$ contains a $(\dd\oplus\sld)$-submodule
isomorphic to $R'$.
Now \coref{cnotomfr} implies that $HR' \subset V$ 
is free as an $H$-module. But 
\begin{equation*}
\dim R' = (\dim\Pi) (\dim U') > (\dim\Pi) (\dim U) = \dim R \,.
\end{equation*}
Therefore, the $H$-submodule $HR' \subset V$ has a larger rank than 
$V=H\tt R$, which is impossible.
\end{proof}

\leref{lss3} is a powerful tool for studying singular vectors, when combined 
with the explicit knowledge of the $\sld$-modules $\gr^p V$. 
It follows from \leref{lftms2}
that $\gr^p V$ is a completely reducible $\sld$-module, all of whose
irreducible $\sld$-components are contained in
$\symp^p \dd \tt U$. In addition, since $\sing V \subset \fil^2 V$ 
by \prref{pss1}, we can restrict our attention to the cases
$p=1$ or $2$.
Our next result shows a typical application of these ideas.

\begin{proposition}\lbb{pss4}
If\/ $V = \VS(\Pi, \Om^n)$, $n \neq 1$, then\/ 
$\sing V \subset \fil^1 V$.
\end{proposition}
\begin{proof}
Recall that $\sing V \subset \fil^2 V$ by \prref{pss1}.
We want to show that $\pi^2(\sing V) = \{0\}$ (see \eqref{pips}).
We know from \leref{lftms2} that
\begin{equation*}
\gr^2 V \simeq (\Pi \tt \kk_{2(\tr \ad - \chi)/N}) \bt (\symp^2
\dd \tt \Om^n) \,. 
\end{equation*}
Thus, any irreducible $\sld$-submodule $U' \subset \gr^2 V$
is contained in $\symp^2 \dd \tt \Om^n$.
One can check (see \leref{lsldm}(iii) below) that all such $U'$ satisfy
$\dim U' > \dim \Om^n$ and $U' \not\simeq \Om^m$ for any $m$.
Hence, we can apply \leref{lss3} to conclude that
$\pi^2(\sing V) \cap U' = \{0\}$.
Therefore, $\sing V \subset \fil^1 V$.
\end{proof}

Note that the above proof does not hold in the case $n=1$, as one of the
irreducible $\sld$-summands in $\symp^2 \dd \tt \Om^1$ is isomorphic to
$\dd \simeq \Om^{N-1}$. To get a complete description of all singular vectors,
we need a detailed study of the $\sld$-modules 
$\symp^2 \dd \tt \Om^1$ and $\dd \tt \Om^n$. 

\begin{lemma}\lbb{lsldm}
{\rm(i)}
For\/ $1 \le n \le N-1$,
we have a direct sum of\/ $\sld$-modules{\rm:}
$\dd \tt \Om^n = \Om^{n-1} \oplus U'$, where\/ $U'$ is irreducible,
$\dim U' > \dim \Om^n$ and\/ $U' \not\simeq \Om^m$ for any\/ $m$.

{\rm(ii)}
We have a direct sum of\/ $\sld$-modules{\rm:}
$\symp^2 \dd \tt \Om^1 = \dd \oplus U''$,
where\/ $U''$ is irreducible,
$\dim U'' > \dim \Om^1 = N$ and\/ $U'' \not\simeq \Om^m$ for any\/ $m$.

{\rm(iii)}
For\/ $2 \le n \le N-1$,
every irreducible\/ $\sld$-submodule\/ $U' \subset \symp^2 \dd \tt \Om^n$
satisfies\/ $\dim U' > \dim \Om^n$ and\/ $U' \not\simeq \Om^m$ for any\/ $m$.
\end{lemma}
\begin{proof}
We will use Table 5 from the Reference Chapter of \cite{OV}.
Following \cite{OV}, we will denote by $R(\La)$ the irreducible
representation of $\sld\simeq\sl_N$ with highest weight $\La$.
We will denote by $\pi_n$ the $n$-th fundamental weight of $\sl_N$,
and we will set $\pi_0=\pi_N=0$. Note that $R(\pi_1) = \dd$
is the vector representation of $\sld$, and $R(0) = \kk$ is the trivial one.
Then we have $\symp^2 \dd \simeq R(2\pi_1)$ and 
\begin{equation*}
\Om^{n} = {\textstyle\bigwedge}^n \dd^* 
\simeq ( {\textstyle\bigwedge}^n \dd )^*  
\simeq {\textstyle\bigwedge}^{N-n} \dd
\simeq R(\pi_{N-n}) \,.
\end{equation*}
Using \cite{OV}, we find:
\begin{align*}
R(\pi_1) \tt R(\pi_{p}) 
& \simeq R(\pi_1 + \pi_{p}) \oplus R(\pi_{p+1}) \,,
\\
R(2\pi_1) \tt R(\pi_{p}) 
& \simeq R(2\pi_1 + \pi_{p}) \oplus R(\pi_1 + \pi_{p+1}) \,,
\end{align*}
and
\begin{align*}
\dim R(\pi_p) &= \binom{N}{p} \,,
\\
\dim R(\pi_1 + \pi_p) &= \frac{p}{p+1} (N+1) \binom{N}{p} \,,
\\
\dim R(2\pi_1 + \pi_p) &= \frac{p}{p+2} \binom{N+2}{2} \binom{N}{p} \,.
\end{align*}
{}From here, it is easy to finish the proof.
\end{proof}

Let us introduce some notation. For a $\dd$-module $\Pi$, we set
\begin{equation}\lbb{pin}
\Pi_n = \Pi \tt \kk_{-\chi + n(\chi - \tr\ad)/N} \,, \qquad
\Pi' = \Pi \tt \kk_{\tr\ad-\chi} \,.
\end{equation}
Then we can restate \eqref{tspuc1} as
\begin{equation}\lbb{tspuc2}
\TS(\Pi, \Om^n) \simeq \VS(\Pi_n, \Om^n) \,,
\end{equation}
while by \eqref{tspuc0} we have an isomorphism of $\Sd$-modules
\begin{equation}\lbb{tspuc3}
\psi\colon \TS(\Pi', \Om^N) \isoto \TS(\Pi, \Om^0) \,.
\end{equation}
Also, recall the $\Pi$-twisted pseudo de Rham complex
of $\Wd$-modules \eqref{domd2}. When we restrict these modules
to $\Sd$, we obtain a complex of $\Sd$-modules
\begin{equation}\lbb{domd3}
0 \to \TS(\Pi,\Om^0) \xrightarrow{\di_\Pi} \TS(\Pi,\Om^1)
\xrightarrow{\di_\Pi} \cdots \xrightarrow{\di_\Pi} \TS(\Pi,\Om^N) \,.
\end{equation}
Note that the isomorphism $\psi$ is compatible with the filtrations
(i.e., it maps each $\fil^p$ to $\fil^p$), while $\di_\Pi$ has degree $1$
(i.e., it maps each $\fil^p$ to $\fil^{p+1}$).

\begin{theorem}\lbb{tsvs}
Let\/ $\Pi$ be an irreducible finite-dimensional\/ $\dd$-module.
Then we have the following equalities and isomorphisms of\/
$(\dd \oplus \sld)$-modules{\rm:}
\begin{align}
\tag{i}
\sing \TS(\Pi, \Om^n) 
&= \fil^0 \TS(\Pi, \Om^n) + \di_\Pi \fil^0 \TS(\Pi, \Om^{n-1})
\\
\notag
&\simeq (\Pi_n \bt \Om^n) \oplus (\Pi_{n-1} \bt \Om^{n-1}) \,,
\qquad 2\le n\le N \,,
\end{align}
\begin{align}
\tag{ii}
\sing \TS(\Pi, \Om^1) 
&= \fil^0 \T(\Pi, \Om^1) + \di_\Pi \fil^0 \TS(\Pi, \Om^0)
+ \di_\Pi \psi \di_{\Pi'} \fil^0 \TS(\Pi', \Om^{N-1})
\\
\notag
&\simeq (\Pi_1 \bt \Om^1) \oplus (\Pi_0 \bt \Om^0) 
\oplus (\Pi_{-1} \bt \Om^{N-1}) \,,
\end{align}
where we use the notation from\/ \eqref{pin}--\eqref{tspuc3}.
\end{theorem}
\begin{proof}
Let $V = \TS(\Pi, \Om^n)$. Then by \eqref{tspuc2},
$V \simeq \VS(\Pi_n, \Om^n)$, and by \leref{lftms2}, we have
an isomorphisms of $(\dd \oplus \sld)$-modules
\begin{equation*}
\gr^p V \simeq \Pi_{n-p} \bt (\symp^p \dd \tt \Om^n) \,.
\end{equation*}
In particular, $\fil^0 V = \gr^0 V \simeq \Pi_{n} \bt \Om^n$.
Note that the latter is an irreducible $(\dd \oplus \sld)$-module.
This implies the isomorphisms in (i) and (ii) above, because
$\di_\Pi$ and $\psi$ are homomorphisms and because
$(\Pi')_{N-1} \simeq \Pi_{-1}$ (see \eqref{pin}).

Recall that $\fil^0 V \subset \sing V \subset \fil^2 V$, and
$\sing V \subset \fil^1 V$ for $n \neq 1$
(see Propositions \ref{pss1} and \ref{pss4}).
Since $\di_\Pi$ and $\psi$ are homomorphisms of $\Sd$-modules,
they map singular vectors to singular vectors. Then it is clear
that the right-hand sides of (i) and (ii) are contained in $\sing V$.

Next, we describe the image of 
$\sing V \cap \fil^1 V$ in $\gr^1 V$ 
under the natural projection \eqref{pips}.
On one hand, we have
\begin{equation*}
\pi^1 (\sing V \cap \fil^1 V) 
\supset \pi^1 ( \di_\Pi \fil^0 \TS(\Pi, \Om^{n-1}) )
\simeq \Pi_{n-1} \bt \Om^{n-1} \,.
\end{equation*}
On the other hand, every irreducible $\sld$-submodule of $\gr^1 V$
is contained in $\dd \tt \Om^n$. 
By \leref{lsldm}(i), we have a direct sum of $\sld$-modules:
$\dd \tt \Om^n = \Om^{n-1} \oplus U'$, where $U'$ is irreducible,
$\dim U' > \dim \Om^n$ and $U' \not\simeq \Om^m$ for any $m$.
Now, by \leref{lss3}, the image $\pi^1 (\sing V \cap \fil^1 V)$
does not intersect $\Pi_{n-1} \bt U'$. 
Therefore, the above inclusion is an equality. 
In particular, we get statement (i).

To finish the proof of (ii), we note that
\begin{equation*}
\pi^2 (\sing V)
\supset \pi^2 ( \di_\Pi \psi \di_{\Pi'} \fil^0 \TS(\Pi', \Om^{N-1}) )
\simeq \Pi_{-1} \bt \Om^{N-1} \,.
\end{equation*}
By the same argument as above, this is an equality, because 
of \leref{lsldm}(ii).
\end{proof}
\begin{remark}\lbb{rsvs0}
It follows from \thref{tsvs} and the isomorphism \eqref{tspuc3} that
\begin{align*}
\sing \TS(\Pi, \Om^0) 
&= \fil^0 \TS(\Pi, \Om^0) + \psi \di_{\Pi'} \fil^0 \TS(\Pi', \Om^{N-1})
\\
&\simeq (\Pi_0 \bt \Om^0) \oplus (\Pi_{-1} \bt \Om^{N-1}) \,.
\end{align*}
\end{remark}

\subsection{Irreducible Finite $\Sd$-Modules}\lbb{sifsm}
We can now complete the classification of irreducible finite
$\Sd$-modules. Our fist result describes all submodules
of the tensor $\Sd$-module $\TS(\Pi, \Om^n)$.

\begin{lemma}\lbb{lsirr2}
Let\/ $\Pi$ be an irreducible finite-dimensional\/ $\dd$-module,
let\/ $T^n=\TS(\Pi, \Om^n)$, and let\/ $M\subset T^n$
be a nontrivial proper\/ $\Sd$-submodule.
Then{\rm:}

{\rm(i)}
$\sing M = \di_\Pi \fil^0 T^{n-1}$, if\/ $2\le n\le N;$

{\rm(ii)}
$M \supset \di_\Pi T^{N-1}$, if\/ $n=N;$

{\rm(iii)}
$M=\di_\Pi T^{n-1}$, if\/ $2\le n\le N-1;$

{\rm(iv)}
$\di_\Pi T^{n-1}$ is irreducible for\/ $2\le n\le N$.
\end{lemma}
\begin{proof}
Let $2\le n\le N$, and let 
$M \subset T^n$ be a nontrivial proper $\Sd$-submodule. 
Then $\sing M \subset \sing T^n$ is a $(\dd \oplus \sld)$-submodule, 
and $M \cap \fil^0 T^n = \{0\}$ by \leref{lm0triv}. 
Now \thref{tsvs}(i) and an argument similar to the one used in the proof 
of \coref{cwcsv3} imply part (i). Then 
\begin{equation*}
M \supset H(\sing M) = \di_\Pi (H (\fil^0 T^{n-1})) = \di_\Pi T^{n-1} \,.
\end{equation*}
The rest of the proof of (iii) is the same as that of \leref{lwirr2},
while (iv) follows from (ii) and (iii).
\end{proof}
\begin{remark}\lbb{rsirr1}
Recall that the $\Wd$-module\/ $\T(\Pi,\Om^1)$ has a unique 
nontrivial proper $\Wd$-submodule, namely $\di_\Pi \T(\Pi,\Om^0)$
(see \leref{lwirr2}).
However, the restriction $\TS(\Pi,\Om^1)$ of $\T(\Pi,\Om^1)$ to $\Sd$
has two nontrivial proper $\Sd$-submodules: 
\begin{equation*}
\di_\Pi \psi \di_{\Pi'} \TS(\Pi', \Om^{N-1})
\subset \di_\Pi \TS(\Pi, \Om^0)
\end{equation*}
(cf.\ \thref{tsvs}(ii)).
Because of \eqref{tspuc3} and the exactness of \eqref{domd3},
these two $\Sd$-modules are isomorphic to the following ones:
\begin{equation*}
\di_{\Pi'} \TS(\Pi', \Om^{N-1}) \subset \TS(\Pi', \Om^N) \,.
\end{equation*}
\end{remark}

Now we can state the main result of this section.

\begin{theorem}\lbb{tsmod}
Any irreducible finite\/ $\Sd$-module is isomorphic to one of the
following\/{\rm:}

{\rm(i)} Tensor modules\/ $\TS(\Pi, U, 0)$, where\/ $\Pi$ is an
irreducible finite-dimensional\/ $\dd$-module, and\/ $U$ is an
irreducible finite-dimensional\/ $\sld$-module not isomorphic to\/
$\Om^n = \bigwedge^n \dd^*$ for any\/ $0\le n \le\dim\dd;$

{\rm(ii)} Images\/ $\di_\Pi \TS(\Pi, \Om^n)$, where\/ $\Pi$ is an
irreducible finite-dimensional\/ $\dd$-module, and\/ 
$1\le n \le\dim\dd-1$ $($see~\eqref{domd3}$)$.
\end{theorem}
\begin{proof}
The proof is similar to that of \thref{twmod}.
Let $V$ be an irreducible finite $\Sd$-module.
Then, by \thref{prstensoruniversal} and \eqref{tspuc},
$V \simeq T/M$, where $T=\TS(\Pi,U)$ is a tensor module and
$M\subset T$ is an $\Sd$-submodule.

If $U \not\simeq \Om^n$
as an $\sld$-module for any $n\ge0$,
then $T$ is irreducible by \thref{trirr2} and \eqref{tspuc}.
In this case, $V\simeq\TS(\Pi,U)$.

Assume that $U \simeq \Om^n$ for some $n\ge0$;
then $T \simeq \TS(\Pi,\Om^n) = T^n$ is not irreducible.
Because of \eqref{tspuc0}, we can assume without loss of generality
that $1\le n\le N=\dim\dd$.
Now if $2\le n\le N-1$, \leref{lsirr2}(iii) implies
that $M=\di_\Pi T^{n-1}$. By the exactness of \eqref{domd3},
we get $V \simeq T^n / \di_\Pi T^{n-1} \simeq \di_\Pi T^n$.

Next, consider the case when $V$ is a quotient of $T^N$. 
Then, by \leref{lsirr2}(ii), we have $M \supset \di_\Pi T^{N-1}$.
Now \prref{pdomd2} implies that $V$ is finite dimensional; 
hence, $\Sd$ acts trivially on it by \exref{elv0}, and $V$ is not
irreducible.

Finally, it remains to consider the case when $V$ is a quotient of $T^1$.
Note that $\di_\Pi M$ is a proper $\Sd$-submodule of $T^2$; hence,
by \leref{lsirr2}(iii), it must be either trivial or equal to $\di_\Pi T^1$.
First, if $\di_\Pi M = \{0\}$, then $M \subset \di_\Pi T^0$ and we have a
surjective homomorphism $T^1/M \to T^1 / \di_\Pi T^0$.
But $T^1/M \simeq V$ is irreducible; therefore,
$V \simeq T^1 / \di_\Pi T^0 \simeq \di_\Pi T^1$.
Second, if $\di_\Pi M = \di_\Pi T^1$, then $M + \di_\Pi T^0 = T^1$
and we have isomorphisms 
$V \simeq T^1/M \simeq (\di_\Pi T^0) / (\di_\Pi T^0 \cap M)$.
Since the map $\di_\Pi \colon T^0 \to T^1$ is injective, we get that
$V \simeq T^0/K$ for some $\Sd$-submodule $K$ of $T^0$.
This case was already considered above, because of \eqref{tspuc0}.
\end{proof}

Finally, for each irreducible finite $\Sd$-module $V$,
we will describe the space $\sing V$ of singular vectors of $V$.

\begin{lemma}\lbb{lrirr1}
Let\/ $R$ be an irreducible finite-dimensional\/ $(\dd\oplus\sld)$-module.
Then\/ $V=\VS(R)$ is an irreducible $\Sd$-module if and only if\/
$\sing V = \fil^0 V$.
\end{lemma}
\begin{proof}
It is clear from \leref{lm0triv} that $V$ is irreducible when
$\sing V = \fil^0 V$. Conversely, assume that $V$ is irreducible.
Consider the grading of $V = \ue(\s_{-1}) \tt R$ 
constructed at the end of \seref{ssfil}.
All homogeneous components of a singular vector are still
singular, so we have to show that the only homogeneous
singular vectors in $V$ are of degree zero. 
If $v \in \sing V$ is a singular vector of negative degree,
then the $\wti \S$-submodule generated by $v$ is 
contained in the negatively graded part of $V$, which
contradicts the irreducibility of $V$. Therefore,
$\sing V = R = \fil^0 V$.
\end{proof}
\begin{theorem}\lbb{tsnoiso}
The irreducible finite\/ $\Sd$-modules listed in \thref{tsmod}
satisfy {\rm(}see \eqref{pin}{\rm):}

{\rm(i)} $\sing\TS(\Pi, U, 0) \simeq \Pi_0 \bt U$
as\/ $(\dd\oplus\sld)$-modules$;$

{\rm(ii)} $\sing(\di_\Pi \TS(\Pi, \Om^n)) \simeq
\Pi_n \bt \Om^n$ as\/
$(\dd\oplus\sld)$-modules.
\\
In particular, no two of them are isomorphic to each other.
\end{theorem}
\begin{proof}
The proof is similar to that of \thref{twnoiso}, and it uses
\eqref{tspuc2}, \thref{tsvs}(i), and 
Lemmas \ref{lsirr2}(i) and \ref{lrirr1}.
\end{proof}

\section*{Acknowledgments}
We acknowledge the hospitality of MSRI (Berkeley) and ESI (Vienna), 
where parts of this work were done. 
We are grateful to the referee for several useful comments which improved 
the exposition. 

\bibliographystyle{amsalpha}


\end{document}